\documentclass[11pt]{extarticle}
\usepackage[utf8]{inputenc}
\usepackage[english]{babel}
\usepackage[a4paper, total={7in, 10in}]{geometry}
\usepackage{parskip}
\usepackage[stretch=10,shrink=10]{microtype}
\usepackage{tgtermes}
\usepackage{mathtools}
\usepackage[round]{natbib}
\usepackage{amsmath,amssymb}
\usepackage{amsthm}
\usepackage{tabularx}
\usepackage{booktabs}
\usepackage{multicol}
\usepackage{multirow}
\usepackage{geometry}
\usepackage{array}
\usepackage{graphicx}
\usepackage{wrapfig}
\usepackage{appendix}
\usepackage[normalem]{ulem}
\usepackage[section]{placeins}
\usepackage{adjustbox}
\usepackage{caption}
\usepackage{subcaption}
\usepackage{float}
\usepackage{authblk}
\usepackage[table]{xcolor}
\usepackage{comment}
\usepackage[colorlinks = true,
            linkcolor = red,
            urlcolor  = blue,
            citecolor = blue,
            anchorcolor = red]{hyperref}

\newcolumntype{P}[1]{>{\centering\arraybackslash}p{#1}}

\numberwithin{equation}{section}

\definecolor{darkred}{rgb}{0.70, 0.10, 0.10}
\definecolor{darkblue}{rgb}{0.0, 0.0, 0.75}

\definecolor{createcolor_blue}{rgb}{0.0, 0.5, 0.8}

\providecommand{\keywords}[1]
{
  \small	
  \textbf{\textit{Keywords---}} #1
}

\title{\textbf{A Decomposition Approach to Last Mile Delivery Using Public Transportation Systems}}

\author[]{Minakshi Punam Mandal}
\author[]{Claudia Archetti}
\affil[]{Department of Information Systems, Decision Sciences and Statistics \\ ESSEC Business School, 95000 Cergy, France}

\date{}

\begin{document}

\maketitle

\begin{abstract}
    This study explores the potential of using public transportation systems for freight delivery, where we intend to utilize the spare capacities of public vehicles like buses, trams, metros, and trains, particularly during off-peak hours, to transport packages within the city instead of using dedicated delivery vehicles. The study contributes {to the growing} literature on innovative strategies for performing sustainable last mile deliveries. We study an operational level problem called the \textit{Three-Tier Delivery Problem on Public Transportation}, where packages are first transported from the Consolidation and Distribution Center (CDC) to nearby public vehicle stations by delivery trucks. From there, public vehicles transport them into the city area. The last leg of the delivery is performed to deliver the packages to their respective customers using green vehicles or eco-friendly systems. We propose mixed-integer linear programming formulations to study the transport of packages from the CDC to the customers, use decomposition approaches to solve them, and provide numerical experiments to demonstrate the efficiency and effectiveness of the system. Our results show that this system has the potential to drastically reduce the length of trips performed by dedicated delivery vehicles, thereby reducing the negative social and environmental impacts of existing last mile delivery systems.
\end{abstract}


% \newpage

\section{Introduction} \label{intro}

The boom in e-commerce over the last decade poses several challenges for last mile deliveries, particularly the need to make them more efficient, less expensive, and more sustainable. It affects not only the shipping companies but urban life as well. The negative impacts caused by delivery vehicles can be three-fold -- economic, social, and environmental \citep{rai2017improving, viu2020impact}. Economic impacts include unreliable delivery times, delivery route inefficiencies, and use of resources, among others. Social impacts comprise traffic and congestion, nuisance caused by vehicles on the road, health risks, and noise pollution. Finally, environmental impacts include emissions of harmful pollutants like CO$_{2}$, Particulate Matter and other Greenhouse Gases, and use of non-renewable resources \citep{rai2017improving, viu2020impact}. A study by the World Economic Forum \citep{deloison2020future} states that there will be a 36\% rise in the number of delivery vehicles on urban roads by 2030 due to the consistent and exponential growth of demands for e-commerce deliveries, leading to an additional emission of 6 million tonnes of CO$_{2}$. It can potentially increase the average commute time of each passenger by 21\% (amounting to an additional 11 minutes) per day due to congestion. Thus, the need of the hour is solutions that are not only sustainable but also competent and economically viable.

To keep up with the growing demands and resulting issues caused by e-commerce delivery systems, companies are looking towards innovative approaches that reduce the costs of social and environmental externalities. Several companies have tested the use of drones for the final leg of the delivery, including Google's Project Wing, Amazon's Amazon PrimeAir, DHL's PaketKopter, and GeoPoste's GeoDrone \citep{ranieri2018review}. Several studies propose using electric and other green vehicles as a viable option for delivery as well \citep{mucowska2021trends, ranieri2018review}. The use of Urban Micro-Consolidation Centers is also being studied extensively for last mile deliveries. These are establishments where goods can be unloaded from delivery trucks, reconsolidated, and transferred to more sustainable means of transport \citep{arrieta2022location}. The utilization of pickup points and lockers are effective, popular, and well-established alternatives \citep{morganti2014final, iwan2016analysis, cleophas2019collaborative} because of convenience and limited restrictions on delivery times for customers. Crowdshipping has gained popularity over the last few years, where parcels are matched with couriers \citep{rai2017crowd, gatta2019public, mucowska2021trends, le2019supply}. The parcels could be dropped off in automated parcel lockers in and around public transit stations or delivered directly to the customers \citep{gatta2019public, gatta2019sustainable}. Examples include Amazon Flex, DHL's Parcel Metro, and Walmart utilizing Uber, Deliv, Lyft, and even their voluntary employees on their way back from work \citep{le2019supply}. Another alternative is the ride-sharing system, where passengers and goods are transported in the same ride, for example, a shared taxi service \citep{cleophas2019collaborative, li2014share}. Different shipping companies can also collaborate to reduce overall delivery costs \citep{cleophas2019collaborative}.

In this paper, we propose to use public transportation systems to transport packages in urban areas. Over recent years, an increasing number of studies have identified and emphasized the opportunities and benefits of utilizing transit networks \citep{trentini2010toward, trentini2012flow, galkina2019investigating, taniguchi2016new}. One of the most significant advantages of using public transit vehicles is that it does not add any new vehicles on the road solely for package delivery, which means no extra congestion and no extra emissions caused by traditional delivery trucks. This would also lead to more reliable deliveries as public transit vehicles run on fixed schedules. The reduced congestion due to the removal of delivery vehicles on roads would further ensure reliable running times of vehicles on the road. The transit companies would be compensated accordingly, and the revenue generated can be used to increase the services for passengers \citep{taniguchi2016new}. Thus it could be a win-win solution for all parties involved. The idea is to use the spare capacities of public vehicles like buses, trams, trains, subways, and other public vehicles, particularly during off-peak hours, or, in some cases, to have dedicated spaces for freight in the vehicles.

The effectiveness of such integrated delivery systems has already been demonstrated in several cities. \citet{marinov2013urban} discuss the case of the French supermarket chain Monoprix which uses urban trains and sustainable vehicles to transport non-perishable products in the city of Paris. This resulted in 10,000 fewer trucks circulating in Paris \citep{AQTr2019FreightonTransit}. Using commuter trains along with natural gas trucks also lowered the emission of CO$_{2}$ by 280 tonnes and Nitrogen Oxide by 19 tonnes. An article by \citet{saito2021nextstop} reports two companies, Nishi Tokyo Bus Co. and Yamato Transport Co., using buses to deliver packages on a route from the city of Akiruno to the village of Hinohara in western Tokyo. This service provided the bus company with a new revenue source and helped it to continue servicing the unprofitable route. For the delivery company, the venture turned out to be profitable as well, as it reduced round-trip traveling distances by 50 kilometers a day. The village has also welcomed the collaboration and acknowledged its benefits. The report says that the companies are further looking toward implementing similar services in the mountainous areas after the success of their initial venture. GB Railfreight, a freight delivery service based in the UK, has explored using old commuter trains to deliver parcels. They found in their trials that standard roll cages can be easily loaded and offloaded in commuter trains at most of the mainline stations \citep{railway2020GBRailfreight}. Amazon is also looking to use public buses for its deliveries and has received a patent that would transform buses into parcel carriers \citep{StarTribune2019amazon, RetailDetail2019amazon}. This service is aimed at customers who do not live near their pickup points or where there is a scarcity of carriers for package delivery. The report mentions that Amazon was looking to invest \$1.5 million in the public transportation system in Seattle. This would also benefit public transport companies, which are threatened by systems like Uber and Lyft. A mobile delivery module would be attached to the buses, and customers would have to choose nearby stops and time windows according to bus schedules and pick up their parcels from the bus stops. An article by Hermes \citep{Hermes2019LastMile} discusses several implementations of the freight-on-transit (FOT) system in practice. For example, the TramFret project in the French city of Saint-{\'E}tienne used decommissioned trams to deliver goods inside the city. Moscow utilizes its metro network to deliver parcels from one end of the city to another. Hermes and the Frankfurt transport authority have partnered together to transport two boxes filled with packages from a hub outside the city to Europaviertel, which is a housing and business district inside Frankfurt \citep{Hermes2019LastMile}. Then they use e-bikes to deliver the packages to their destinations. The researchers involved identify that such delivery projects have great potential to be impactful in the cities, particularly in pedestrian zones or areas around the transit stations. The public transportation network has also been used for waste disposal in several cities; for example, the tram network in Zurich \citep{Hermes2019LastMile}, garbage trains in New York and Toronto \citep{AQTr2019FreightonTransit}. \citet{galkina2019investigating} analyze the effectiveness of FOT by conducting a study in the city of Bratislava and find that it has the potential to reduce overall transportation costs by 8-12 times. All these examples demonstrate the viability and the advantages of shared transportation systems. 

Inspired by the real-life examples of FOT, particularly the type of delivery systems proposed by Amazon and the one implemented by Hermes, we wish to explore the operational implementation of the system and optimize the delivery of packages on the network in a cost and time-effective manner. We consider a delivery company that wants to utilize the public transit network for its last mile deliveries. The problem can be divided into three tiers. In the first tier (also referred to as T1), delivery trucks belonging to the company carry the packages from the Consolidation and Distribution Center (CDC) to nearby stops, called \textit{drop-in} stops, of the public transportation network. The second tier (T2) of the delivery is the one that occurs on board public vehicles, which have pre-determined schedules, itineraries, and stops. The vehicles pick up the packages from the drop-in stops and transport them to some other stops on their route, which we call \textit{drop-out} stops. These stops are within the city and typically close to customer locations. Finally, the city freighters pick up the packages from the drop-out stops and deliver them to the customers using sustainable and green modes of transport, like electric vehicles, drones, bikes, or even freighters simply walking for the delivery. These freighters could either belong to the company's fleet or be comprised of autonomous crowdsourced drivers. This constitutes the third and final tier (T3) of the system. We wish to find a delivery plan such that the delivery costs of the first and third tiers, and thereby distance traversed using dedicated delivery means, are minimized. We name this problem the \textit{Three-Tier Delivery Problem on Public Transportation} (3T-DPPT hereafter).

The contributions of the paper are the following. We introduce the 3T-DPPT and provide a mixed-integer linear programming formulation for it. We propose a decomposition methodology to solve it. 
Based on the natural three tiers of the problem, we provide three different approaches for decomposing and integrating the model. For the second tier of delivery on public vehicles, we formulate and analyze three objective functions that support the primary intent of the system, which is to minimize delivery distances using dedicated vehicles. We generate instances that mimic real-life public transportation networks and package demands, and we implement our models on them. Finally, we inspect our solutions and make recommendations for implementing such a delivery system.

The remainder of the paper is organized as follows. In Section \ref{lit_review}, we provide a review of the literature with regard to our problem. Section \ref{problem_formulation} describes the setting of the problem and introduces a formulation for it. We dedicate Section \ref{decomposition_approaches} to describing different approaches for solving the problem. Section \ref{numerical_experiments} provides numerical studies, and we conclude in Section \ref{Conclusion} with some suggestions for future research avenues.

\section{Literature Review} \label{lit_review}

In this section, we review the literature that studies deliveries using public transportation networks in some capacity, either together with passengers in the same vehicle or in isolation: only making use of the infrastructure. The works can broadly be classified into two categories from the perspective of our study. The first category is where the conceptual models of such a system are introduced, and the effectiveness of the integration is considered either theoretically or through small-scale implementations in practice. The second is the group of papers that study similar systems and propose mathematical models to optimize freight transport on public transit networks. We briefly discuss the first category and focus our attention primarily on reviewing the works in the second category.

\citet{trentini2010toward} and \citet{trentini2012flow} conducted one of the earliest studies on the topic. Apart from providing a survey of existing examples of cities that have adapted strategies for shared passenger and goods flows, they provide a conceptual model for the same, which they then show how to implement using a case study in the French city of La Rochelle. \citet{gatta2019public, gatta2019sustainable} combine crowdshipping with public transport, where the packages are dropped off or picked up by people traveling by these transit systems, specifically in or around public transit stations. On a survey in the city of Rome, \citet{gatta2019public} estimate a reduction of 239kg of particulate matter each year. \citet{villa2021metro} study the potential of using metro networks and their existing capacities and lockers in metro stations for parcel delivery. They investigate two kinds of scenarios: making use of the spare capacities of vehicles or utilizing dedicated runs of freight trains on the existing lines. They study the system from the perspective of costs and impacts and also provide a case study on the city of Madrid. They find that a shared system has 11.16\% to 14.72\% lower operating costs than current systems, and the average external delivery cost per parcel is 8.2 to 9.8 times lower. While several other studies analyze the system's viability in real life, due to the scope of our study, we concentrate on works that employ mathematical programming models to optimize integrated delivery systems.

\citet{crainic2009models} provided one of the first modeling frameworks for tactical and operational strategies of two-tier city logistics systems, with the possibility of using public vehicles like trams in the first tier and electric vehicles in the second tier. They propose generalized two-tier models where the movement of urban vehicles and freight is integrated, and consider routing and scheduling decisions. 

\citet{cheng2018planning, cheng2018packages} study the distribution of packages using Crowdsourced Public Transportation Systems. \citet{cheng2018packages} model the transport of packages using the idle capacities of public transportation systems, where the packages can be loaded at a starting node, unloaded and reloaded at intermediary nodes, and finally unloaded at their destination stop, using a multi-commodity flow model. They also propose a heuristic to solve the problem. In \citet{cheng2018planning}, the authors explore the problem further. They divide the study into two parts -- the first being a Passenger Transit Model that estimates the number of passengers at each station and, thereby, calculates the under-utilized capacity. In the second part, the decisions about the actual assignment and delivery of the packages are made. They propose two approaches for the second part -- the Minimum Limitation Delivery Method, which uses only the minimum under-utilized capacity, and the Adaptive Limitation Delivery Method (ALD), which utilizes the entire under-utilized capacity at each trip. They find ALD to perform better, with only slightly higher risks of affecting the quality of passenger experience. The main focus of their study is on the transfer using public vehicles, which constitutes a part of the problem studied in our paper.

\citet{ghilas2013integrating, ghilas2016adaptive, ghilas2016pickup, ghilas2016scenario} utilize scheduled public transportation lines for freight delivery. They study the Pickup and Delivery Problem with Time Windows and Scheduled Lines (PDPTW-SL), where freight requests are transported via public vehicle lines from one end to another. The authors provide an arc-based mixed-integer programming formulation for the problem and perform computational studies to demonstrate the benefits of using such a system \citep{ghilas2013integrating, ghilas2016pickup}. In \citet{ghilas2016adaptive}, the authors use an Adaptive Large Neighborhood Search (ALNS) heuristic algorithm to solve the problem on several synthetic instances, and on an instance generated based on the metro system in Amsterdam. The results on instances of sizes up to 100 indicate significant benefits in terms of cost savings, which range from 0 to 30\%, reduced driving times ranging from 0 to 31\%, and, proportionally, reduced CO$_{2}$ emissions. The authors also study a stochastic version of the same problem in \citet{ghilas2016scenario}, called the Pickup and Delivery Problems with Time Windows, Scheduled Lines and Stochastic Demands, where the demands of the requests are considered uncertain. They solve the problem using ALNS embedded into a sample average approximation method. Their computational studies show up to 16\% reduced costs compared to a traditional PDPTW on instances of sizes up to 40 requests (each request being a pickup and delivery destination pair). In this paper, we assume that the trucks pick up the packages from only one CDC instead of multiple pickup points and then decompose the problem to solve it. We do not associate time windows with the origin of the packages. We also assume that the packages can be picked up and dropped off at several points on each line, not just end-to-end transfer on the lines. 

\citet{mourad2021integrating} also study a version of the PDPTW-SL, where the scheduled lines consist of shuttles that can transport passengers along with robots carrying packages. In their setting, robots replace the pickup and delivery vehicles, and these robots perform the task of transporting packages between the stations and their origins and destinations. They study a stochastic version of the problem where the capacity of the shuttles is unknown and use an ALNS algorithm to solve it. Their method solves instances with up to 60 requests and finds solutions within 0.6\% of the optimal solutions. \citet{huang2020new, huang2020scheduling} investigate the deliveries of parcels that use drones interacting with existing public transit vehicles like trains, trams, etc. They propose algorithms and simulations to demonstrate the performance of their models. \citet{fatnassi2015planning} also explore the idea of transporting goods and passengers in a shared system, where they propose to use personal rapid transit and freight rapid transit (FRT) alternatively in a shared transportation network. They propose a dynamic or on-demand model that minimizes the waiting time of the passengers and goods, along with the movement of empty vehicles.

\citet{masson2017optimization} introduce the Mixed Urban Transportation Problem. They consider the transportation of packages from the CDC to the city on a bus route and subsequent deliveries from the bus stations by city freighters simultaneously. They categorized the problem into two specific classes of Vehicle Routing Problems (VRP)-- the Two-Echelon VRP (2E-VRP) and the Pickup and Delivery Problem with Transfers. To solve the problem, they also use an ALNS methodology. They evaluate their model in the city of La Rochelle. Our problem is an extension of the setting considered by \citet{masson2017optimization}, as we consider a more generalized delivery system with multiple public vehicle routes. We also have the option of delivering packages by trucks to the drop-in stops, where they can be collected, instead of a bus passing through the CDC and collecting them there. Their primary focus is on the last leg of the delivery, where most of the decisions about package movement, particularly on public transportation, are already made beforehand.

\citet{azcuy2021designing} also study the potential of using public transit for last mile deliveries. They have a two-tier system with goods being moved on a public transit line to an intermediate transfer location, from where the packages can be delivered to the end customers. Their problem includes the location decision of the transfer station, apart from the routing decisions of the last mile vehicles. They consider line network and circular network configurations for the transit networks, and customers are assumed to be uniformly distributed around them. The problem is solved using an ALNS heuristic using a Greedy Randomized Adaptive Search Procedure, and they find savings of up to 7.1\% for line networks and 5.4\% for circular networks.

\citet{loy2021combined} address a three-tier transportation problem of integrating passengers and packages. They first study a strategic level problem to select interesting and well-located lines and stops for transporting packages using a modified facility location problem. Then they study the three-tier problem, with the first tier composed of trucks or vans delivering packages from a hub to a departure station. In the second tier, the public transit vehicles transport them to other stops near the customers, and finally, the city freighters use electric cargo bikes to make the final delivery. They also integrate the reverse flow of packages from the urban area to the consolidation centers. They study the problem using three optimization models-- the Uncapacitated Facility Location Problem, the VRP with Simultaneous Delivery and Pickup, extended with time windows, and the Capacitated VRP with Pickup and Deliveries and Time Windows. They use a routing engine called Open Source Routing Machine for routing decisions. While their setting is more general than ours, our study differs in providing a full formulation for the entire problem and then using decomposition techniques to solve it.

In \citet{donne2021freight}, the authors study the FOT problem at the strategic level. They discuss decisions about the public vehicle lines and the drop-in and drop-out stops to be included in the delivery system to maximize the demand covered. They propose several formulations for the problem and employ column generation-based heuristic approaches to solve them. They find that the topology of the public transportation network, along with the demand distribution and density, play the most significant role in achieving their objective.

Finally, a recent detailed survey has appeared on freight delivery on urban public transportation systems \citep{elbert2021freight}. We refer the reader to this exhaustive review for a complete view of the literature on the topic.

\section{Problem Setting} \label{problem_formulation}

Let $\mathcal{C}$ be the set of all the customers where a package has to be delivered. Let $o$ denote the CDC or the parcel depot, and $o'$ denote a copy of the CDC. Let the set of dedicated delivery trucks at the CDC be denoted by $\mathcal{D}$. Let $\mathcal{P}$ be the set of public vehicles that can be used for delivering packages. Let $\mathcal{S}$ be the set of all the stops that have been equipped for the delivery system. For each vehicle $p \in \mathcal{P}$, let $\mathcal{S}_{p}$ denote the set of its stops. In our notation, each public vehicle has a unique representation, and it can be identified by the route (or line) it serves and the time it starts from its depot. Furthermore, let $\mathcal{S}_{in}$ and $\mathcal{S}_{out}$ denote the drop-in and the drop-out stops on the public vehicle network, respectively, so we have $\mathcal{S} = \mathcal{S}_{in} \cup \mathcal{S}_{out}$. For each customer $i$, we pre-assign a set of drop-out stops that can potentially serve the customers based on the distance from the stop to them. This is denoted by the set $\mathcal{S}_{i}^{out}$. These drop-out stops lead to a set of drop-in stops from where customer $i$'s package can be transported (the drop-in stops that are on the line of the public vehicles that serve the drop-out stops in $\mathcal{S}_{i}^{out}$) and is denoted by  $\mathcal{S}^{in}_{i}$. On the other hand, the set of customers that can be served from each drop-out stop $v$ and each drop-in stop $u$ are denoted by $C^{out}_{v}$ and $C^{in}_{u}$, respectively. Let $\mathcal{P}_{s}$ be the set of all public vehicles that visit stop $s$. Let $\mathcal{K}$ be the set of all city freighters, and for each stop, $s\in \mathcal{S}$, let $\mathcal{K}_{s}$ denote the set of freighters that serve stop $s$. In our setting, each freighter serves exactly one drop-out stop. $\mathcal{N}$ denotes the set of all nodes or locations in the delivery system -- the customers, the stops, and the CDC. Finally, we introduce a hypothetical node $\tilde{o}$ that represents the depot for all public vehicles. For simplicity of modeling, we assume that all public vehicles start from the same imaginary depot $\tilde{o}$ before reaching their first drop-in stop. This is because we are not concerned about the route of a public vehicle before it reaches its first drop-in stop.

We consider a daily planning horizon here. Each customer $i$ has to be delivered a package, which consumes a capacity $q_{i}$. The packages are delivered from the depot $o$ by a delivery truck $d$, and dropped off at a stop $u \in \mathcal{S}_{in}$. The capacity of each delivery truck is given by $Q^{1}_{d}$. Then a public vehicle $p$ collects the packages from the drop-in stop $u$ and drops them off at the drop-out stop $v \in \mathcal{S}_{out}$. Let $Q^{2}_{p}$ be the capacity of each vehicle $p$. Let $T_{sp}$ denote the time at which a public vehicle $p$ visits a stop $s \in \mathcal{S}_{p}$ on its line. Each drop-out stop $s$ is served by a group of freighters allotted to that stop ($\in \mathcal{K}_{s}$), who pick up the packages dropped off at the stop and deliver them to the respective customers. The capacities of the freighters are given by $Q^{3}_{k}$. Each customer $i$ must be served within their time window $\left[ \underline{T}_{i} , \overline{T}_{i} \right]$. Additionally, we consider that there is some service time, denoted by $T^{'}_{s}$, associated with each of the stops. This includes the time it takes to load the packages onto the vehicles or unload them, the time for the freighters to start their journey, and so on. Each customer $i$ also has a service time $\widehat{T}_{i}$, which refers to the time required to deliver a package at the customer location. It includes the time taken by freighters to find a parking spot and locate the exact apartment or building once they reach the customer, among others.

Let $D_{ij}$ denote the distance between any two locations of the system, where $i, j \in \mathcal{N}$. Let $ T^{1}_{uvd} $ denote the time taken by delivery truck $d$ to traverse arc $(u,v)$, where $u \in \mathcal{S}_{in} \cup \{o\}$ and $v \in \mathcal{S}_{in} \cup \{o'\}$, by a delivery truck $d$, and $C^{1}_{uvd}$ be the cost of using a delivery truck to traverse the arc. Similarly, let $ T^{3}_{ijk} $ denote the travel time for traversing arc $(i,j)$, where, $i, j \in \mathcal{C} \cup \mathcal{S}_{out}$, by freighter $k$, and $C^{3}_{ijk}$ be the corresponding cost. We use the parameter $\alpha_{uvp}$ to identify the public vehicle routes. It takes the value 1 if the vehicle $p$ goes from a stop $u$ to a stop $v$ $\in \mathcal{S}_{p}$, and 0 otherwise. We also have a parameter $V_{s}$ that ensures that a package cannot be left at a stop $s \in \mathcal{S}$ for more than that amount of time. For example, if we need to deliver a package to a customer in the evening, this parameter guarantees that we do not deliver it to a drop-in or a drop-out stop early in the morning. The stops used are usually equipped with some storage facilities, and this parameter ensures that no package consumes storage for more than a certain amount of time.

We have some assumptions on the setting of our problem. Firstly, we assume that the delivery system has already been set up and we are only dealing with the day-to-day operational problem and, thus, the contracts between parties involved have been established. Moreover, the necessary infrastructure like storage units for parcels, dedicated spaces in public vehicles (if required), equipment for transferring packages to and from the vehicles, and the personnel handling these operations have all been established.

We study a deterministic version of the problem here, so the orders for the day are known. The residual capacities of public vehicles is assumed to have been estimated from observing previous occupancy data of the vehicles at the stops. We assume that there are enough delivery trucks and freighters, and enough capacity on the public vehicles to deliver all the packages. Finally, we assume that there is no transfer of packages within the same tier.

\begin{table}[htbp]
  \centering
  \caption{List of sets and parameters used}
  \vspace{-1em}
    \resizebox{16cm}{!}{
    \begin{tabular}{p{3em}p{27em}rrr}
\cmidrule{1-2}\cmidrule{4-5}    \textbf{Sets}  & \multicolumn{1}{r}{} &       & \multicolumn{1}{p{3em}}{\textbf{Parameters}} &  \\
\cmidrule{1-2}\cmidrule{4-5}    $o$ & the CDC or parcel depot &       &     \multicolumn{1}{p{3em}}{$ C^{1}_{uvd} $} & \multicolumn{1}{p{23em}}{ cost of traversing arc $(u,v)$ using a delivery truck} \\
    $o'$ & dummy CDC (a copy of the CDC) &       & \multicolumn{1}{p{3em}}{$ C^{3}_{ijk} $} & \multicolumn{1}{p{23em}}{ cost of traversing arc $(i,j)$ by a freighter}  \\
    $\mathcal{C}$ & set of customers &       & \multicolumn{1}{p{3em}}{$D_{ij}$} & \multicolumn{1}{p{23em}}{distance between locations $i,j \in \mathcal{N}$} \\
    $\mathcal{P}$ & set of public transportation vehicles &       & \multicolumn{1}{p{3em}}{$T_{sp}$} & \multicolumn{1}{p{24em}}{time when public vehicle $p$ is scheduled to reach stop $s$} \\
    $\mathcal{S}$ & set of all stops for public transportation vehicles &       & \multicolumn{1}{p{3em}}{$Q^{1}_{d}$} & \multicolumn{1}{p{23em}}{capacity of delivery truck $d$} \\
    $\mathcal{S}_{in}$ & the set of \textit{drop-in} stops &       & \multicolumn{1}{p{3em}}{$Q^{2}_{p}$} & \multicolumn{1}{p{23em}}{capacity of public vehicle $p$} \\
    $\mathcal{S}_{out}$ & the set of \textit{drop-out} stops &       & \multicolumn{1}{p{3em}}{$Q^{3}_{k}$} & \multicolumn{1}{p{23em}}{capacity of freighter $k$} \\
    $\mathcal{S}_{p}$ & set of stops traversed by public transportation vehicle $p$ &       & \multicolumn{1}{p{3em}}{$q_{i}$} & \multicolumn{1}{p{23em}}{capacity consumed by the package of customer $i$} \\
    $\mathcal{P}_{s}$ & set of public transportation vehicles that visit stop &       & \multicolumn{1}{p{3em}}{$ T^{1}_{uvd} $} & \multicolumn{1}{p{23em}}{time taken by delivery truck $d$ to traverse arc $(u , v)$} \\
    $\mathcal{D}$   & set of delivery trucks at the CDC &       & \multicolumn{1}{p{3em}}{$ T^{3}_{ijk} $} & \multicolumn{1}{p{23em}}{time taken by freighter $k$ to traverse arc $(i , j)$} \\
    $\mathcal{K}$ & set of all city freighters &       & \multicolumn{1}{p{3em}}{$\underline{T}_{i}$} & \multicolumn{1}{p{23em}}{earliest time that customer $i$ can be served} \\
    $\mathcal{K}_{s}$ & set of freighters that serve drop-out stop $s$ &       & \multicolumn{1}{p{3em}}{$\overline{T}_{i}$} & \multicolumn{1}{p{23em}}{latest time that customer $i$ can be served} \\
    $\mathcal{S}_{i}^{out}$ & set of drop-out stops from where a customer $i$ can be served &       & \multicolumn{1}{p{3em}}{$ \widehat{T}_{i} $} & \multicolumn{1}{p{23em}}{service time (delivery time) required at customer $i$} \\
    $\mathcal{S}_{i}^{in}$ & set of drop-in stops where customer $i$'s package can be loaded onto a public vehicle &       & \multicolumn{1}{p{3em}}{$T^{'}_{s}$} & \multicolumn{1}{p{23em}}{service time required at stop $s$} \\
    $C_{s}^{in}$ & the set of customers that can be served via drop-in stop $s$ & &  \multicolumn{1}{p{3em}}{$ V_{s} $} & \multicolumn{1}{p{23em}}{maximum time that a package can be left at stop $s$}\\
    $C_{s}^{out}$ & the set of customers that can be served from drop-out stop $s$ & & \multicolumn{1}{p{3em}}{$\alpha_{uvp}$} & \multicolumn{1}{p{23em}}{equals 1 if public vehicle $p$ goes from stop $u$ to stop $v$, 0 otherwise} \\
    $\mathcal{N}$ & the set of all nodes in the system, equals $\mathcal{C} \cup \mathcal{S} \cup \{ o \} \cup \{ o' \}$ &       & \\
    $\Tilde{o}$ & bus depot (dummy node) &       & \\
\cmidrule{1-2}\cmidrule{4-5}    
\end{tabular}%
}
  \label{tab:notation_summary_1}%
\end{table}%

\begin{table}[htbp]
  \centering
  \caption{List of decision variables for the problem}
  \vspace{-1em}
   \resizebox{16cm}{!}{
     \begin{tabular}{ll}
    \toprule
    \multicolumn{2}{p{55em}}{\textbf{Decision Variables}} \\
    \midrule
    $r_{isd}$ & equals 1 if the package for customer $i$ is delivered by truck $d$ to the drop-in stop $s$, 0 otherwise \\
    $w_{uvd}$ & equals 1 if a delivery truck $d$ traverses arc $(u,v)$, 0 otherwise \\
    $t_{ud}^{1}$ & time when delivery truck $d$ starts from (or leaves) node $u \in \mathcal{S}_{in} \cup \{ o \}$ \\
    $y_{isp}^{1}$ & equals 1 if the package for customer $i$ is picked up by public vehicle $p$ from drop-in stop $s$, 0 otherwise \\
    $y_{isp}^{2}$ & equals 1 if the package for customer $i$ is dropped off by public vehicle $p$ at drop-out stop $s$, 0 otherwise \\
    $l_{up}^{2}$ & load of a public vehicle $p$ as it leaves stop $u \in \mathcal{S}_{p}$ \\
    $z_{ik}$ & equals 1 if freighter $k \in K_{s}$ picks up customer $i$'s package, 0 otherwise \\
    $x_{ijk}$ & equals 1 if freighter $k$ traverses arc $(i,j)$, 0 otherwise \\
    $t_{ik}^{3}$ & time when freighter $k$ starts from node $i \in \mathcal{C} \cup \mathcal{S}_{out}$ \\
    \bottomrule
    \end{tabular}%
    }
  \label{tab:notation_summary_2}%
\end{table}%

Next, we describe the decision variables for the problem. $r_{isd}$ denotes binary variables that take the value 1 if delivery truck $d$ transports the package $i$ from the CDC to the drop-in stop $s$, and 0 otherwise. Binary variables $w_{uvd}$ take the value 1 if delivery truck $d$ traverses arc $(u,v)$, where $u \in \mathcal{S}_{in} \cup \{o\}$ and $v \in \mathcal{S}_{in} \cup \{o'\}$, and 0 otherwise. $t^{1}_{ud}$ is a continuous variable that notes the time as the truck $d$ leaves the CDC (if $u = o$) or visits each drop-in stop $u \in \mathcal{S}_{in}$. These three sets of variables comprise the decisions of the first tier of the problem. Then, we have binary variables $y^{1}_{isp}$ that take value 1 if the public vehicle $p$ picks package $i$ from the drop-in stop $s$, and 0 otherwise. Similarly, binary variables $y^{2}_{isp}$ take value 1 if the public vehicle $p$ drops package $i$ at the drop-out stop $s$. Variables $l^{2}_{up}$ are continuous variables that update the load of the public vehicle $p$ as it visits each of its stops. This load is computed as the initial load of the vehicle $p$ when visiting stop $u$, minus the volume of packages dropped off at $u$, plus the volume of packages picked up from stop $u$. These are the decision variables associated with the second tier of the problem. Finally, we have the variables that are related to the third tier. Binary variables $z_{ik}$ take the value 1 if package $i$ is assigned to freighter $k$, and 0 otherwise. Variables $x_{ijk}$ are also binary and take the value 1 if freighter $k$ traverses arc $(i,j)$, where $i, j \in \mathcal{S}_{out} \cup \mathcal{C}$. Finally, we have continuous variables $t^{3}_{ik}$ that update the time as each freighter $k$ leaves a drop-out stop (if $i \in \mathcal{S}_{out}$) or visits the customers (if $i \in \mathcal{C}$).

Our problem is to find a path for each package from the CDC via a delivery truck, followed by a public vehicle, and finally a freighter to the customer, such that the total distance traveled by the delivery trucks and the freighters is minimized. Tables \ref{tab:notation_summary_1} and \ref{tab:notation_summary_2} summarize the different notations and decision variables for the 3T-DPPT problem.

\begin{figure}[ht]
\centering
\includegraphics[width=11cm]{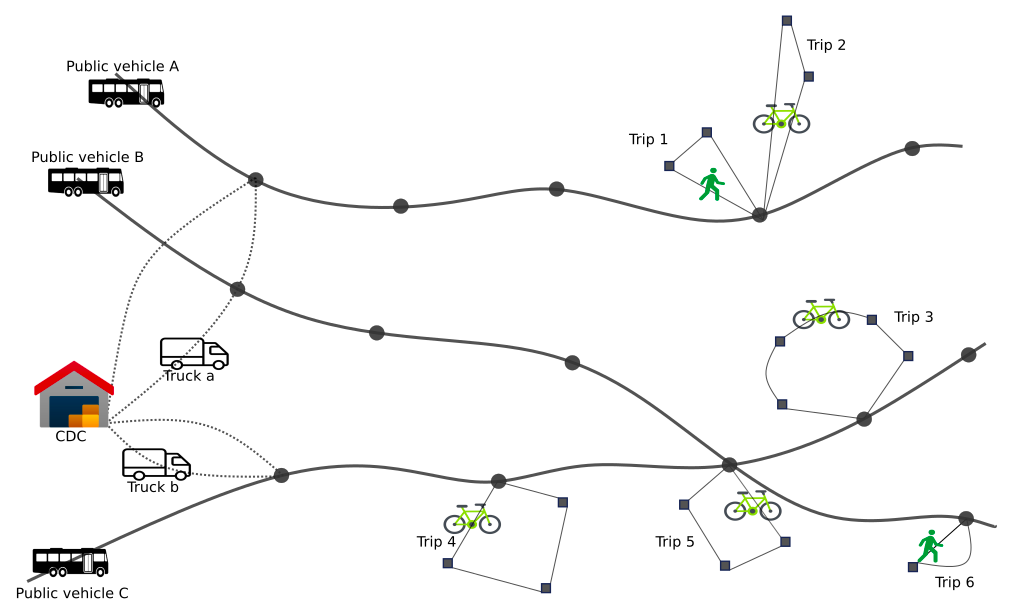}
\vspace{-0.5em}
\caption{\textit{An example of the freight delivery structure}}
\label{flow_example}
\end{figure}

Figure \ref{flow_example} shows a small example of the problem that we study. It consists of a CDC, where the packages originate. The packages are then delivered from the CDC using two delivery trucks, denoted by a and b, to public vehicle stops represented by the small circles. This example consists of three public vehicle lines-- A, B, and C, with corresponding vehicles on them. These lines might have common stops, for example, for lines B and C. Delivery truck $a$ serves lines A and B, while truck $b$ serves line C. The time when the trucks perform their trips depends on when the packages need to be picked up by the vehicles. The vehicles then drop the packages at several stops on their route, from where they are picked up and delivered by the freighters. The packages may arrive on different public vehicles. Once all the packages that need to be delivered in a trip are available at the drop-out stop, the freighters start their route. The freighters then deliver the packages to the customers, denoted by black squares, using green means like bikes or walking.

\FloatBarrier

\subsection{Model Formulation}

The formulation of the 3T-DPPT can be broken down into three parts based on the three tiers of the approach, with constraints linking each part.

The process starts with the delivery trucks transporting packages from the CDC to the drop-in stops in tier 1. Each package is assigned to precisely one delivery truck which delivers it to a drop-in stop so that the truck's capacity is not exceeded. We model this using the following constraints:
\begin{allowdisplaybreaks}
\begin{align}
    & \sum_{s \in \mathcal{S}_{in}} \sum_{d \in \mathcal{D}} r_{isd} = 1 \qquad \forall i \in \mathcal{C}, \label{cons:unique_customer_truck} \\
    & \sum_{i \in C} \sum_{s \in \mathcal{S}_{in}} q_{i} r_{isd} \leq Q^{1}_{d} \qquad \forall d \in \mathcal{D} \label{cons:truck_capacity}.
\end{align}
\end{allowdisplaybreaks}
The next set of constraints defines the route of the trucks. Each truck starts from the CDC, follows a route to deliver the packages to different drop-in stops, and returns to the CDC.
\begin{allowdisplaybreaks}
\begin{align}
    & \sum_{v \in \mathcal{S}_{in} \cup \{o'\}} w_{ovd} = 1 \qquad \forall d \in \mathcal{D}, \label{cons:truck_starting_node} \\
    & \sum_{v \in \mathcal{S}_{in} \cup \{o\}} w_{vo'd} = 1 \qquad \forall d \in \mathcal{D}, \label{cons:truck_ending_node} \\
    & \sum_{\mathclap{\substack{v \in \mathcal{S}_{in} \cup \{o'\} \\ v \neq u}}} \; w_{uvd} = \sum_{\mathclap{\substack{v \in \mathcal{S}_{in} \cup \{o\} \\ v \neq u}}} \; w_{vud} \qquad \forall u \in \mathcal{S}_{in}, \; d \in \mathcal{D}, \label{cons:truck_node_balance} \\
    & \sum_{\mathclap{\substack{u \in \mathcal{S}_{in} \cup \{o\} \\ v \neq u}}} \; w_{uvd} \geq \frac{1}{M} \sum_{i \in \mathcal{C}} r_{ivd} \qquad \forall v \in \mathcal{S}_{in}, \; d \in \mathcal{D}. \label{cons:truck_route_stop_customer_link}
\end{align}
\end{allowdisplaybreaks}
Constraints \eqref{cons:truck_starting_node} and \eqref{cons:truck_ending_node} ensure that the delivery trucks start and end at the CDC. Constraints \eqref{cons:truck_node_balance} are the flow conservation constraints for the delivery trucks at each drop-in stop. Constraints \eqref{cons:truck_route_stop_customer_link} ensure that if a package is assigned to a drop-in stop via a delivery truck, that truck must visit the stop. $M$ denotes a large number.

Then, we use constraints to link tier 1 to tier 2, i.e., if a delivery truck drops off a package at a drop-in stop, a public vehicle must pick it up from that stop, and the constraints are given by:
\begin{allowdisplaybreaks}
\begin{align}
    & \sum_{d \in \mathcal{D}} r_{isd} = \sum_{p \in \mathcal{P}_{s}} y_{isp}^{1} \qquad \forall i \in \mathcal{C}, \; s \in \mathcal{S}_{in}. \label{cons:truck_stop_vehicle_customer_link}
\end{align}
\end{allowdisplaybreaks}
Moreover, we have time constraints on when the packages can be dropped off and picked up at the drop-in stops.
\begin{allowdisplaybreaks}
\begin{align}
    & t_{vd}^{1} \geq t_{ud}^{1} + T^{1}_{uvd} + T^{'}_{v} - M \left( 1 - w_{uvd} \right) \qquad \forall u \in \mathcal{S}_{in} \cup \{o\}, \; v \in \mathcal{S}_{in}, \; u \neq v, \; d \in \mathcal{D}, \label{cons:truck_time_balance_linearized} \\
    & t_{sd}^{1} \leq \sum_{p \in \mathcal{P}_{s}} T_{sp} y_{isp}^{1} + M \left( 1 - r_{isd} \right) \qquad \forall i \in \mathcal{C}, \; s \in \mathcal{S}_{in}, \; d \in \mathcal{D}, \label{cons:truck_time_balance_upper_limit_linearized} \\
    & \sum_{p \in \mathcal{P}_{s}} T_{sp} y_{isp}^{1} - t_{sd}^{1} \leq V_{s} + M \left( 1 - r_{isd} \right) \qquad \forall i \in \mathcal{C}, \; s \in \mathcal{S}_{in}, \; d \in \mathcal{D}. \label{cons:package_pickup_stop_timelimit_linearized}
\end{align}
\end{allowdisplaybreaks}
Constraints \eqref{cons:truck_time_balance_linearized} update the time for each delivery truck as it visits each drop-in stop. Constraints \eqref{cons:truck_time_balance_upper_limit_linearized} ensure that a delivery truck visits the drop-in stops before the packages that it carries are scheduled to be picked up by the public vehicles, and constraints \eqref{cons:package_pickup_stop_timelimit_linearized} provide a limit on the time that a package spends at a drop-in stop.

Next, we model the flow of the packages on the public vehicle network. We have the following constraints:
\begin{allowdisplaybreaks}
\begin{align}
    & \sum_{p \in \mathcal{P}} \; \sum_{s \in \mathcal{S}_{p} \cap \mathcal{S}_{in}} y_{isp}^{1} = 1 \qquad \forall i \in \mathcal{C}, \label{cons:vehicle_package_pick}\\
    & \sum_{p \in \mathcal{P}} \; \sum_{s \in \mathcal{S}_{p} \cap \mathcal{S}_{i}^{out}} y_{isp}^{2} = 1 \qquad \forall i \in \mathcal{C}, \label{cons:vehicle_package_drop} \\
    & y_{isp}^{1} + y_{isp}^{2} \leq 1 \qquad \forall i \in \mathcal{C}, \; s \in \mathcal{S}_{in} \cap \mathcal{S}_{out}, \; p \in \mathcal{P}_{s}, \label{cons:drop-in_stop_different_drop-out_stop} \\
    & \sum_{s \in \mathcal{S}_{p} \cap \mathcal{S}_{in}} y_{isp}^{1} = \sum_{s \in \mathcal{S}_{p} \cap \mathcal{S}_{i}^{out}} y_{isp}^{2} \qquad \forall i \in \mathcal{C}, \; p \in \mathcal{P}, \label{cons:pick_drop_vehicle_stop_link} \\
    % & \text{If} \; y_{ivp}^{2} = 1, \; \text{then} \;\; T_{up} y_{iup}^{1} \leq T_{vp} y_{ivp}^{2} \qquad \forall i \in \mathcal{C}, \; p \in \mathcal{P}, \; u \in \mathcal{S}_{p} \cap \mathcal{S}_{in}, \; v \in \mathcal{S}_{p} \cap \db{\mathcal{S}^{out}_{i}}, \label{cons:pickup_time_less_dropoff_time} \\
    & T_{up} y_{iup}^{1} \leq T_{vp} y_{ivp}^{2} + M \left( 1 - y_{ivp}^{2} \right) \qquad \forall i \in \mathcal{C}, \; p \in \mathcal{P}, \; u \in \mathcal{S}_{p} \cap \mathcal{S}_{in}, \; s \in \mathcal{S}_{p} \cap \mathcal{S}^{out}_{i}, \label{cons:pickup_time_less_dropoff_time_linearized} \\
    & l_{\Tilde{o} p}^{2} = 0 \qquad \forall p \in \mathcal{P}, \label{cons:public_vehicle_load@depot} \\
    & \alpha_{uvp} \left[ l_{vp}^{2} = l_{up}^{2} + \sum_{i \in \mathcal{C}^{in}_{v} } q_{i} y_{ivp}^{1} - \sum_{i \in \mathcal{C}^{out}_{v} } q_{i} y_{ivp}^{2} \right] \qquad \forall p \in \mathcal{P}, \; u \in \mathcal{S}_{p} \cup \{ \tilde{o} \}, \; v \in \mathcal{S}_{p}, \; u \neq v, \label{cons:vehicle_capacity_1} \\
    & l_{vp}^{2} \leq Q^{2}_{p} \qquad \forall p \in \mathcal{P}, \; v \in \mathcal{S}_{p}. \label{cons:vehicle_capacity_2}
\end{align}
\end{allowdisplaybreaks}
Constraints \eqref{cons:vehicle_package_pick} and \eqref{cons:vehicle_package_drop} state that the package for each customer $i$ is picked up by exactly one public vehicle $p$ from a drop-in stop and dropped off at a drop-out stop. Constraints \eqref{cons:drop-in_stop_different_drop-out_stop} ensure that the drop-in stop of a package is different from its drop-out stop. Constraints \eqref{cons:pick_drop_vehicle_stop_link} make certain that the vehicle $p$ that picks up a package for customer $i$ from a drop-in stop is the same that drops off the package at a drop-out stop. Constraints \eqref{cons:pickup_time_less_dropoff_time_linearized} ensure that a package is picked up by a public vehicle before it is dropped off at one of the drop-out stops on its route. Constraints \eqref{cons:public_vehicle_load@depot} set the load for each public vehicle to be zero initially, constraints \eqref{cons:vehicle_capacity_1} update the load of each public vehicle at each stop, and \eqref{cons:vehicle_capacity_2} guarantee that the capacity of the public vehicles is respected.

Once the packages reach the drop-out stops, we have the third and final tier of the problem. We have some constraints that link T2 and T3, similar to what we had for T1 and T2. Specifically, the following constraints say that if a package is dropped off at a drop-out stop by a public vehicle, then a freighter serving that stop must pick it up:
\begin{allowdisplaybreaks}
\begin{align}
    & \sum_{k \in \mathcal{K}_{s}} z_{ik} = \sum_{p \in \mathcal{P}_{s}} y_{isp}^{2} \qquad \forall i \in \mathcal{C}, \; s \in \mathcal{S}_{i}^{out}. \label{cons:vehicle_stop_freighter_customer_link}
\end{align}
\end{allowdisplaybreaks}
Then, we have  constraints that are exclusive to the third tier and that describe the routes of the freighters:
\begin{allowdisplaybreaks}
\begin{align}
    & \sum_{k \in \mathcal{K}} \sum_{ j \in \mathcal{C} \cup \mathcal{S}_{out} } x_{ijk} = 1 \qquad \forall i \in \mathcal{C}, \label{cons:customer_visited_exactly_once} \\
    & \sum_{ \mathclap{ \substack{ j \in \mathcal{C} \cup \{ s \} \\ j \neq i } } } x_{jik} = z_{ik} \qquad \forall i \in \mathcal{C}, \; s \in \mathcal{S}_{i}^{out}, \; k \in \mathcal{K}_{s}, \label{cons:freighter_customer_link} \\
    & \sum_{i \in \mathcal{C}} q_{i} z_{ik} \leq Q^{3}_{k} \qquad \forall k \in \mathcal{K}, \label{cons:freighter_capacity} \\
    & \sum_{j \in \mathcal{C} \cup \mathcal{S}_{out}} x_{sjk} = 1 \qquad \forall s \in \mathcal{S}_{out}, \; k \in \mathcal{K}_{s}, \label{cons:freighter_starting_stop} \\
    & \sum_{j \in \mathcal{C} \cup \mathcal{S}_{out}} x_{jsk} = 1 \qquad \forall s \in \mathcal{S}_{out}, \; k \in \mathcal{K}_{s}, \label{cons:freighter_ending_stop} \\
    & \sum_{ \mathclap{ \substack{ j \in \mathcal{C} \cup \mathcal{S}_{out} \\ j \neq i } } } \; x_{ijk} = \sum_{ \mathclap{ \substack{ j \in \mathcal{C} \cup \mathcal{S}_{out} \\ j \neq i } } } \; x_{jik} \qquad \forall i \in \mathcal{C}, \; k \in \mathcal{K}. \label{cons:freighter_node_balance}
\end{align}
\end{allowdisplaybreaks}
Constraints \eqref{cons:customer_visited_exactly_once} ensure every customer is visited exactly once. Constraints \eqref{cons:freighter_customer_link} establish the link between variables $x$ and $z$. Constraints \eqref{cons:freighter_capacity} make sure that each freighter's capacity is respected. Constraints \eqref{cons:freighter_starting_stop}-\eqref{cons:freighter_ending_stop} state that the freighters start and end at their own pre-specified drop-out stop. Constraints \eqref{cons:freighter_node_balance} are flow conservation constraints for the freighters at each customer. The freighters start from a drop-out stop, follow a route to visit the customers, and come back to the same drop-out stop.

Finally, we have constraints that note the time that the freighters visit the customers and ensure that the packages reach the customers within their time windows:
\begin{allowdisplaybreaks}
\begin{align}
    & t_{sk}^{3} \geq \sum_{p \in \mathcal{P}_{s}} T_{sp} y_{isp}^{2} + T^{'}_{s} - M \left( 1 - z_{ik} \right) \qquad \forall i \in \mathcal{C}, \; s \in \mathcal{S}_{i}^{out}, \; k \in \mathcal{K}_{s}, \label{cons:freighter_time_balance_1stnode_linearized} \\
    & t_{jk}^{3} \geq t_{ik}^{3} + T^{3}_{ijk} + \widehat{T}_{j} - M \left(1 - x_{ijk} \right) \qquad \forall i \in \mathcal{C} \cup \mathcal{S}_{out}, \; j \in \mathcal{C}, \; i \neq j, \; k \in \mathcal{K}, \label{cons:freighter_time_balance_linearized} \\
    & t_{ik}^{3} \geq \underline{T}_{i} - M \left( 1 - z_{ik} \right) \qquad \forall i \in \mathcal{C}, \; k \in \mathcal{K}, \label{cons:customer_time_limits_linearized_1} \\
    & t_{ik}^{3} \leq \overline{T}_{i} + M \left( 1 - z_{ik} \right) \qquad \forall i \in \mathcal{C}, \; k \in \mathcal{K}, \label{cons:customer_time_limits_linearized_2} \\
    & t_{sk}^{3} - \sum_{p \in \mathcal{P}_{s}} y_{isp}^{2} T_{sp} \leq V_{s} + M \left( 1 - z_{ik} \right) \qquad \forall i \in \mathcal{C}, \; s \in \mathcal{S}_{out}, \; k \in \mathcal{K}_{s}. \label{cons:package_dropoff_stop_timelimit_linearized}
\end{align}
\end{allowdisplaybreaks}
Constraints \eqref{cons:freighter_time_balance_1stnode_linearized}-\eqref{cons:freighter_time_balance_linearized} update the time of the freighters as each customer is visited, and constraints \eqref{cons:customer_time_limits_linearized_1} and \eqref{cons:customer_time_limits_linearized_2} ensure that the time windows of the customers are satisfied. We use constraints \eqref{cons:package_dropoff_stop_timelimit_linearized} to make sure that we do not leave the packages at the drop-out stops for more than a certain amount of time. Constraints \eqref{cons:truck_time_balance_linearized} and \eqref{cons:freighter_time_balance_linearized} also help prevent subtours for the trucks and the freighters, respectively.

Apart from the above constraints, we use the following standard symmetry-breaking constraints:
\begin{allowdisplaybreaks}
\begin{align}
    & \sum_{v \in \mathcal{S}_{in}} w_{ovd} \geq \sum_{v \in \mathcal{S}_{in}} w_{ovd+1} \qquad \forall d \in |\mathcal{D}|-1, \label{cons:symmetry_breaking_tier1_1} \\
    & \sum_{u \in \mathcal{S}_{in} \cup \{ o \} } \sum_{ \substack{ v \in \mathcal{S}_{in} \cup \{ o' \} \\ u \neq v } } w_{uvd} \geq \sum_{u \in \mathcal{S}_{in} \cup \{ o \} } \sum_{ \substack{ v \in \mathcal{S}_{in} \cup \{ o' \} \\ u \neq v } } w_{uvd+1} \qquad \forall d \in |\mathcal{D}|-1, \label{cons:symmetry_breaking_tier1_2} \\
    & \sum_{j \in \mathcal{C}} x_{sjk} \geq \sum_{j \in \mathcal{C}} x_{sjk+1} \qquad \forall s \in \mathcal{S}_{out}, \; k \in |\mathcal{K}_{s}|-1, \label{cons:symmetry_breaking_tier3_1} \\
    & \sum_{i \in \mathcal{C} \cup \{ s \} } \sum_{j \in \mathcal{C} } x_{ijk} \geq \sum_{i \in \mathcal{C} \cup \{ s \} } \sum_{j \in \mathcal{C} } x_{ijk+1} \qquad \forall s \in \mathcal{S}_{out}, \; k \in |\mathcal{K}_{s}|-1. \label{cons:symmetry_breaking_tier3_2} 
\end{align}
\end{allowdisplaybreaks}
Constraints \eqref{cons:symmetry_breaking_tier1_1} say that delivery trucks with smaller indices have to be used first and \eqref{cons:symmetry_breaking_tier1_2} state that trucks with smaller indices have to be assigned larger routes. Constraints \eqref{cons:symmetry_breaking_tier3_1} and \eqref{cons:symmetry_breaking_tier3_2} establish the same, respectively, for the freighters.

The objective of our problem is to minimize the cost of using delivery trucks and freighters in the first and the third tier of the system, respectively and is given by:
\begin{allowdisplaybreaks}
\begin{align}
    \textit{Minimize} \quad & \sum_{d \in \mathcal{D}} \; \sum_{u \in \mathcal{S}_{in} \cup \{o \} } \sum\limits_{ \substack{ v \in \mathcal{S}_{in} \cup \{o' \} \\ u \neq v } } C^{1}_{uvd} w_{uvd} \; + \; \sum_{k \in \mathcal{K}} \; \sum_{i \in \mathcal{S}_{out} \cup \mathcal{C} } \sum\limits_{ \substack{ j \in \mathcal{S}_{out} \cup \mathcal{C} \\ i \neq j } } C^{3}_{ijk} x_{ijk}. \label{obj:original_model}
\end{align}
\end{allowdisplaybreaks}
We denote by \textit{FULL} the complete mathematical formulation of the 3T-DPPT given by the objective function \eqref{obj:original_model} subject to constraints \eqref{cons:unique_customer_truck}-\eqref{cons:symmetry_breaking_tier3_2}.

% \sout{Finally, note that we can explicitly set some variables to zero due to the setting of our problem. The following constraints describe them:}
\begin{comment}
\begin{allowdisplaybreaks}
\begin{align}
    & y^{1}_{isp} = 0 \qquad \forall i \in \mathcal{C}, \; s \in \mathcal{S}_{in}, \; p \in \mathcal{P} \setminus \mathcal{P}_{s}, \label{cons:vehicle_package_pick_redundants} \\
    & y^{2}_{isp} = 0 \qquad \forall i \in \mathcal{C}, \; s \in \mathcal{S}_{out} \setminus \mathcal{S}_{i}^{out}, \; p \in \mathcal{P}, \label{cons:vehicle_package_drop_redundants_1} \\
    & y^{2}_{isp} = 0 \qquad \forall i \in \mathcal{C}, \; s \in \mathcal{S}_{i}^{out}, \; p \in \mathcal{P} \setminus \mathcal{P}_{s}, \label{cons:vehicle_package_drop_redundants_2} \\
    & z_{ik} = 0 \qquad \forall i \in \mathcal{C}, \; s \in \mathcal{S}_{out} \setminus \mathcal{S}_{i}^{out}, \; k \in \mathcal{K}_{s}, \label{cons:freighter_customer_link_redundants} \\
    & x_{sik} = 0 \qquad \forall i \in \mathcal{C}, \; s \in \mathcal{S}_{out}, \; k \in \mathcal{K} \setminus \mathcal{K}_{s}, \label{cons:freighter_serve_respective_drop-off} \\
    & x_{ijk} = 0 \qquad \forall i \in \mathcal{C}, j \in \mathcal{S}_{out} \cup \mathcal{C}, \; s \in \mathcal{S}_{out} \setminus \mathcal{S}_{i}^{out}, \; k \in \mathcal{K}_{s}. \label{cons:customer_route_freighter_redundants}
\end{align}
\end{allowdisplaybreaks}
\end{comment}
% \sout{These constraints can be added to \textit{FULL} to set the redundant variables to zero or removed from the definition of the problem completely.}

Variables $r_{isd}$, $w_{uvd}$, $y^{1}_{isp}$, $y^{2}_{isp}$, $z_{ik}$, and $x_{ijk}$ are binary, while $t^{1}_{sd}$, $l^{2}_{up}$, and $t^{3}_{sk}$ are non-negative continuous variables.

The full model is huge and computationally challenging to solve. The commercial solver used in the computational experiments struggles to find a feasible solution for instances with as few as 40 customers within the time limit. In order to solve larger-sized instances, we employ a heuristic methodology where the problem is decomposed into its natural three tiers and solved individually and sequentially. We describe this solution approach in detail in the next section.

\section{Decomposition Methodology} \label{decomposition_approaches}

In order to solve the 3T-DPPT by decomposition, we first formulate the problems for T1, T2, and T3 individually. We obtain three decomposition approaches, each of which prioritizes one of the tiers and solves it first. We aim to analyze them and identify when one performs better than the others. We could start from T1, then solve T2, and finally T3; we call this approach \textit{D1} (the \textit{D} stands for decomposition, and the number represents which tier is solved first). Alternatively, we could do the reverse and solve T3 first, followed by T2 and T1 (called \textit{D3}). Finally, we could start with T2 and then solve T1 and T3 simultaneously (called \textit{D2}). Each of these methods has its own benefits and challenges, which we describe in detail in the following subsections.

\subsection{Decomposition Starting with Tier 2} \label{decomposition:D2}

We first discuss the case when we handle the problem using a decomposition approach that first solves T2. Since T2 is the tier that links the deliveries in T1 and T3, and our focus is on the use of public transport services, we first solve the T2 problem in this case. We use constraints from \textit{FULL} primarily concerning this tier, along with some additional linking constraints. These additional constraints take into account the delivery time windows of the customers, which are addressed directly only in T3 of the problem. The idea is to solve the problem in T2 by assigning a drop-in stop, a drop-out stop, and a public vehicle that traverses them, to each package, along with the pickup and drop-off time. Then we feed these decisions to the problems in T1 and T3 to get the final complete solution. 

\subsubsection{Tier 2 Problem} \label{decomposition:D2_T2}

The first step of this decomposition technique is to model the flow of packages on public vehicles. For each package $i \in \mathcal{C}$, we consider the task of sending it from a drop-in stop $u$ to a drop-out stop $v$ via a public vehicle $p$. We use variables $y^{1}_{isp}$, $y^{2}_{isp}$, and $l^{2}_{up}$ from the original formulation. The majority of the constraints used for this tier are also the same as in the \textit{FULL} formulation and are given by \eqref{cons:vehicle_package_pick}-\eqref{cons:vehicle_capacity_2}.

Apart from the constraints mentioned above, we also need to make sure that the timing of the assignment of the packages to the public vehicles is such that there are no time inconsistencies for deliveries in T1 and T3. In other words, firstly, the packages have to be assigned to the drop-in stops and the public vehicles so that the delivery trucks have enough time to deliver them from the CDC to the drop-in stop before a public vehicle is scheduled to pick them up. Secondly, we must also ensure that the freighters have enough time to make the deliveries in T3 to satisfy the customers' time windows. We achieve these by using the following constraints:
\begin{align}
    & T_{sp} y_{isp}^{1} \geq  \left[ T^{1avg}_{os} + T^{'}_{s} \right] y_{isp}^{1} \qquad \forall i \in \mathcal{C}, \; s \in \mathcal{S}_{in}, \; p \in \mathcal{P}_{s}, \label{nf1:t1_truck_delivery_lower_bound_time} \\
    & \sum_{p \in \mathcal{P}} \sum_{s \in \mathcal{S}_{p} \cap \mathcal{S}_{i}^{out}} \left[ T_{sp} + T^{3avg}_{si} + T^{'}_{s} + \widehat{T}_{i} \right] y_{isp}^{2} \leq \overline{T}_{i} \qquad \forall i \in \mathcal{C} \label{nf1:package_timing}.
\end{align}
$T^{1avg}_{os}$ denotes the average time taken by a delivery truck to go from the CDC to the drop-in stop $s$, and $T^{3avg}_{si}$ denotes the average time that a freighter requires to travel from drop-out stop $s$ to customer $i$. The first set of constraints ensures that the time when the delivery trucks deliver a package to a drop-in station does not exceed the scheduled time when it is supposed to be picked up by a public vehicle. The second set of constraints ensures that the package can be delivered to its customer before their time window ends.

Note here that we do not introduce an (approximate) constraint concerning the lower bound on the time window for a customer. This is not only because it is hard to estimate the time when a freighter visits the customer, taking into account the maximum waiting time $V_{s}$ and the time taken by the freighter to visit other potential customers before that particular customer, but also because our formulation allows waiting time at the destinations. This means that if a package is available at a drop-out stop at a time such that it is too early to visit the customer, even after waiting for $V_{s}$ amount of time at the drop-out stop taking into consideration his or her lower bound on the time window, the freighter may pick up the package, travel to the customer's location and wait there before it is time to drop the package off.

Next, we discuss the objective functions for the problem in T2. Recall that in the \textit{FULL} formulation, we only had objectives corresponding to minimizing the travel distances in T1 and T3, and nothing specific to T2. We develop three objective functions for the T2 problem, keeping in mind our objective of the original formulation- to minimize the routing costs. The three objectives we propose are: 1) to minimize the number of drop-in stops or drop-out stops or both used, 2) to minimize the approximate routing distances in T1 and T3, and 3) to minimize the approximate number of freighters used in T3. For the third objective function, we focus only on T3 because the third tier is the most complicated one: the distance traveled for delivery in T3 is much higher than T1, firstly, due to a significantly larger number of customers than drop-in stops, and secondly, due to lower capacities of freighters. We describe the objectives in detail in the following paragraphs.

\textbf{Minimizing the number of drop-in and drop-out stops (\textit{Obj1}):} The first objective is to minimize the number of drop-in and drop-out stops used. We develop this objective to encourage high consolidation of packages on the delivery trucks and freighters, thereby using a lesser number of each vehicle. This objective function is given by
\begin{align}
    \textit{Minimize} \; \sum_{s \in \mathcal{S}_{in}} \min \Big\{ \sum_{i \in \mathcal{C}} \sum_{p \in \mathcal{P}} y_{isp}^{1} , 1 \Big\} + \sum_{s \in \mathcal{S}_{out}} \min \Big\{ \sum_{i \in \mathcal{C}} \sum_{p \in \mathcal{P}} y_{isp}^{2} , 1 \Big\} \label{nf1:obj_min_pick_drop_stops}.
\end{align}
This objective function is non-linear. To linearize it, we use binary variables $\phi_{s}^{1}$ and $\phi_{s}^{2}$, which represents if some package is picked up or dropped off at stop $s$, respectively. Then the objective function can be written linearly as
\begin{align}
    \textit{Minimize} \; \sum_{s \in \mathcal{S}_{in}} \phi_{s}^{1} + \sum_{s \in \mathcal{S}_{out}} \phi_{s}^{2} \label{nf1:obj_min_pick_drop_stops_linearized},
\end{align}
with two additional constraints:
\begin{align}
    & \phi_{s}^{1} \geq \frac{\sum_{i \in \mathcal{C}} \sum_{p \in \mathcal{P}} y_{isp}^{1}}{M} \qquad \forall s \in \mathcal{S}_{in}, \label{nf1:min_pick_drop_stops_linearized_cons1} \\
    & \phi_{s}^{2} \geq \frac{\sum_{i \in \mathcal{C}} \sum_{p \in \mathcal{P}} y_{isp}^{2}}{M} \qquad \forall s \in \mathcal{S}_{out}. \label{nf1:min_pick_drop_stops_linearized_cons2}
\end{align}

A drawback of this objective is that it might lead to the selection of very few drop-in and drop-out stations, and assign customers to drop-out stops that might not be the closest to them, and as a result, increase routing distances.

\textbf{Minimizing the approximate routing cost (distance traveled) in T1 and T3 (\textit{Obj2})}: We construct the second objective function by taking into account the distances of the drop-in stops from the CDC  and the drop-out stops from the customers. Though these distances will not be precisely equal to the routing distances obtained by solving the T1 and T3 problems, they serve as a proxy for the objectives of T1 and T3. The second objective function is given by:
\begin{align}
    \textit{Minimize} \; \sum_{i \in \mathcal{C}} \sum_{s \in \mathcal{S}_{in}} D_{os} \sum_{p \in P} y^{1}_{isp} + \sum_{i \in \mathcal{C}} \sum_{s \in \mathcal{S}_{out}} D_{is} \sum_{p \in P} y^{2}_{isp} \label{nf1:obj_min_routing_cost}.
\end{align} 

\textbf{Minimizing the approximate number of freighters used in T3 (\textit{Obj3}):} Finally, we introduce the objective of minimizing the number of freighters used in T3, taking into account the average capacity, denoted by $\overline{Q}_{F}$. This is aimed at reducing the overall routing cost of the third tier. The approximation is made in the following way. First, we divide the entire day into time periods over which the public vehicles operate, denoted by $\mathcal{T} = \{ \tau_{0}, \tau_{1}, \ldots, \tau_{n} \}$. Let $\mathcal{P}^{\tau_{j}}_{s}$ denote the set of public vehicles that visit drop-out station $s$ during period $\tau_{j}$.

To estimate the number of freighters required, we introduce new decision variables $h_{s}^{\tau_{j}}$, which give the number of freighters required during any period to deliver the parcels dropped off in that period. We determine the value of $h_{s}^{\tau_{j}}$ through the following constraints:
\begin{align}
    \sum_{i \in \mathcal{C}} \sum_{p \in \mathcal{P}_{\tau_{j}}^{s}} q_{i} y_{isp}^{2} \leq \overline{Q}_{F} h^{\tau_{j}}_{s} \qquad \forall s \in \mathcal{S}_{out}, \; \tau_{j} \in \mathcal{T} \label{nf1:freighter_estimate},
\end{align}
and the objective function is:
\begin{align}
    \textit{Minimize} \; \sum_{\tau_{j} \in \mathcal{T}} \sum_{s \in \mathcal{S}_{out}}  h_{s}^{\tau_{j}} \label{nf1:obj_min_freighters}.
\end{align}

We solve the tier 2 problem with each of the objectives discussed above. From the solution, we obtain a drop-in stop, a drop-in time, a drop-out stop, and a drop-out time associated with each package. Then, we use these as inputs for the T1 and T3 problems. In the following subsections, we describe the formulations of the problems in T1 and T3.

\subsubsection{Tier 1 Problem} \label{decomposition:D2_T1}

For the tier 1 problem, the task is to assign packages to the delivery trucks and determine the routing decisions of the delivery trucks from the CDC to the drop-in stops. Let $B^{in}_{is}$ be a parameter that takes the value 1 if package $i$ is picked up by a public vehicle from drop-in station $s$, and $T_{i}^{in}$ be the corresponding time. We obtain these values from the solution of T2. Let $r_{id}$ be a binary variable that takes the value 1 if package $i$ is assigned to delivery truck $d$, and $g_{dv}$ be a binary variable that denotes if a truck $d$ visits drop-in stop $v$ or not. Apart from these, we use the variables $w_{uvd}$ and $t^{1}_{vd}$, denoting the arc traversal variables and the time variables, respectively, from the \textit{FULL} model. The objective of the T1 model is 
\begin{align}
    \textit{Minimize} \; \; & \sum\limits_{d \in \mathcal{D}} \; \sum\limits_{u \in \mathcal{S}_{in} \cup \{o \} } \sum\limits_{ v \in \mathcal{S}_{in} \cup \{o' \} } C^{1}_{uvd} w_{uvd} \label{t1:obj},
\end{align}
and the constraints are:
\begin{allowdisplaybreaks}
\begin{align}
    & \sum\limits_{d \in \mathcal{D}} r_{id} = 1 \qquad \forall i \in \mathcal{C}, \label{t1:unique_customer_truck} \\
    & \sum\limits_{i \in \mathcal{C}} q_{i} r_{id} \leq Q^{1}_{d} \qquad \forall d \in \mathcal{D}, \label{t1:truck_capacity} \\
    & g_{dv} \geq B_{iv}^{in} r_{id} \qquad \forall i \in \mathcal{C}, d \in \mathcal{D}, v \in \mathcal{S}_{in} \label{t1:customer_truck_stop_link}, \\
    & g_{dv} = \sum\limits_{\mathclap{\substack{u \in \mathcal{S}_{in} \cup \{o\} \\ u \neq v}}} w_{uvd} \qquad \forall d \in D, \; v \in \mathcal{S}_{in}, \label{t1:truck_stop_visit_link} \\
    & t_{vd}^{1} \leq T^{in}_{i} + M ( 1 - r_{id} ) \qquad \forall i \in \mathcal{C}, \; v \in \mathcal{S}_{in}, \; s.t. \; B_{iv}^{in} = 1, \; d \in \mathcal{D}, \label{t1:truck_visit_time_limit_linearized} \\
    & B^{in}_{is} T^{in}_{i} g_{ds} - t_{sd}^{1} \leq V_{s} + M ( 1 - r_{id} ) \qquad \forall i \in \mathcal{C}, \; s \in \mathcal{S}_{in}, \; d \in \mathcal{D}, \label{t1:package_pickup_stop_timelimit_linearized} 
\end{align}
\end{allowdisplaybreaks}
plus constraints \eqref{cons:truck_starting_node}-\eqref{cons:truck_node_balance}, \eqref{cons:truck_time_balance_linearized}, and \eqref{cons:symmetry_breaking_tier1_1}-\eqref{cons:symmetry_breaking_tier1_2} from the \textit{FULL} model. Constraints \eqref{t1:unique_customer_truck} and \eqref{t1:truck_capacity} are analogous to constraints \eqref{cons:unique_customer_truck} and \eqref{cons:truck_capacity}, respectively. Constraints \eqref{t1:customer_truck_stop_link} say that if a delivery truck is assigned to a package, the truck must visit the corresponding drop-in station for the package. Constraints \eqref{t1:truck_stop_visit_link} link variables $g_{dv}$ to the arc variables $w_{uvd}$. Constraints \eqref{t1:truck_visit_time_limit_linearized} guarantee that a package is delivered to its respective drop-in station before a public vehicle is scheduled to pick it up. Finally, constraints \eqref{t1:package_pickup_stop_timelimit_linearized} ensure that any package is not left at the drop-in stops for more than $V_{s}$ amount of time. 

\subsubsection{Tier 3 Problem} \label{decomposition:D2_T3}

For the tier 3 problem, we need to determine the route of the freighters to deliver the packages from the drop-out stations to the respective customers. Let $B_{is}^{out}$ denote a parameter that takes the value 1 if package $i$ is dropped off at station $v$ by a public vehicle, and let $T_{i}^{out}$ be the corresponding time. These parameters, once again, are determined from the output of T2. For writing this formulation, we use variables $z_{ik}$, $x_{ijk}$, and $t_{jk}^{3}$ from the \textit{FULL} model.

An interesting fact about the problem in tier 3 for the decomposition approach is that it separates into individual problems according to each drop-out stop. Once $B_{is}^{out}$ is known, we know the customers that need to be delivered from a drop-out stop $s$ precisely and can, thus, solve the routing problem for each stop individually. Let $\mathcal{C}_{s}$ be the set of customers scheduled to be served from drop-out stop $s$. The tier 3 model, for each drop-out stop $s$, is given by:
\begin{allowdisplaybreaks}
\begin{align}
    \textit{Minimize} \; \; & \sum_{k \in \mathcal{K}_{s}} \sum\limits_{i \in \mathcal{C}_{s} \cup \{s\} } \sum\limits_{ \substack{ j \in \mathcal{C}_{s} \cup \{s\} \\ j \neq i } } \; C^{3}_{ijk} x_{ijk}, \label{t3_stopwise:obj} \\
    \textit{subject to:} \; \; & \sum_{k \in \mathcal{K}_{s}} z_{ik} = 1 \qquad \forall i \in \mathcal{C}_{s}, \label{t3_stopwise:customer_unique_freighter} \\
    & \sum\limits_{ \substack{ j \in \mathcal{C}_{s} \\ i \neq j } } x_{jik} = z_{ik} \qquad \forall i \in \mathcal{C}_{s} \cup \{ s \}, \; \forall k \in \mathcal{K}_{s}, \label{t3_stopwise:freighter_customer_link} \\
    & \sum_{i \in \mathcal{C}_{s}} q_{i} z_{ik} \leq Q^{3}_{k} \qquad \forall k \in \mathcal{K}_{s}, \label{t3_stopwise:freighter_capacity} \\
    & \sum_{j \in \mathcal{C}_{s} \cup \{s\} } x_{sjk} = 1 \qquad \forall k \in \mathcal{K}_{s}, \label{t3_stopwise:freighter_starting_stop} \\
    & \sum_{j \in \mathcal{C}_{s} \cup \{s\} } x_{jsk} = 1 \qquad \forall k \in \mathcal{K}_{s}, \label{t3_stopwise:freighter_ending_stop} \\
    & \sum\limits_{ \substack{ j \in \mathcal{C}_{s} \cup \{s\} \\ j \neq i } } x_{ijk} = \sum\limits_{ \substack{ j \in \mathcal{C}_{s} \cup \{s\} \\ j \neq i } } x_{jik} \qquad \forall i \in \mathcal{C}_{s}, \; \forall k \in \mathcal{K}_{s}, \label{t3_stopwise:freighter_node_balance} \\
    & t_{sk}^{3} \geq B_{is}^{out} T_{i}^{out} z_{ik} + T^{'}_{s}, \; \forall k \in \mathcal{K}_{s} \qquad \forall i \in \mathcal{C}_{s} , \label{t3_stopwise:freighter_time_balance_1stnode}\\
    & t_{jk}^{3} \geq t_{ik}^{3} + T^{3}_{ijk} + \widehat{T}_{j} - M (1 - x_{ijk}) \qquad \forall i \in \mathcal{C}_{s} \cup \{s\}, \; j \in \mathcal{C}_{s}, i \neq j, \; \forall k \in \mathcal{K}_{s}, \label{t3_stopwise:freighter_time_balance_linearized} \\
    & t^{3}_{ik} \geq \underline{T}_{i} - M (1 - z_{ik}) \qquad \forall i \in \mathcal{C}_{s}, \; \forall k \in \mathcal{K}_{s}, \label{t3_stopwise:customer_time_limits_linearized_1} \\
    & t^{3}_{ik} \leq \overline{T}_{i} + M (1 - z_{ik}) \qquad \forall i \in \mathcal{C}_{s}, \; \forall k \in \mathcal{K}_{s}, \label{t3_stopwise:customer_time_limits_linearized_2} \\
    & t_{sk}^{3} - B^{out}_{is} T^{out}_{i} \leq V_{s} + M (1 - z_{ik}) \qquad \forall i \in \mathcal{C}_{s}, \; k \in \mathcal{K}_{s}, \label{t3_stopwise:package_dropoff_stop_timelimit_linearized} \\
    & \sum_{j \in \mathcal{C}_{s}} x_{sjk} \geq \sum_{j \in \mathcal{C}_{s}} x_{sjk+1} \qquad \forall k \in |\mathcal{K}_{s}|-1, \label{t3_stopwise:symmetry_breaking_t3_1} \\
    & \sum_{i \in \mathcal{C}_{s} \cup \{s\} } \sum\limits_{ \substack{ j \in \mathcal{C}_{s} \\ j \neq i } } x_{ijk} \geq \sum_{i \in \mathcal{C}_{s} \cup \{s\} } \sum\limits_{ \substack{ j \in \mathcal{C}_{s} \\ j \neq i } } x_{ijk+1} \qquad \forall k \in |\mathcal{K}_{s}|-1. \label{t3_stopwise:symmetry_breaking_t3_2}
\end{align}
\end{allowdisplaybreaks}
The above formulation is very similar to the constraints of \textit{FULL} that pertain to T3, decomposed stop-wise. Constraints \eqref{t3_stopwise:freighter_time_balance_1stnode} ensure that the freighters start their journey only after all the packages to be delivered by them have been dropped off.

We denote by \textit{D2-Obj1} the decomposition approach that solves the T2 problem first, followed by T1 and T3, and uses \textit{Obj1} for the problem in T2. When the objective functions used in T2 is \textit{Obj2} or \textit{Obj3}, we denote it by \textit{D2-Obj2} and \textit{D2-Obj3}, respectively.

\subsection{Decomposition Starting with Tier 1} \label{decomposition:D1}

Next, we discuss the decomposition approach when we solve the 3T-DPPT starting from T1. Since we propose to use environmentally and economically sustainable transportation means in T3, whereas T1 uses dedicated delivery vehicles to transport packages, most of the harmful emissions occur in the first tier. Thus, if our objective is to minimize the emissions caused, it would be more beneficial to prioritize the solution of T1. Here, we solve tier 1 first and make decisions about assigning the packages to drop-in stops and delivery trucks and the routes of the trucks. The T1 solution also guides us about when the packages will arrive at the drop-in stops. We use this as input to solve the T2 problem, where we assign the packages to public vehicles and drop-out stops. Finally, we solve T3 using the solution of T2. Here, we assign the packages to freighters and develop the routes of the freighters.

\subsubsection{Tier 1 Problem} \label{decomposition:D1_T1}

For the tier 1 model, similar to \textit{FULL}, we use variables that involve T1, viz. $r_{isd}$, $w_{uvd}$, and $t^{1}_{sd}$. Variables $y^{1}_{isp}$ from \textit{FULL} are replaced by variables $\gamma^{1}_{is}$ which take the value 1 if a customer $i$ is assigned to drop-in stop $s$, and 0 otherwise. The model is given by:
\begin{allowdisplaybreaks}
\begin{align}
    \textit{Minimize} \quad & \sum_{d \in \mathcal{D}} \; \sum_{u \in \mathcal{S}_{in} \cup \{o \} } \sum\limits_{ \substack{ v \in \mathcal{S}_{in} \cup \{o' \} \\ u \neq v } } C^{1}_{uvd} w_{uvd}, \label{t1_starting:obj} \\
    \textit{subject to:} \quad & \sum_{s \in \mathcal{S}_{i}^{in}} \sum_{d \in \mathcal{D}} r_{isd} = 1 \qquad \forall i \in \mathcal{C}, \label{t1_starting:unique_customer_truck} \\
    & \sum_{i \in C} \sum_{s \in \mathcal{S}_{in}^{i}} q_{i} r_{isd} \leq Q^{1}_{d} \qquad \forall d \in \mathcal{D}, \label{t1_starting:truck_capacity} \\
    & \sum_{d \in \mathcal{D}} r_{isd} = \gamma_{is}^{1} \qquad \forall i \in \mathcal{C}, \; s \in \mathcal{S}_{i}^{in}, \label{t1_starting:truck_stop_vehicle_customer_link} \\
    & \sum_{s \in \mathcal{S}_{i}^{in}} \gamma_{is}^{1} = 1 \qquad \forall i \in \mathcal{C}, \label{t1_starting:vehicle_package_pick}
\end{align}
\end{allowdisplaybreaks}
plus constraints \eqref{cons:truck_starting_node}-\eqref{cons:truck_route_stop_customer_link}, \eqref{cons:truck_time_balance_linearized}, and \eqref{cons:symmetry_breaking_tier1_1}-\eqref{cons:symmetry_breaking_tier1_2} from the \textit{FULL} model. Constraints \eqref{t1_starting:unique_customer_truck} and \eqref{t1_starting:truck_capacity} are analogous to constraints \eqref{cons:unique_customer_truck}, and \eqref{cons:truck_capacity}, restricted to specific drop-in stops based on the drop-out stops from where the customers can be served. Constraints \eqref{t1_starting:truck_stop_vehicle_customer_link} and \eqref{t1_starting:vehicle_package_pick} say that if a package is delivered to a drop-in stop, it must be picked up by a unique public vehicle.

Additionally, we need to ensure that there is enough time to deliver the packages to the customers within their time windows. We do so by adding the following constraints:
\begin{align}
    & t_{sd}^{1} \leq T^{avg}_{i} - 1.3 \cdot D_{is} + M (1 - r_{isd}) \qquad \forall i \in \mathcal{C}, \; s \in \mathcal{S}_{i}^{in}, \; d \in \mathcal{D} \label{t1_starting:truck_time_balance_upper_limit_linearized},
\end{align}
where, $T^{avg}_{i} = \frac{\underline{T}_{i} + \overline{T}_{i}}{2}$. To approximate the time, the idea is to assume that we go directly from the CDC to the customer and reach them in the middle of their delivery window. However, since, in practice, the package would be transported first by a public vehicle and then by a freighter, and could also potentially wait in between, the time taken would be greater. To account for the longer journey, we multiply the time taken by a factor greater than 1. From preliminary experiments, we find 1.3 to be a suitable multiplicative factor.

A drawback of our modeling approach is that all the trucks leave the CDC at the beginning of the distribution operations. They deliver the packages to the drop-in stops in the time it takes for them to make the journey from the CDC to the stops, irrespective of when the packages must be delivered. So, even packages that should be delivered to the customers in the evening reach their drop-in stops quite early. However, since we do not allow packages to be left at any stop for more than $V_{s}$ time, we might incur infeasibilities when solving T2. This is because there might not be enough public vehicles to pick up all the packages dropped off by the trucks within time $V_{s}$. To overcome this, we divide the day into periods and perform a pre-processing that assigns each package to the period of its delivery.

Here, we split the delivery period into two halves, with the mid-point denoted by $t_{mid-day}$. Then, we define parameters, $\tau_{is}^{1}$ and $\tau_{is}^{2}$, such that $\tau_{is}^{1}$ takes the value 1 if we need to deliver package $i$ to drop-in stop $s$ during the first half of the day, and zero otherwise. Similarly, $\tau_{is}^{2}$ equals 1 if we have to deliver the package to the drop-in stop during the second half. Next, we discuss how we determine the value of the parameters regarding which half of the day the package needs delivery. We decide this based on the following rule:
\begin{allowdisplaybreaks}
\begin{align*}
    & \text{If} \; \overline{T}_{i} - V_{s} - 1.3 \cdot T_{si} \leq t_{mid-day}, \; \text{then} \; \tau_{is}^{1} = 1 \qquad \forall i \in \mathcal{C}, s \in \mathcal{S}_{in}, \\
    & \text{If} \; \overline{T}_{i} - V_{s} - 1.3 \cdot T_{si} \geq t_{mid-day}, \; \text{then} \; \tau_{is}^{2} = 1 \qquad \forall i \in \mathcal{C}, s \in \mathcal{S}_{in}.
\end{align*}
\end{allowdisplaybreaks}
Essentially, $\tau_{is}^{2} = 1 - \tau_{is}^{1}$, whenever $r_{isd} = 1$. The above criteria say the following: if the closing time window for a customer is such that taking into account the (approximate) time to reach the customer, including the waiting time at one of the stops, the package has to start its journey before time $t_{mid-day}$, and we set $\tau_{is}^{1}$ to be 1. Otherwise, we assign $\tau_{is}^{2}$ to 1, because it is enough for the package to start its journey from the CDC during the second half of the delivery day.

Then, we add the following constraints to the Tier 1 problem:
\begin{align}
    & t_{sd}^{1} \leq t_{mid-day} + M (1 - \tau_{is}^{1} r_{isd}) \qquad \forall i \in \mathcal{C}, s \in \mathcal{S}_{in}, d \in \mathcal{D}, \label{t1_starting:midday_1_linearized} \\
    & t_{sd}^{1} \geq t_{mid-day} - M (1 - \tau_{is}^{1} r_{isd}) \qquad \forall i \in \mathcal{C}, s \in \mathcal{S}_{in}, d \in \mathcal{D}. \label{t1_starting:midday_2_linearized}
\end{align}
This approach can also be extended to more than two periods, which could be required for a large number of packages.

\subsubsection{Tier 2 Problem} \label{decomposition:D1_T2}

From the solution of T1, we get the drop-in stop for each package $i$, and the time when the package is delivered to the drop-in stop $s$. We denote them by $B_{is}^{in}$ and $T_{i}^{in}$, respectively. These parameters are used as inputs for the T2 problem. Let $\overline{s}_{i}$ denote the drop-in stop that serves customer $i$. In T2, we make decisions about allocating each package to a public vehicle and a drop-out stop, and consequently, the time when the package reaches a drop-out stop. For the problem described here, we replace variables $y^{1}_{isp}$ from \textit{FULL} with $\gamma^{1}_{ip}$ to denote the assignment of packages to public vehicles. We also use variables $y^{2}_{isp}$ and $l^{2}_{up}$ from the original formulation. The constraints are:
\begin{allowdisplaybreaks}
\begin{align}
    & \sum_{p \in \mathcal{P}_{\overline{s}_{i}}} \gamma^{1}_{ip} = 1 \qquad \forall i \in \mathcal{C}, \label{t1_starting:vehicle_package_pick_unique}\\
    & \sum_{\mathcal{P} \setminus \mathcal{P}_{\overline{s}_{i}}} \gamma^{1}_{ip} = 0 \qquad \forall i \in \mathcal{C}, \label{t1_starting:vehicle_package_pick_redundants} \\
    & B^{in}_{is} \gamma_{ip}^{1} + y_{isp}^{2} \leq 1 \qquad \forall i \in \mathcal{C}, \; s \in \mathcal{S}_{in} \cap \mathcal{S}_{out}, \; p \in \mathcal{P}_{s}, \label{t1_starting:drop-in_stop_different_drop-out_stop} \\
    & \gamma_{ip}^{1} = \sum_{v \in \mathcal{S}_{p} \cap \mathcal{S}_{i}^{out}} y_{ivp}^{2} \qquad \forall i \in \mathcal{C}, \; p \in \mathcal{P}, \label{t1_starting:pick_drop_stop_vehicle_stop_link} \\
    & T_{up} B^{in}_{iu} \gamma_{ip}^{1} \leq T_{vp} y_{ivp}^{2} + M (1 - y^{2}_{ivp}) \qquad \forall i \in \mathcal{C}, \; p \in \mathcal{P}, \; u \in \mathcal{S}_{p} \cap \mathcal{S}_{in}, \; v \in \mathcal{S}_{p} \cap \mathcal{S}_{i}^{out}, \label{t1_starting:pickup_time_less_dropoff_time_linearized} \\
    & \alpha_{uvp} \left[ l_{vp}^{2} = l_{up}^{2} + \sum_{i \in \mathcal{C}^{in}_{v} } B^{in}_{iv} q_{i} \gamma_{ip}^{1} - \sum_{i \in \mathcal{C}^{out}_{v} } q_{i} y_{ivp}^{2} \right] \qquad \forall p \in \mathcal{P}, \; u \in \mathcal{S}_{p} \cup \{ \tilde{o} \}, \; v \in \mathcal{S}_{p}, \label{t1_starting:vehicle_capacity_bound} \\
    & \sum_{p \in \mathcal{P}_{s}} T_{sp} B^{in}_{is} \gamma^{1}_{ip} \geq T^{in}_{i} B^{in}_{is} \qquad \forall i \in \mathcal{C}, \; s \in \mathcal{S}_{i}^{in}, \label{t1_starting:time_restriction_buses_from_t1} \\
    & \text{If} \; B^{in}_{is} = 1, \; \text{then} \; \sum_{p \in P_{s}} T_{sp} \gamma^{1}_{ip} \leq T^{in}_{i} + V_{s} \qquad \forall i \in \mathcal{C}, \; s \in \mathcal{S}_{in}, \label{t1_starting:max_time_package_at_stop}
\end{align}
\end{allowdisplaybreaks}
plus constraints \eqref{cons:vehicle_package_drop}, \eqref{cons:public_vehicle_load@depot}, and \eqref{cons:vehicle_capacity_2}. Constraints \eqref{t1_starting:vehicle_package_pick_unique}-\eqref{t1_starting:vehicle_capacity_bound} are similar to that of \textit{FULL}, with variables $y^{1}_{isp}$ replaced by variables $\gamma^{1}_{ip}$. Constraints \eqref{t1_starting:time_restriction_buses_from_t1} ensure that the public vehicles pick up the packages only after the trucks have delivered them to the drop-in stops. Constraints \eqref{t1_starting:max_time_package_at_stop} restrict the maximum time a package can be left at a stop before it needs to be picked up by the public vehicles.

In addition, we need to ensure that there is enough time for the freighters to deliver the packages once the public vehicles have dropped them off at the drop-out stops so that the customers receive the packages before their closing time window. To do so, we add the following constraint:
\begin{align}
    & \sum_{p \in \mathcal{P}} \sum_{s \in \mathcal{S}_{p} \cap \mathcal{S}_{out}} \left[ T_{sp} + T^{3avg}_{si} + T^{'}_{s} + \widehat{T}_{i} \right] y_{isp}^{2} \leq \overline{T}_{i} \qquad \forall i \in \mathcal{C}. \label{t1_starting:time_window_upper_limit}
\end{align}
The objective functions for T2, as described in section \ref{decomposition:D2_T2}, could be to minimize the number of drop-out stops used in T3, or to minimize the approximate distance traveled by the freighters in T3, or to minimize the approximate number of freighters used. Since we have already solved T1 here, the objective functions only correspond to the decisions of tier 3.

\textbf{Minimizing the number of drop-out stops (\textit{Obj1}) :} The objective function is:
\begin{align}
    \textit{Minimize} \; \sum_{s \in \mathcal{S}_{in}} \sum_{s \in \mathcal{S}_{out}} \phi_{s}^{2}, \label{t1_starting:obj_min_pick_drop_stops_linearized}
\end{align}
along with the additional constraint \eqref{nf1:min_pick_drop_stops_linearized_cons2}.

\textbf{Minimizing the approximate routing cost of T3 (\textit{Obj2}) :} The objective function in this case is given by:
\begin{align}
    \textit{Minimize} \; \sum_{i \in \mathcal{C}} \sum_{s \in \mathcal{S}_{out}} D_{is} \sum_{p \in P} y^{2}_{isp} \label{t1_starting:obj_min_routing_cost}.
\end{align}

\textbf{Minimizing the approximate number of freighters used in T3 (\textit{Obj3}) :}
Here, the objective function is given by \eqref{nf1:obj_min_freighters}, along with constraints \eqref{nf1:freighter_estimate}.

The tier 3 model takes as input $B^{out}_{is}$ and $T^{out}_{is}$, which denote the drop-out stop and the drop-out time and come from the output of T2. The T3 model here is the same as in the decomposition approach described in Section \ref{decomposition:D2_T3}. As done previously, we denote by \textit{D1-Obj1}, \textit{D1-Obj2}, and \textit{D1-Obj3} the decomposition approaches where the objectives used in T2 are \textit{Obj1}, \textit{Obj2}, and \textit{Obj3}, respectively.

\subsection{Decomposition Starting with Tier 3} \label{decomposition:D3}

The third and final decomposition technique starts with solving the problem from T3. Here, apart from making decisions about assigning packages to freighters and the routes of the freighters, we also decide the drop-out station for each package. This guides the decisions of assigning the packages to the drop-in stops and public vehicles in T2. Finally, the output of tier 2 can be used to make the decisions in T1. This sequence of solving the problem is beneficial when most of the delivery distance (and hence delivery costs) arise in T3 due to a significantly larger number of customers.

\subsubsection{Tier 3 Problem} \label{decomposition:D3_T3}

For the T3 problem, we introduce additional variables $\gamma_{is}^{2}$, which take the value 1 if package $i$ is to be dropped off at drop-out stop $s$ by a public vehicle, partially analogous to variables $y^{2}_{isp}$. We also use variables $x_{ijk}$, $z_{ik}$ and $t^{3}_{ik}$. The T3 model is given by:
\begin{allowdisplaybreaks}
\begin{align}
    \textit{Minimize} \; \; & \sum_{k \in \mathcal{K}} \; \sum_{i \in \mathcal{S}_{out} \cup \mathcal{C} } \sum\limits_{ \substack{ j \in \mathcal{S}_{out} \cup \mathcal{C} \\ i \neq j } } C^{3}_{ijk} x_{ijk}, \label{t3_starting:obj} \\
    \textit{subject to:} \; \; & \sum_{s \in \mathcal{S}_{i}^{out}} \gamma^{2}_{is} = 1 \qquad i \in \mathcal{C}, \label{t3_starting:package_unique_drop_off} \\
    & \sum_{k \in \mathcal{K}_{s}} z_{ik} = \gamma^{2}_{is} \qquad \forall i \in \mathcal{C}, \; s \in \mathcal{S}_{i}^{out}, \label{t3_starting:vehicle_stop_freighter_customer_link} \\
    & t_{sk}^{3} \geq T^{first}_{s} + T^{'}_{s} - M (1 - x_{sik}) \qquad \forall i \in \mathcal{C}, \; s \in \mathcal{S}_{i}^{out}, \; k \in \mathcal{K}_{s}, \label{t3_starting:freighter_time_balance_1stnode_linearized}
\end{align}
\end{allowdisplaybreaks}
plus constraints \eqref{cons:customer_visited_exactly_once}-\eqref{cons:freighter_node_balance}, \eqref{cons:freighter_time_balance_linearized}-\eqref{cons:customer_time_limits_linearized_2}, and \eqref{cons:symmetry_breaking_tier3_1}-\eqref{cons:symmetry_breaking_tier3_2} from \textit{FULL}. Constraints \eqref{t3_starting:package_unique_drop_off} assign each customer to a drop-out station, constraints \eqref{t3_starting:vehicle_stop_freighter_customer_link} ensure that if a package is dropped at a drop-out stop, then it is assigned to a freighter. Finally, constraints \eqref{t3_starting:freighter_time_balance_1stnode_linearized} update the time as the freighters leave the drop-out stops. $T_{s}^{first}$ in constraint \eqref{t3_starting:freighter_time_balance_1stnode_linearized} is an estimated time when the freighters can start making deliveries. We introduce this parameter to make sure that the trucks and the public vehicles have enough time to deliver the packages before the freighters start their journey. We define $T_{s}^{first}$ as the time when the first public vehicle visits stop $s$.

We obtain the drop-out stops and the time by which the packages must reach them from the solution of tier 3. These are then used as inputs for T2, which we describe below.

\subsubsection{Tier 2 Problem} \label{decomposition:D3_T2}

For tier 2, we use as inputs the solutions from T3-- $B^{out}_{is}$ and $T^{out}_{i}$. $B^{out}_{is}$ denotes the drop-out stop $s$ for each package $i$. $T^{out}_{i}$ is the time when package $i$ is supposed to leave the drop-out stop for its delivery, i.e., the time by which the package must reach the drop-out stop on a public vehicle. Let $\overline{s}_{i}$ denote the drop-out stop for package $i$. The model is similar to tier 2 of \textit{FULL}, with variables $y^{1}_{isp}$ and $l^{2}_{up}$, and variables $\gamma^{2}_{ip}$ instead of $y^{2}_{isp}$. At this step, we decide the assignment of drop-in stops and public vehicles for each package. We introduce variables $\gamma^{2}_{ip}$ that is equal to 1 if vehicle $p$ transports package $i$, and 0 otherwise. The constraints of the problem are the following:
\begin{allowdisplaybreaks}
\begin{align}
    & \sum_{p \in \mathcal{P} \cap \mathcal{P}_{\overline{s}_{i}} } \gamma^{2}_{ip} = 1 \qquad \forall i \in \mathcal{C}, \label{t3_starting_t2:vehicle_package_drop} \\
    & \sum_{p \in \mathcal{P} \setminus \mathcal{P}_{\overline{s}_{i}} } \gamma^{2}_{ip} = 0 \qquad \forall i \in \mathcal{C}, \label{t3_starting_t2:vehicle_package_drop_redundants} \\
    & y_{isp}^{1} + \gamma_{ip}^{2} \leq 1 \qquad \forall i \in \mathcal{C}, \; s \in \mathcal{S}_{in} \cap \mathcal{S}_{out}, \; p \in \mathcal{P}_{s}, \; \text{if } B^{out}_{is} = 1, \label{t3_starting_t2:drop-in_stop_different_drop-out_stop} \\
    & \sum_{s \in \mathcal{S}_{p} \cap \mathcal{S}_{in}} y_{isp}^{1} = \left[ \sum_{s \in \mathcal{S}_{p} \cap \mathcal{S}_{out}} B_{is}^{out} \right] \gamma_{ip}^{2} \qquad \forall i \in \mathcal{C}, \; p \in \mathcal{P}, \label{t3_starting_t2:pick_drop_vehicle_stop_link} \\
    & T_{up} y_{iup}^{1} \leq T_{vp} B_{iv}^{out} \gamma_{ip}^{2} + M (1 - B_{iv}^{out} \gamma_{ip}^{2}) \qquad \forall i \in \mathcal{C}, \; p \in \mathcal{P}, \; u \in \mathcal{S}_{p} \cap \mathcal{S}_{in}, \; v \in \mathcal{S}_{p} \cap \mathcal{S}_{i}^{out}, \label{t3_starting_t2:pickup_time_less_dropoff_time_linearized} \\
    & \alpha_{uvp} \left[ l_{vp}^{2} = l_{up}^{2} + \sum_{i \in \mathcal{C}^{in}_{v} } q_{i} y_{iup}^{1} - \sum_{i \in \mathcal{C}^{out}_{v} } q_{i} B_{iv}^{out} \gamma_{ip}^{2} \right] \qquad \forall p \in \mathcal{P}, \; u \in \mathcal{S}_{p} \cup \{\tilde{o}\}, \; v \in \mathcal{S}_{p}, \label{t3_starting_t2:vehicle_capacity_1} \\
    & T_{sp} B_{is}^{out} \gamma^{2}_{ip} \leq T^{out}_{i} \qquad \forall i \in \mathcal{C}, \; s \in \mathcal{S}_{out}, \; p \in \mathcal{P}_{s}, \label{t3_starting_t2:time_constraint} \\
    & B^{out}_{is} T^{out}_{i} - \sum_{p \in \mathcal{P}_{s}} T_{sp} B_{is}^{out} \gamma^{2}_{ip} \leq B^{out}_{is} V_{s} \qquad \forall i \in \mathcal{C}, \; s \in \mathcal{S}_{i}^{out}, \label{t3_starting_t2:max_time_of_package_at_stop} \\
    & T_{sp} y_{isp}^{1} \geq \left[ T^{1avg}_{os} + T^{'}_{s} \right] y_{isp}^{1} \qquad \forall i \in \mathcal{C}, \; s \in \mathcal{S}_{in}, \; p \in \mathcal{P}_{s}, \label{t3_starting_t2:t1_truck_delivery_lower_bound_time}
\end{align}
\end{allowdisplaybreaks}
plus constraints \eqref{cons:vehicle_package_pick}, \eqref{cons:public_vehicle_load@depot}, and \eqref{cons:vehicle_capacity_2}. Constraints \eqref{t3_starting_t2:vehicle_package_drop}-\eqref{t3_starting_t2:time_constraint} are quite similar to the constraints of \textit{FULL}. Constraints \eqref{t3_starting_t2:max_time_of_package_at_stop} give a bound on the maximum time a package can stay at a drop-out stop. Finally, \eqref{t3_starting_t2:t1_truck_delivery_lower_bound_time} ensure that the trucks have enough time to deliver the packages to the drop-in stops before the public vehicles would pick them up. $T^{1avg}_{os}$ denotes the average time that a truck takes to travel from the origin to a drop-out stop $s$.

However, when we start the decomposition technique from T3, the freighters that deliver at least one package start from their respective drop-out stop at time $T_{s}^{first}$. Thus the parameter $T^{out}_{i}$ takes the value $T_{s}^{first}$ for each package. It is easy to see that this results in infeasibilities in T2 whenever the capacity of the first public vehicle to visit that stop is lower than the sum of the capacities consumed by all the packages to deliver. Merely changing the time $T_{s}^{first}$ to a higher value is not sufficient to resolve the issue. This is because the problem would persist whenever the sum of the volumes of all the packages to be delivered from the stop exceeds the sum of the capacities of the public vehicles visiting the stop before the time mentioned above. Because of how we formulated the tier 3 model, the freighters always start their journey from the stop at whatever time we provide as $T_{s}^{first}$. We do the following to find a way around this issue without complicating the problem further. For each freighter $k$, we take the list of customers it visits in order. Then, we start with the last customer visited by the freighter and set the time of delivery to their closing time window. For the other customers, say $i$, we consider the customer visited after $i$, say $j$, and denote the time when $j$ is visited by $t^{*}_{j}$. Then, the time when the freighter $k$ delivers the package to customer $i$ is given by $ \min \{ \overline{T}_{i} , t^{*}_{j} - T^{3}_{ijk} \} $. Though this approximates the delivery time, the routing cost remains the same because we keep the order of visits intact and only play with the time.

Once again, we can have multiple objectives for T2 similar to the previous decomposition approaches. The first objective is to minimize the number of drop-in stops utilized in T1. The second is to minimize the approximate distance traversed in T1.

\textbf{Minimizing the number of drop-in stops (\textit{Obj1}):} The objective function is:
\begin{align}
    \textit{Minimize} \; \sum_{s \in \mathcal{S}_{in}} \phi_{s}^{1}, \label{t3_starting:obj_min_drop-in_stops_linearized}
\end{align}
along with constraints \eqref{nf1:min_pick_drop_stops_linearized_cons1}.

\textbf{Minimizing the approximate routing cost of T1 (\textit{Obj2}):} Finally, in this case, the objective function is given by:
\begin{align}
    \textit{Minimize} \; \sum_{i \in \mathcal{C}} \sum_{s \in \mathcal{S}_{in}} D_{os} \sum_{p \in P} y^{1}_{isp} \label{t3_starting:obj_min_routing_cost}.
\end{align}

We do not need the third objective function, minimizing the approximate number of freighters being used, here, like in the previous approaches, since \textit{Obj3} corresponds exclusively to T3, and we have already solved it.

The T1 problem here is the same as the T1 model in \textit{D2} as described in Section \ref{decomposition:D2_T1}. It takes as input the same parameters as described there and gives the same outputs of assigning the packages to the delivery trucks and their routes. Depending on the objective used for T2, we refer to the decomposition approaches as \textit{D3-Obj1} and \textit{D3-Obj2}, respectively.

Let us also fix a quick notation here to easily refer to the models of a specific tier for the different decomposition approaches. We call them as \textit{decomposition approach number-tier number}. For example, the T3 model from \textit{D2} would be referred to as \textit{D2-T3}.

\section{Numerical Experiments} \label{numerical_experiments}

This section describes how we build a typical instance for the 3T-DPPT problem comprising the public vehicle network, the customers, and the CDC. We then solve our models on these instances, and provide computational results and analysis. We implemented the models on the Spyder IDE of Python 3.8.5 and solved them using the CPLEX 22.1.0 standard solver. The tests were conducted on an Intel(R) Xeon(R) W-2255 CPU with a clock speed of 3.70 GHz and 128 GB RAM.

\subsection{Instance Generation} \label{instance_generation}

\begin{figure}[ht!]
\centering
\vspace{-1em}
\includegraphics[width=11cm]{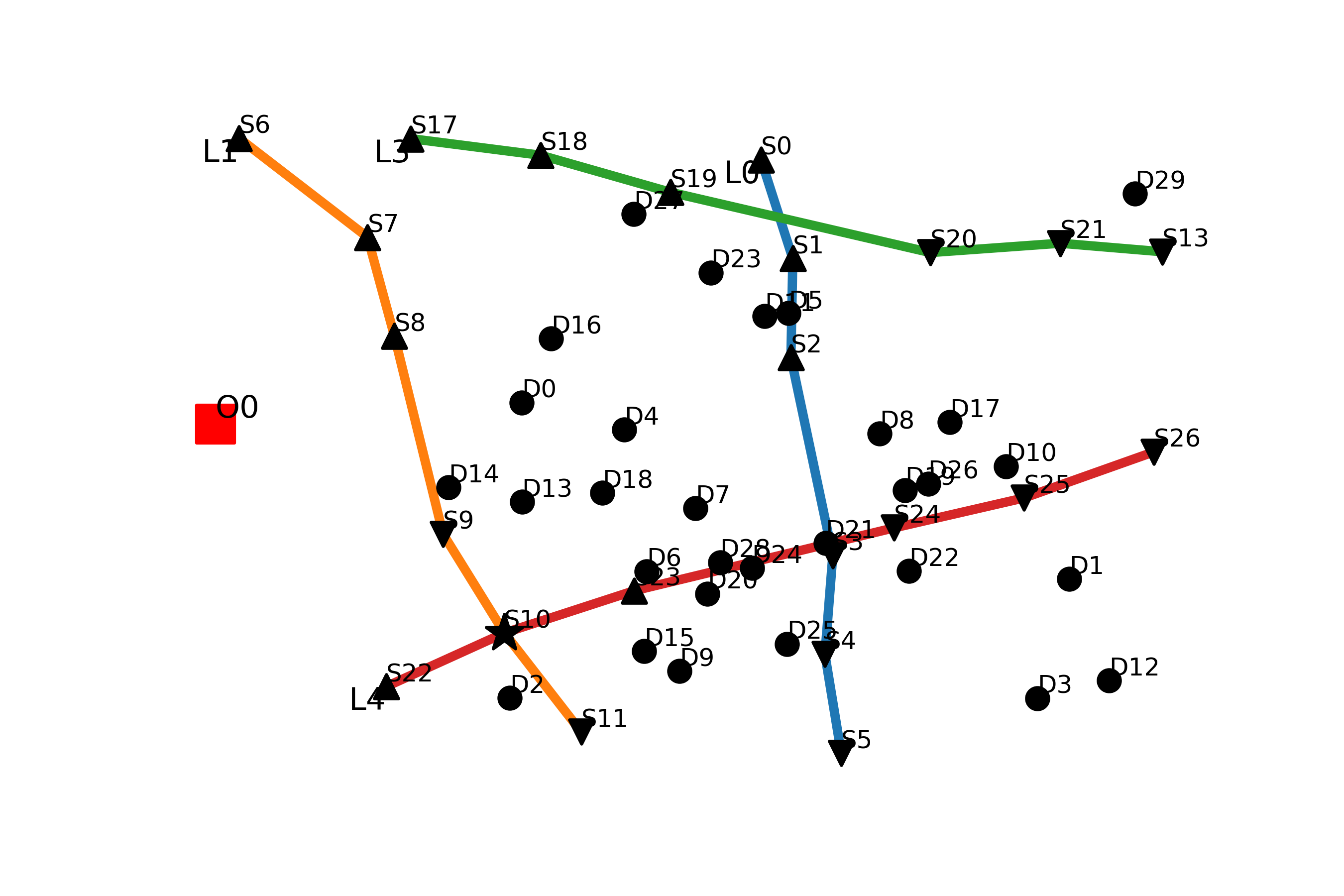}
\vspace{-3em}
\caption{\textit{An example of a synthetic instance}}
\label{instance_example}
\end{figure}

We test our models on artificial instances, generated randomly, to mimic real-life transportation networks. We use transportation networks similar to the ones described in the paper by \citet{donne2021freight}, with some modifications to suit our setting. In their paper, the authors handle a strategic problem associated with the one studied in this paper, so we adapt some of the features for our case. Figure \ref{instance_example} shows a typical example of a synthetic instance. The red square denoted by O0 denotes the CDC. The thick lines, along with the triangles, depict the public transportation network. Each line refers to a single bus or metro line, which consists of several public vehicles throughout the day. We represent the stops on the public vehicle network with triangles; upward triangles denote drop-in stops, while inverted triangles represent drop-out stops. Stars denote stops that serve as both drop-in and drop-out. Finally, the circles represent the customers. The figure shows an instance with four public vehicle lines, twelve drop-in stops, twelve drop-out stops (three drop-in stops and three drop-out stops on each line), and thirty customers.

We divide the whole day into 30 periods, each 30 minutes long. Thus, the entire time horizon during which the delivery system is active is 15 hours in our implementations. We assume that the first period starts with the $0^{th}$ hour of the delivery system. We reserve the first two and a half hours of the system exclusively for the delivery trucks to start transporting packages from the CDC to the drop-in stops. Hence, we do not include buses or metros in our delivery system during this period. The public vehicles that are supposed to carry packages start operating from time point 150 (or after two and a half hours have passed). We use 15 public vehicles per line for up to 50 customers and 18 vehicles per line for instances with a higher number of customers. They are assumed to run every 30 minutes starting from time point 150. We require the definition of periods primarily for \textit{Obj3}, where we estimate the number of freighters required per period.

We generate the customer time windows such that each time window is at least 3 hours long, and it can potentially be as long as the entire delivery horizon. Here, we have not restricted the delivery period to certain hours of the day. If we want to restrict the delivery system to specific periods during the day, like the early hours of the morning or late at night, we need to consider the public vehicles during those hours. We set the value of $t_{mid-day}$ to 400. Finally, we allow a package to stay at any stop for at most 300 time units or 5 hours, and we assume it to be the same for all the stops in our implementations. We have considered the service time at each stop to be 10 minutes. The service time at each of the customers is assumed to be 0. We do so because a positive service time does not change our analyses; it only adds a fixed value to the delivery times.

For the numerical experiments, we use a homogeneous fleet of delivery vehicles. The capacities of the delivery trucks are assumed to be 160 units, and the ones of the freighters are assumed to be 20 units. The capacities of the delivery trucks are significantly higher than that of the freighters because, usually, traditional delivery trucks are quite large in practice. We randomly generate the capacity of public vehicles between 60 and 80 units. However, all the public vehicles that run on the same line are assumed to have the same capacity. This is a very simplistic assumption, but it can be changed easily by manipulating the capacities as a percentage of the total available capacity in public vehicles. These numbers could be higher during the off-peak hours and very small or even zero during the busy hours. These numbers could also come from a distribution that models passenger traffic on public vehicles. The demands of the customers lie between 5 and 20 units. We chose such values of demands and capacities to have some consolidation during the deliveries but to prevent all packages from being assigned to the same vehicles at any tier.

We measure the distance between any two nodes (the CDC and the drop-in stops, any two stops, the drop-out stops and the customers, and any two customers) by the Euclidean distance between them. For our tests, we assume that the speed of each vehicle in the system is the same, or in other words, the time taken by any transportation mode, be it delivery trucks, public vehicles, or freighters, is the same for traversing between any two locations. Even though this is a crude assumption, we can modify it by simply changing the data without any change in the formulations. The time taken for delivery and the corresponding cost is assumed to be proportional to the distance. To convert the distance into the time taken for traversal between any two points in the network, we multiply the distance by 0.2. 

The cost of using the delivery trucks in tier 1 is assumed to be equal to the distance traversed, i.e., $C^{1}_{uvd} = D_{uv}, u,v \in \mathcal{S}_{in} \cup \{o\}, d \in \mathcal{D}$. For the third tier, we assume that the cost of using freighters is 50\% of the costs of using trucks, i.e., $C^{3}_{ijk} = 0.5*D_{ij}, i, j \in \mathcal{S}_{out} \cup \mathcal{C}, k \in \mathcal{K}$. We use these values because we intend to use vehicles of lower costs, both economic and environmental, here.

We created 24 instances, each with a different number of customers, lines, drop-in stops, and drop-out stops. Among these, we have three instances with 10 customers, three with 20 customers, and up to three instances with 80 customers. We did not go beyond 80 customers because not all solution methods could solve the instances beyond this size. The number of public vehicle lines ranges between 1 and 6. The number of drop-in stops varies between 2 and 24, and the number of drop-out stops varies between 2 and 20. Table \ref{tab:instance_description} in the Appendix \ref{instance_generation_appendix} summarizes the description of the instances. Often, the total number of stops is not equal to the number of drop-in stops plus the number of drop-out stops because a stop can serve as both a drop-in and a drop-out stop. 
We do not define the size of an instance due to the different elements in it. An instance with 20 customers, 4 lines, 10 drop-in stops, and 12 drop-out stops would be very different from an instance with 20 customers, 3 lines, 8 drop-in stops, and 9 drop-out stops. We could consider the instance size to be the total number of nodes in the network. However, the number of individual elements in the instances would cause them to behave differently owing to the structure of the network and the problem itself. Hence, we refrain from formally defining a so-called instance size.

For solving the models, we provide a time limit to CPLEX. We set the time limit to 3 hours for \textit{FULL}. For \textit{D1}, the time limit for solving T1 was 2 hours, and 1 hour each for T2 and T3. Similarly, for \textit{D3}, the time limit for solving T3 was 2 hours and 1 hour for the consecutive tiers. For \textit{D2}, we had a time limit of 1 hour for solving each tier. However, whenever T3 is not the first tier to be solved, since the problem decomposes for each drop-out stop, the 1-hour time limit is for solving the problem for each individual drop-out stop. Thus, the total time limit for solving the third tier is actually higher. Finally, the value of $M$, used to linearize the constraints, was set to 1000.
\FloatBarrier

\subsection{Computational Results} \label{computational_results}

\begin{figure}[ht]
\centering
\includegraphics[width=9cm]{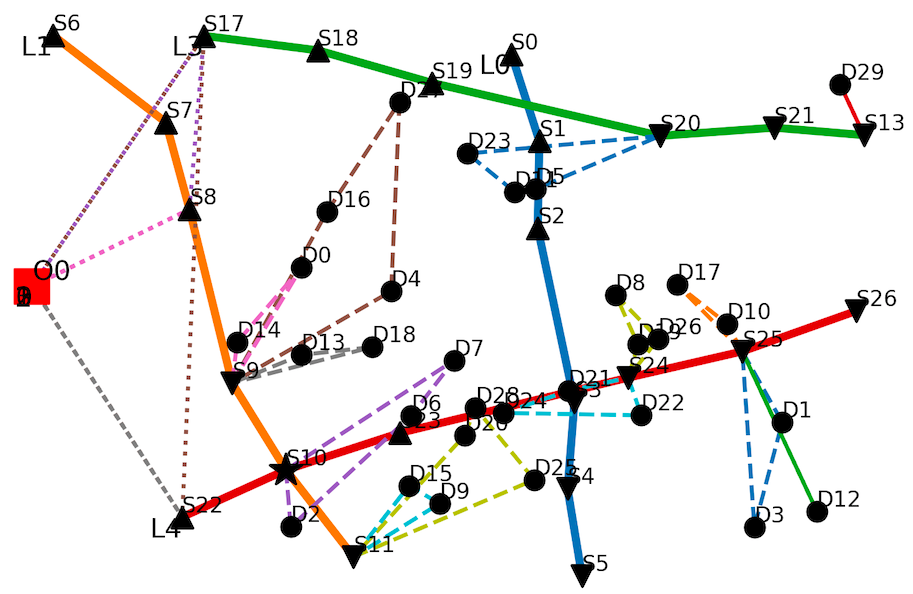}
\vspace{-1em}
\caption{\textit{An example of the solution of an instance}}
\label{fig:a_solution_example}
\end{figure}

We solve our instances using the formulations described by \textit{FULL}, \textit{D1-Obj1}, \textit{D1-Obj2}, \textit{D1-Obj3}, \textit{D2-Obj1}, \textit{D2-Obj2}, \textit{D2-Obj3}, \textit{D3-Obj1}, and \textit{D3-Obj2}. An example of the solution of an instance is given in figure \ref{fig:a_solution_example}. It shows the routes obtained in T1 and T3. The dotted lines show the routes of the delivery trucks, and the dashed lines show the routes of the freighters.

% \FloatBarrier

% \vspace{-1.5em}
\begin{table}[htbp]
  \centering
  \caption{\textit{Solution using formulation \textit{FULL}}}
  \vspace{-1em}
  \resizebox{8cm}{!}{
    \begin{tabular}{ccccc}
    \toprule
    \multicolumn{1}{p{4em}}{ \centering \textbf{Instance number}} & \multicolumn{1}{p{5.5em}}{ \centering \textbf{Routing cost T1}} & \multicolumn{1}{p{5.4em}}{ \centering \textbf{Routing cost T3}} & \multicolumn{1}{p{6.7em}}{ \centering \textbf{Total routing cost}} & \multicolumn{1}{p{4.2em}}{ \centering \textbf{Running time}} \\
    \midrule
    1     & 588.27 & 1706.75 & 2295.02 & 7.94 \\
    2     & 668.86 & 791.38 & 1460.23 & 5.30 \\
    3     & 401.62 & 848.40 & 1250.02 & 85.69 \\
    4     & 767.22 & 1415.48 & 2182.70 & 9751.70 \\
    5     & 1031.46 & 795.40 & 1826.87 & 11153.71 \\
    6     & 428.72 & 1500.46 & 1929.18 & 11256.77 \\
    7     & 1252.22 & 1624.38 & 2876.60 & 14922.37 \\
    8     & 1319.49 & 2083.48 & 3402.97 & 13903.76 \\
    9     & 981.62 & 1649.33 & 2630.94 & 12533.18 \\
    10    & ---    & ---    & ---    & --- \\
    11    & 1068.76 & 2226.12 & 3294.88 & 39396.09 \\
    12    & 1845.23 & 3006.79 & 4852.02 & 23229.73 \\
    \bottomrule
    \end{tabular}%
    }
  \label{tab:solution_full_formulation}%
\end{table}%

Table \ref{tab:solution_full_formulation} provides the solution of the \textit{FULL} formulation of the problem on our instances. It could solve all the instances with 30 customers and only two of the instances with 40 customers. \textit{FULL} could not find a feasible solution for the other instances within a reasonable time. The computational time increases dramatically when moving from 10 to 20 customers and blows up for instances with 40 customers.

In tables \ref{tab:number_of_best_solutions_found} and \ref{tab:avg_percentage_deviation_of_each_method_wrt_best_soln}, we provide a summary of the solutions obtained through each of the solution approaches. Table \ref{tab:number_of_best_solutions_found} reports the number of best solutions found using each decomposition method. Each row and column combination refers to a decomposition method. The number in the cells gives the number of instances for which the technique found the best solution, and the number in brackets provides the number of instances that the method could solve. For \textit{D1}, since tier 1 is the first tier to be solved, we have a unique T1 routing cost. The same happens for the routing costs of T3 for \textit{D3}. For \textit{FULL}, we have one value in each column reporting the routing costs of T1, T2, and the total routing cost. \textit{D1} and \textit{D2} could find a feasible solution for all the 24 instances within the time limit, while \textit{D3} could solve only up to instance 15. In table \ref{tab:avg_percentage_deviation_of_each_method_wrt_best_soln}, we report the average percentage deviation of the solutions found by each decomposition method from the best solution found for that instance, the average taken only over the instances that the approach could solve. 
\begin{table}[htbp]
  \centering
  \caption{\textit{Number of best solutions found by each solution method}}
  \vspace{-1em}
  \resizebox{12cm}{!}{  
  \begin{tabular}{cccccccccc}
    \toprule
          \multicolumn{1}{p{14em}}{\textbf{\textit{Decomposition technique}}} & \multicolumn{3}{c}{\textbf{Routing cost T1}} & \multicolumn{3}{c}{\textbf{Routing cost T3}} & \multicolumn{3}{c}{\textbf{Total routing cost}} \\
    \cmidrule(lr){2-4} \cmidrule(lr){5-7} \cmidrule(lr){8-10}
    \multicolumn{1}{p{14em}}{\textbf{\textit{with Objective function}}} & \textit{Obj1}  & \textit{Obj2}  & \textit{Obj3}  & \textit{Obj1}  & \textit{Obj2}  & \textit{Obj3}  & \textit{Obj1}  & \textit{Obj2}  & \textit{Obj3} \\
    \cmidrule(lr){1-1} \cmidrule(lr){2-4} \cmidrule(lr){5-7} \cmidrule(lr){8-10}
    \textit{FULL}  & \multicolumn{3}{c}{5 (11)} & \multicolumn{3}{c}{3 (11)} & \multicolumn{3}{c}{10 (11)} \\
    \textit{D1}    & \multicolumn{3}{c}{15 (24)} & 0 (24) & 0 (24) & 0 (24) & 0 (24) & 0 (24) & 0 (24) \\
    \textit{D2}    & 1 (24) & 10 (24) & 2 (24) & 1 (24) & 9 (24) & 0 (24) & 0 (24) & 11 (24) & 0 (24) \\
    \textit{D3}    & 0 (15) & 1 (15) & --    & \multicolumn{3}{c}{15 (15)}    & 0 (15) & 5 (15) & -- \\
    \bottomrule
    \end{tabular}%
    }
  \label{tab:number_of_best_solutions_found}%
\end{table}%

% Table generated by Excel2LaTeX from sheet 'Sheet3'
\begin{table}[htbp]
  \centering
  \caption{\textit{Average percentage deviation of the solutions wrt the best solution found for each instance}}
  \vspace{-1em}
  \resizebox{12.9cm}{!}{
  \begin{tabular}{cccccccccc}
    \toprule
          \multicolumn{1}{p{14em}}{\textbf{\textit{Decomposition technique}}} & \multicolumn{3}{c}{\textbf{Routing cost T1}} & \multicolumn{3}{c}{\textbf{Routing cost T3}} & \multicolumn{3}{c}{\textbf{Total routing cost}} \\
    \cmidrule(lr){2-4} \cmidrule(lr){5-7} \cmidrule(lr){8-10}
    \multicolumn{1}{p{14em}}{\textbf{\textit{with Objective function}}} & \textit{Obj1}  & \textit{Obj2}  & \textit{Obj3}  & \textit{Obj1}  & \textit{Obj2}  & \textit{Obj3}  & \textit{Obj1}  & \textit{Obj2}  & \textit{Obj3} \\
    \cmidrule(lr){1-1} \cmidrule(lr){2-4} \cmidrule(lr){5-7} \cmidrule(lr){8-10}
    \textit{FULL}  & \multicolumn{3}{c}{16.77\%} & \multicolumn{3}{c}{4.01\%} & \multicolumn{3}{c}{0.62\%} \\
    \textit{D1}    & \multicolumn{3}{c}{9.09\%} & \multicolumn{1}{c}{50.09\%} & 32.63\% & 54.40\% & 28.64\% & 17.53\% & 31.63\% \\
    \textit{D2}    & 69.75\% & 15.08\% & 87.58\% & \multicolumn{1}{c}{24.76\%} & 12.85\% & 31.50\% & 28.42\% & 6.17\% & 37.46\% \\
    \textit{D3}     & 71.52\% & 49.74\% & --    & \multicolumn{2}{c}{0.00\%} & --    & 14.15\% & 8.35\% & -- \\
    \bottomrule
    \end{tabular}%
    }
  \label{tab:avg_percentage_deviation_of_each_method_wrt_best_soln}%
\end{table}%

The average running time using each method is provided in table \ref{tab:average_running_time}, with the average, again, considered only over the instances that the method could find a solution for within the time limit. The running time reported is the sum of the running times of all the tiers for any approach. Also note that, often, particularly for larger instances, the solution time taken is far greater than the time limit. This is because of the presolve time taken by CPLEX before actually solving the problem.

% \vspace{-0.5em}
% Table generated by Excel2LaTeX from sheet 'Sheet3'
\begin{table}[htbp]
  \centering
  \caption{\textit{Average computational time in seconds}}
  \vspace{-1em}
  \resizebox{7cm}{!}{
    \begin{tabular}{cccc}
    \toprule
    \multicolumn{1}{p{14em}}{\textbf{\textit{Decomposition technique with Objective function}}} & \textit{Obj1}  & \textit{Obj2}  & \textit{Obj3} \\
    \midrule
    \textit{FULL}  & \multicolumn{3}{c}{12386.02} \\
    \textit{D1}    & 5140.18 & 5839.04 & 6272.99 \\
    \textit{D2}    & 5597.89 & 3827.48 & 5429.45 \\
    \textit{D3}    & 20242.13 & 19450.94 & -- \\
    \bottomrule
    \end{tabular}%
    }
  \label{tab:average_running_time}%
\end{table}%

\begin{figure}[ht]
\centering
\includegraphics[width=10cm, keepaspectratio]{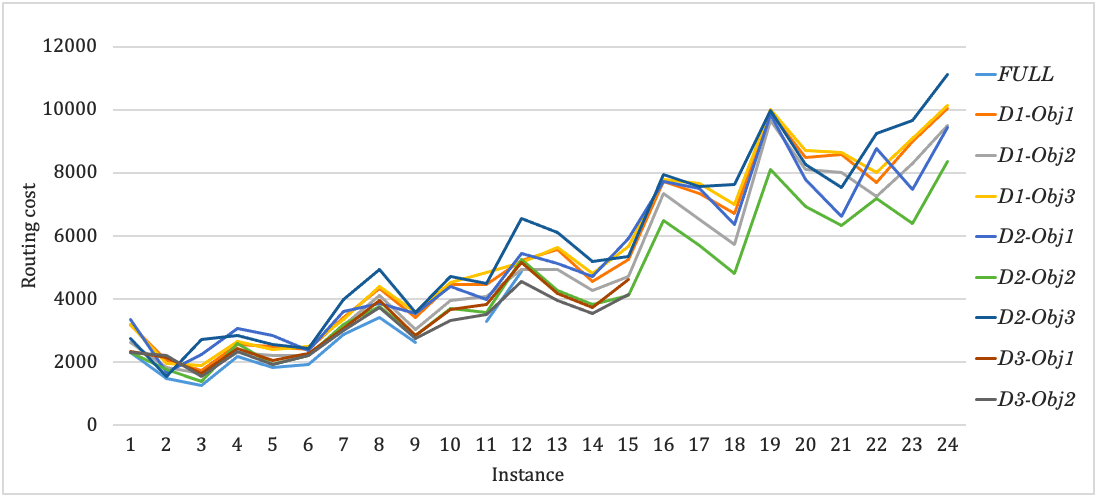}
\vspace{-0.5em}
\caption{\textit{Comparison of the total routing cost obtained by each of the solution approaches}}
\label{fig:solution_comparison_total_routing_cost}
\end{figure}

\begin{figure}[ht]
    \centering
    \begin{subfigure}{.5\textwidth}
      \centering
      \includegraphics[width=1\linewidth]{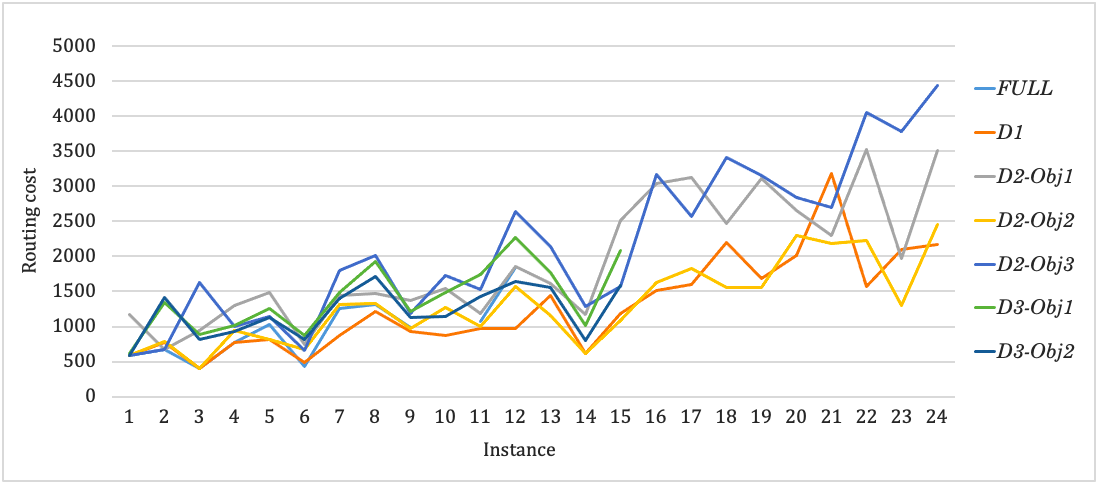}
      \caption{\textit{T1 routing costs}}
      \label{fig:solution_comparison_T1_routing_cost}
    \end{subfigure}%
    \begin{subfigure}{.5\textwidth}
      \centering
      \includegraphics[width=1\linewidth]{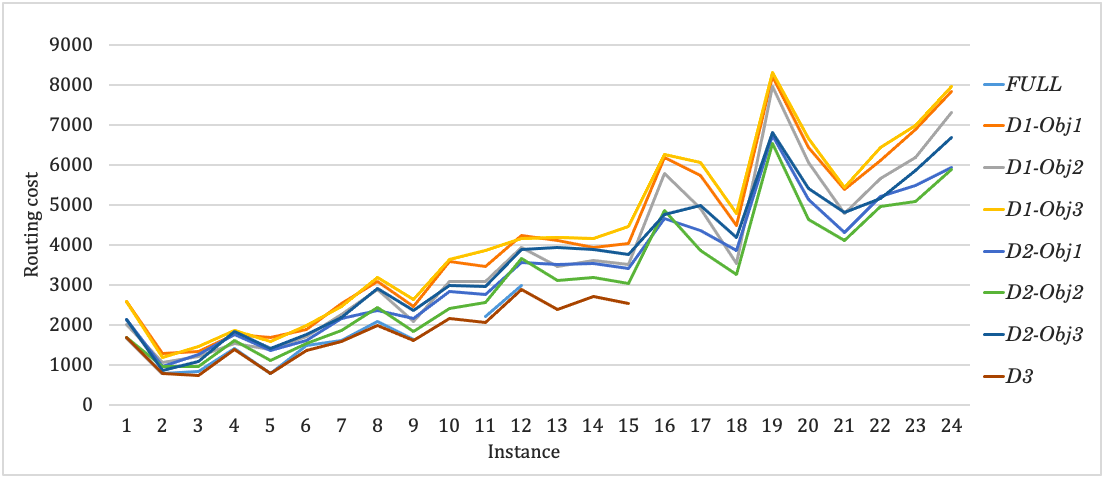}
      \caption{\textit{T3 routing costs}}
      \label{fig:solution_comparison_T3_routing_cost}
    \end{subfigure}
    \vspace{-0.5em}
    \caption{\textit{Comparison of the routing cost of T1 and T3 obtained by each of the solution approaches}}
    \label{fig:comparison_T1_T3_routing costs}
\end{figure}

The performance of the approaches is depicted in figures \ref{fig:solution_comparison_total_routing_cost} and \ref{fig:comparison_T1_T3_routing costs}. Each line in the graphs shows a solution approach-- either the full formulation or a decomposition technique with one of the three objectives for T2. The x-axis represents the instances, and the y-axis gives the value of the solutions, i.e., the routing costs. When we solve \textit{D1}, we have only one set of solutions for tier 1, regardless of the objective function used in T2, since T1 is solved first. Similarly, for \textit{D3}, we have one set of solutions for the routing costs of tier 3, irrespective of whether \textit{Obj1} or \textit{Obj2} is used later in T2. This explains why we have only one line showing the routing costs of T1 (\textit{D1}) and T3 (\textit{D3}) in figures \ref{fig:solution_comparison_T1_routing_cost} and \ref{fig:solution_comparison_T3_routing_cost}, respectively. Below, we discuss the performance of the approaches in detail.

Among all the approaches, the best-performing objective function for T2 is \textit{Obj2}, i.e., to minimize the approximated routing costs of T1 and T3. This is related to the fact that our true objective is to minimize the routing costs. \textit{Obj1} and \textit{Obj3} are able to find the best solution in only a handful of cases. The approaches that use \textit{Obj3}, which is to minimize the approximate number of freighters in T3, generally perform the worst within each decomposition method, especially as the instances grow larger. This is also reflected in table \ref{tab:number_of_best_solutions_found}, where we see that the average percentage deviations from the best solution found are the highest for \textit{Obj3}.

The full formulation beats the decomposition approaches in terms of the total routing costs, but it can only solve very small instances, and even then, it takes a long time. Between the three decomposition approaches, \textit{D3-Obj2} performs the best initially. However, compared to the other decomposition techniques, it can solve relatively smaller-sized instances-- only instances with up to 50 customers. Moreover, the running times of \textit{D3} also escalate quickly after instances with 20 customers. For larger instances, almost all of the best total routing costs were found by \textit{D2-Obj2}.

For the routing costs of T1, \textit{D1} performs the best initially, but, later on, for larger-sized instances, it is often beaten by \textit{D2-Obj2}. For T3 routing costs, \textit{D3} finds the best solutions for instances that it can solve. After that, \textit{D2-Obj2} obtains the best T3 routing costs as well. Thus, overall, \textit{D2-Obj2} dominates the other decomposition techniques as the instance size increases.

These findings pertain partially to the particularity of our setting and implementations. We have considered one CDC, up to 80 customers, equal delivery speeds for trucks and freighters, and the majority of the distance covered comes from T3. Moreover, the cost of covering a certain distance in T3 is assumed to be 50\% of the cost in T1 using trucks. If we have multiple CDCs, and the cost of using delivery trucks is even higher, or our objective is to minimize emissions in particular, and low-cost green means are used in T3, it would be beneficial to use \textit{D1}. If the cost of vehicles used in T3 is significantly higher and greatly depends on the individual trips performed, we might want to start the decomposition technique from T3. Finally, if the setting is such that the routing costs of T1 and T3 are comparable, we would be better off using \textit{D2}. For the objective functions, we can conclusively say that \textit{Obj2} is the best performing one because our actual goal is to minimize routing costs. On average, we find the decomposition approach \textit{D2-Obj2} to perform the best, both in terms of the objective value of the solutions and the running times.

%%%%%%%%%%%%%%%%%%%%%%%%%%%%%%%%%%%%%%%%%%%%%
%%%%%%%%%%%%%%%%%%%%%%%%%%%%%%%%%%%%%%%%%%%%%
% \FloatBarrier

\begin{figure}
    \centering
    \begin{subfigure}{.5\textwidth}
      \centering
      \includegraphics[width=1\linewidth]{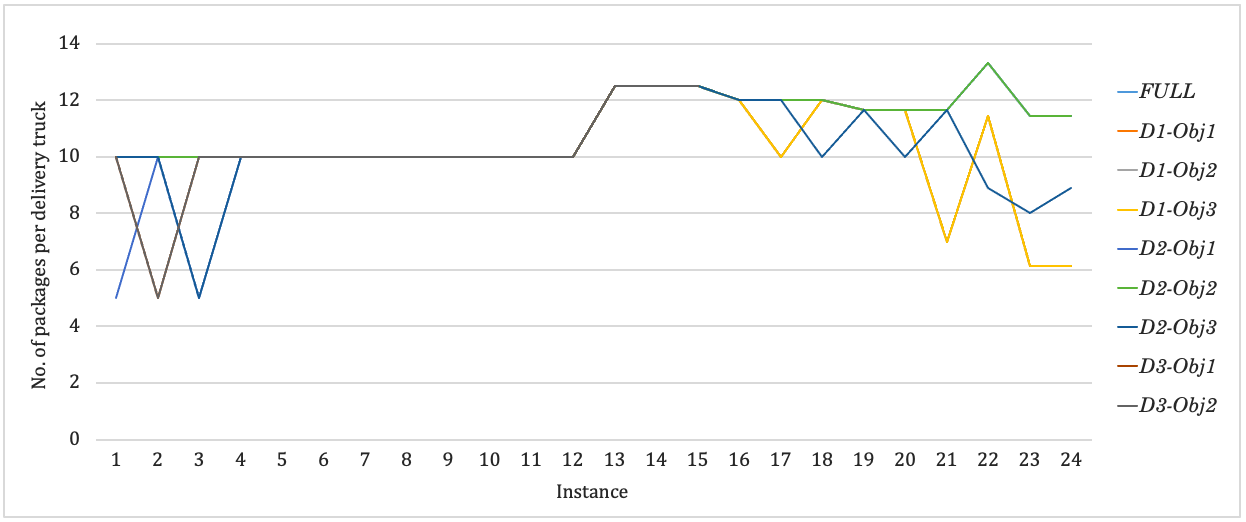}
      \caption{\textit{Average number of packages per truck}}
      \label{fig:customers_per_truck}
    \end{subfigure}%
    \begin{subfigure}{.5\textwidth}
      \centering
      \includegraphics[width=1\linewidth]{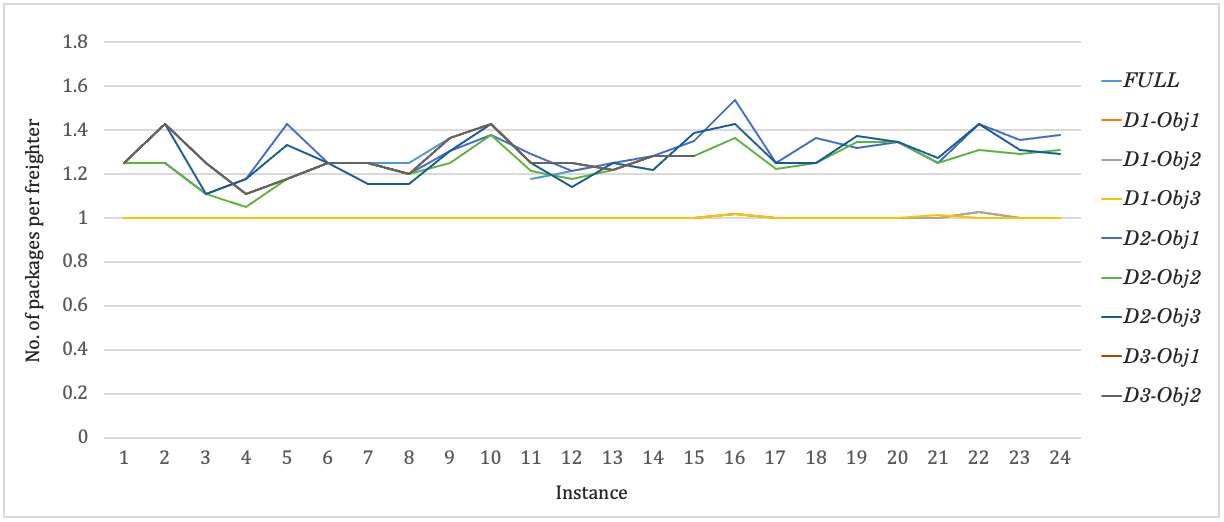}
      \caption{\textit{Average number of packages per freighter}}
      \label{fig:customers_per_freighter}
    \end{subfigure}
    \vspace{-0.5em}
    \caption{\textit{Comparing the usage of delivery trucks and freighters by the different solution approaches}}
    \label{fig:comparison_trucks_freighters_required}
\end{figure}

Next, we discuss some features of the solutions that we obtain. Figure \ref{fig:comparison_trucks_freighters_required} shows the average number of packages assigned to the delivery trucks and freighters for each instance for each solution approach. There does not seem to be a significant difference between the different decomposition techniques in the utilization of delivery trucks (figure \ref{fig:customers_per_truck}). This is because we have one CDC, and the number of trucks used is far fewer than the number of freighters to demonstrate any significant difference.

It is interesting to note that \textit{D1} uses a higher number of freighters than any other solution method. This means that there is little to no consolidation of packages in \textit{D1}, regardless of the objective function used in T2, with the ratio of packages to freighters rarely exceeding 1 (figure \ref{fig:customers_per_freighter}). Though \textit{D1} would be beneficial if our objective is to minimize the emissions caused by large delivery vehicles, it would not be ideal if the routing costs in T3 are high as well. In that case, we would need to make more trips and recruit a higher number of freighters. However, we can use it in some cases, for example, if we use drones to make deliveries to the end customers, and our sole aim is to control the use of trucks.

\begin{figure}
    \centering
    \begin{subfigure}{.5\textwidth}
      \centering
      \includegraphics[width=1\linewidth]{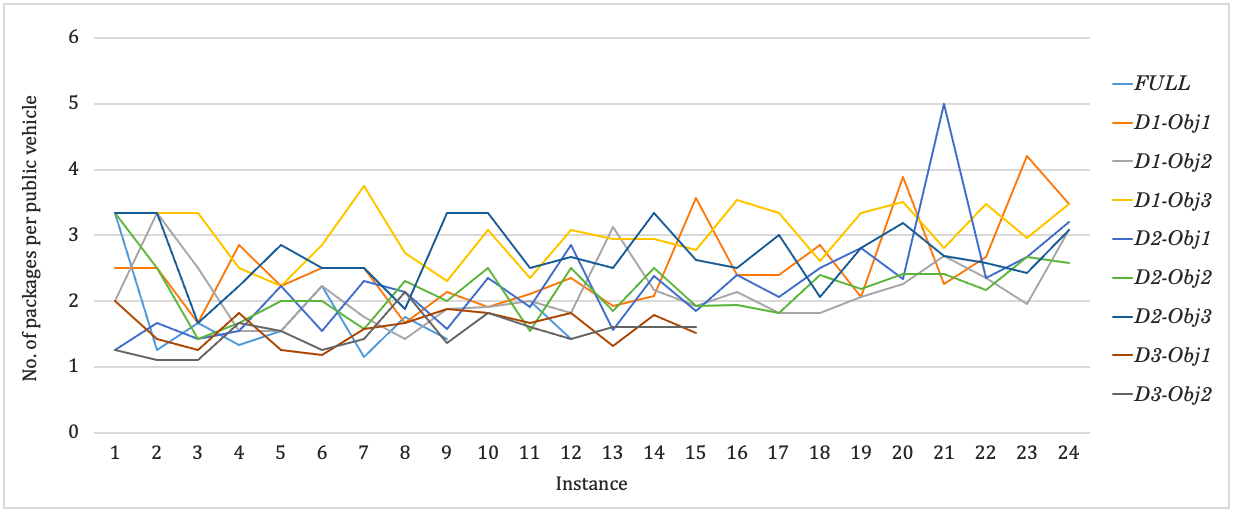}
      \caption{\textit{Average number of packages per public vehicle}}
      \label{fig:customers_per_public_vehicles}
    \end{subfigure}%
    \begin{subfigure}{.5\textwidth}
      \centering
      \includegraphics[width=1\linewidth]{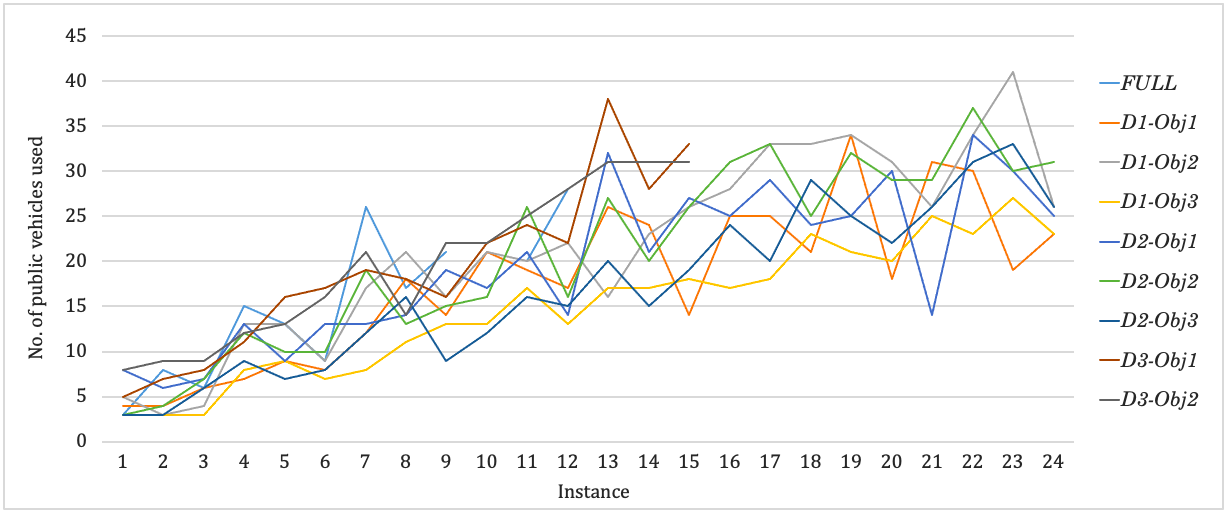}
      \caption{\textit{Number of public vehicles used}}
      \label{fig:percentage_public_vehicles_used}
    \end{subfigure}
    \vspace{-1em}
    \caption{\textit{Comparing the usage of public vehicles by the different solution approaches}}
    \label{fig:comparison_public_vehicles_required}
\end{figure}

We also find from our solutions that decomposition approaches \textit{D3} and \textit{D2-Obj2} consistently use a higher number of drop-out stops and, thus, have a lower ratio of customers per drop-out stop. This also results in smaller distances traveled in T3. \textit{D2-Obj1} utilizes the lowest number of drop-in and drop-out stops among all the approaches, particularly as the instance size increases. Figure \ref{fig:comparison_public_vehicles_required} shows the average number of customers per public vehicle and the number of public vehicles used by each solution technique. On average, \textit{Obj2} utilizes the highest number of public vehicles, and \textit{Obj3} the lowest. Thus, \textit{Obj3} has a higher percentage of packages per vehicle, as the consolidation of packages in T3 is encouraged here.

%%%%%%%%%%%%%%%%%%%%%%%%%%%%%%%%%%%%%%%%%%%%%
%%%%%%%%%%%%%%%%%%%%%%%%%%%%%%%%%%%%%%%%%%%%%

%%%%%%%%%%%%%%%%%%%%%%%%%%%%%%%%%%%%%%%%%%%%%
%%%%%%%%%%%%%%%%%%%%%%%%%%%%%%%%%%%%%%%%%%%%%
% \FloatBarrier

\FloatBarrier

\subsection{Comparison with Traditional Methods}

In order to compare the effectiveness of using public vehicles, we compare the system with a standard delivery method, where trucks are sent directly from the CDC to the customer locations. This corresponds to a Vehicle Routing Problem with Time Windows (VRPTW, henceforth). We use a standard formulation of the VRPTW and present the model in appendix \ref{Appendix_VRPTW}.

Figure \ref{fig:comparison_VRPTW_heterogeneous_0.5} shows the comparison of the routing costs obtained by the simple VRPTW compared to the different approaches for the 3T-DPPT. Figure \ref{fig:comparison_VRPTW_total_routing_costs_heterogeneous_0.5} gives the total routing costs, and \ref{fig:comparison_VRPTW_T1_routing_costs_heterogeneous_0.5} gives the costs obtained using dedicated delivery vehicles; thus, we consider only the cost associated with tier 1 for 3T-DPPT in this case. Even though initially, the routing costs associated with VRPTW are lower, it is easy to see that the costs explode as the number of customers increases. In our approaches, since the majority of the distances are being covered on public vehicles or sustainable freight systems, we keep a check on the usage of delivery trucks. The difference is even more pronounced in figure \ref{fig:comparison_VRPTW_T1_routing_costs_heterogeneous_0.5}. However, one must remember that these differences depend highly on the instances. The difference will be higher if the longest part of the delivery occurs on board using public vehicles. If the CDC is located close to the drop-in stops, and these drop-in stops are also very few, compared to the number of customers and the distances between them and the drop-out stops, then also it would be highly beneficial to use the public transportation system.

%%%%%%%%%%%%%%%%%%%%%%%%%%%%%%%%%%%%%%%%%%%%%
\begin{figure}[ht!]
    \centering
    \begin{subfigure}{.5\textwidth}
      \centering
      \includegraphics[width=1\linewidth]{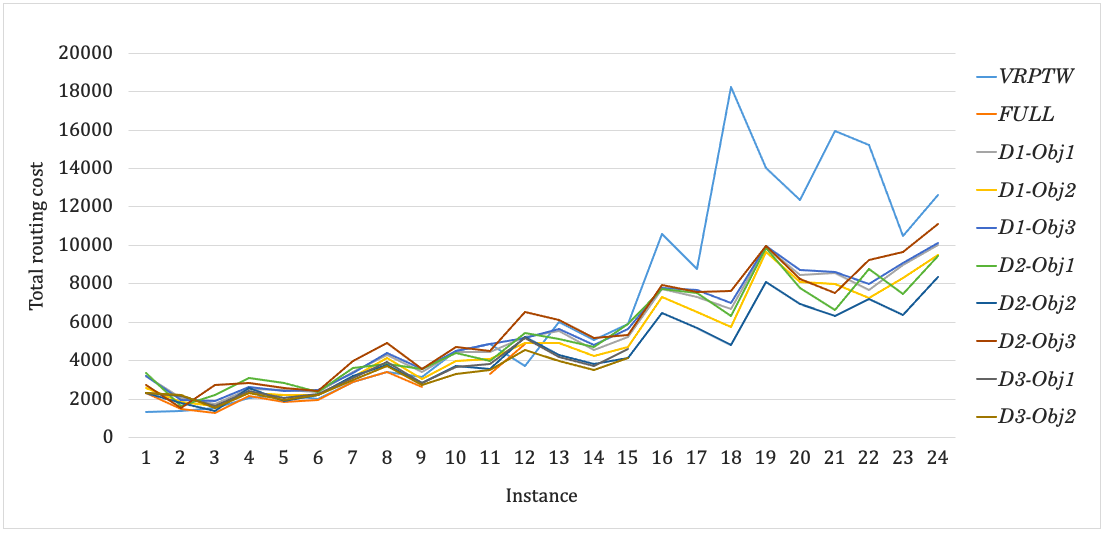}
      \caption{\textit{Total routing costs}}
      \label{fig:comparison_VRPTW_total_routing_costs_heterogeneous_0.5}
    \end{subfigure}%
    \begin{subfigure}{.5\textwidth}
      \centering
      \includegraphics[width=1\linewidth]{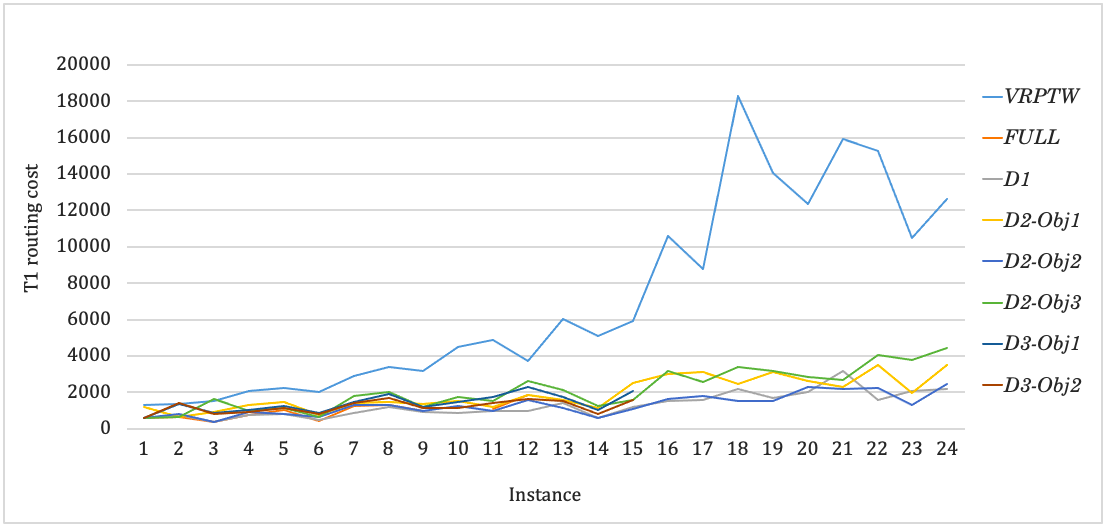}
      \caption{\textit{Routing costs for delivery trucks}}
      \label{fig:comparison_VRPTW_T1_routing_costs_heterogeneous_0.5}
    \end{subfigure}
    \vspace{-0.5em}
    \caption{\textit{Comparison of a  VRPTW with the different solution approaches for 3T-DPPT}}
    \label{fig:comparison_VRPTW_heterogeneous_0.5}
\end{figure}
%%%%%%%%%%%%%%%%%%%%%%%%%%%%%%%%%%%%%%%%%%%%%

Our results show that using public vehicles could lead to a 67.64\% reduction in the distances traversed using heavy trucks on average, consequently reducing the emissions caused by these vehicles proportionally. Even though the service may not be profitable for a small number of customers, we find significant cost reductions as the size of the instances increases. For instances with 40 customers or more, we find a 21.24\% reduction in total costs on average. One should note here that, for our instances, the average distances from the CDC to the customers are comparable to the average distances from the CDC to the drop-in stops plus the distances from the drop-out stops to the customers. This is because we implement our models on instances with relatively fewer customers, representing much smaller-sized towns than in reality. Even then, we demonstrate significant cost savings. In larger cities, we intend to use drop-in stops located close to the CDC, which, in turn, is usually built on the outskirts of the city; and the drop-out stops would be chosen in the vicinity of the customers. Thus, the CDC would, in practice, serve customers much farther away from it, such that the above distances would not always be comparable. Therefore, the leg of delivery using public transportation would be much longer than in our instances, and we hypothesize that it would lead to higher cost savings in the day-to-day implementation of the service.

\FloatBarrier

\subsection{Managerial Insights} \label{managerial_insights}

\subsubsection{Impact of Heterogeneity in Delivery Costs}

In this section, we analyze the behavior of the solutions on small-sized instances--instances with up to 30 customers, for different proportions of costs in T3 compared to costs in T1. We consider $C^{1}_{uvd} = D_{uv}, u,v \in \mathcal{S}_{in} \cup \{o\}, d \in \mathcal{D}$, and $C^{3}_{ijk} = \beta \cdot D_{ij}, i, j \in \mathcal{S}_{out} \cup \mathcal{C}$, where $\beta$ is a parameter that determines the proportion of costs between T1 and T3. In our computations, we have $\beta = 0.1, 0.25, 0.5, 0.75,$ and $1$. We have solved the formulations \textit{D1-Obj1}, \textit{D1-Obj2}, \textit{D1-Obj3}, \textit{D2-Obj1}, \textit{D2-Obj2}, \textit{D2-Obj3}, \textit{D3-Obj1}, and \textit{D3-Obj2} on instances with up to 30 customers to study how the solutions, particularly the usage of delivery trucks, vary.

%%%%%%%%%%%%%%%%%%%%%%%%%%%%%%%%%%%%%%%%%%%%%%%%%%%%%%%%%%%%%
\begin{figure}[ht]
\centering
\includegraphics[width=12cm]{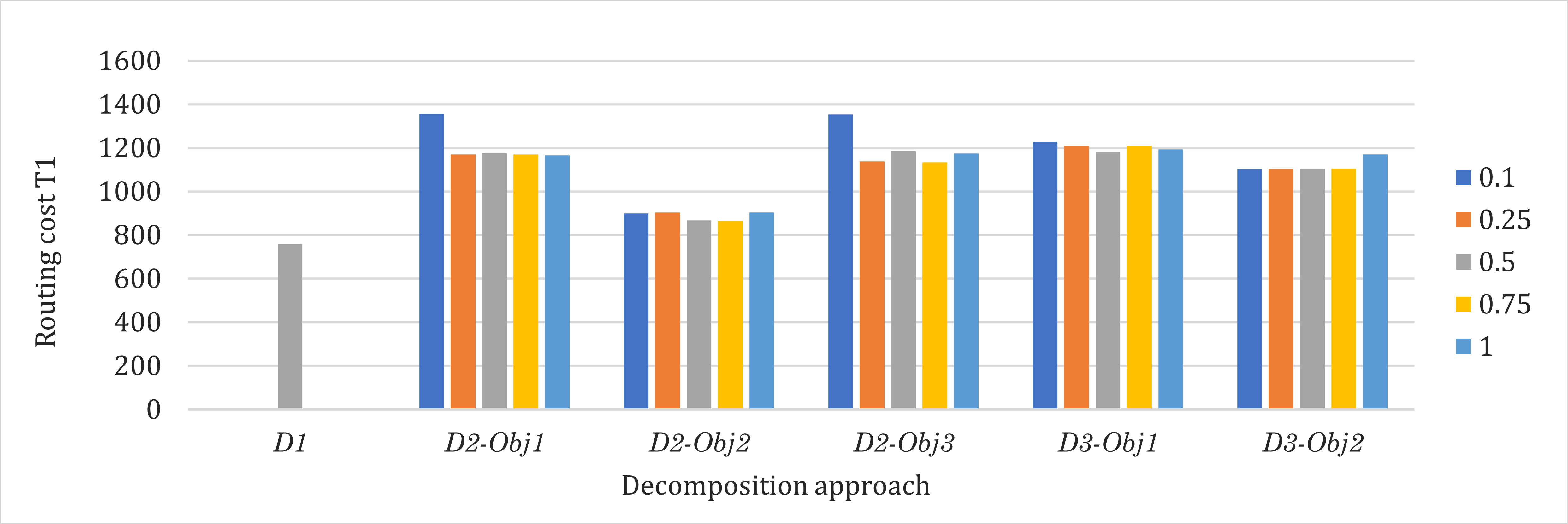}
\vspace{-0.5em}
\caption{\textit{Comparison of T1 routing costs for different values of $\beta$}}
\label{fig:beta_comparison_T1_routing_cost}
\end{figure}
%%%%%%%%%%%%%%%%%%%%%%%%%%%%%%%%%%%%%%%%%%%%%%%%%%%%%%%%%%%%%

Figure \ref{fig:beta_comparison_T1_routing_cost} shows how the routing costs in T1 change when we use different values of $\beta$. For \textit{D1}, we have only one column representing the cost because: firstly, T1 is the first tier to be solved here, and the routing costs obtained here do not depend on which objective function is used in T2; and secondly, the value of $\beta$ does not change T1 routing costs in \textit{D1}, they impact only the costs obtained in T3. Apart from $\textit{D1}$, which is expected to perform well for tier 1, \textit{D2-Obj2} provides quite competitive solutions.

\begin{figure}[ht]
\centering
\includegraphics[width=12cm]{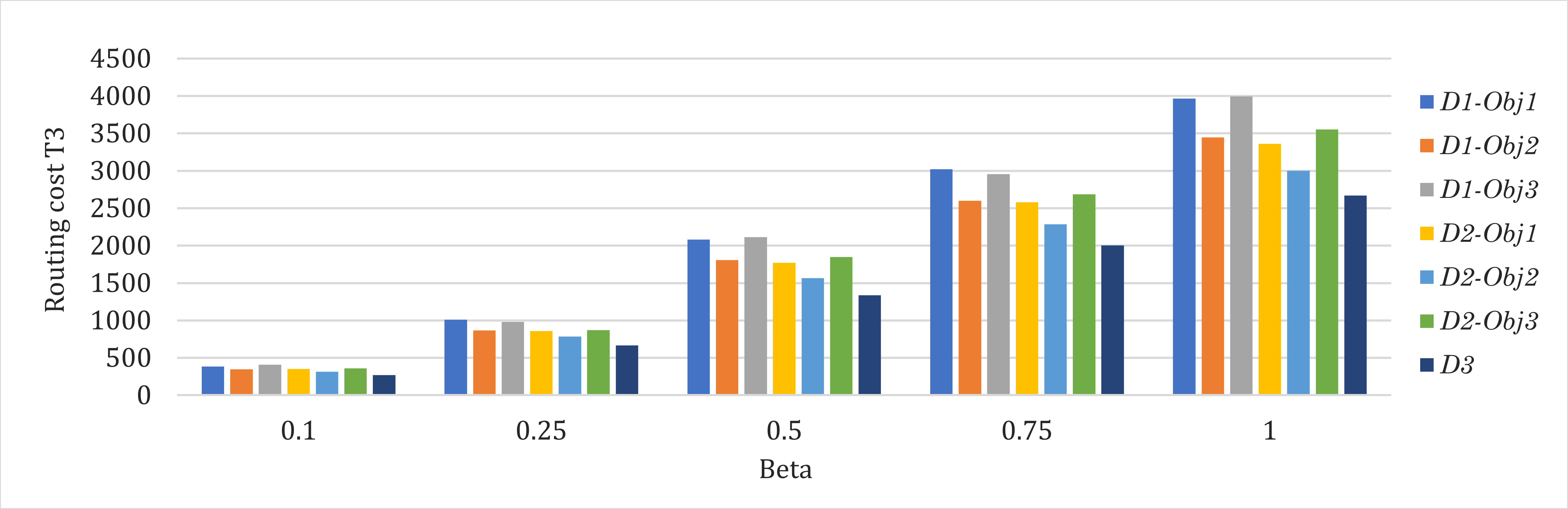}
\vspace{-0.5em}
\caption{\textit{Comparison of T3 routing costs for different values of $\beta$}}
\label{fig:beta_comparison_T3_routing_cost}
\end{figure}

\begin{figure}[ht]
\centering
\includegraphics[width=12cm,keepaspectratio]{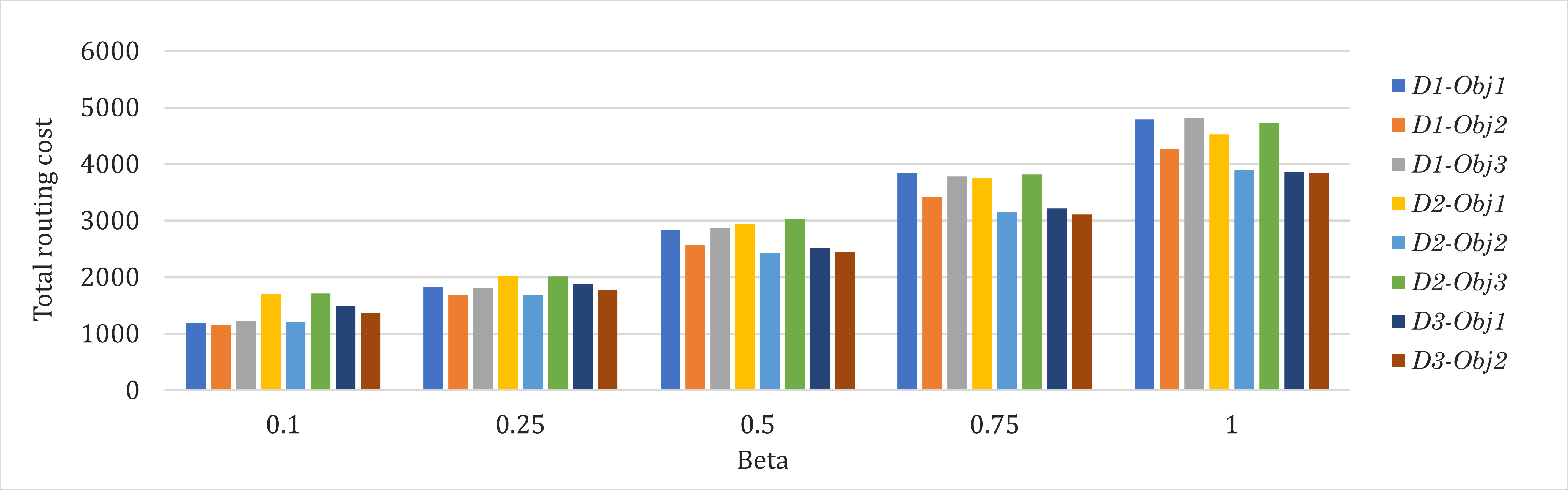}
\vspace{-0.5em}
\caption{\textit{Comparison of total routing costs for different values of $\beta$}}
\label{fig:beta_comparison_total_routing_cost}
\end{figure}

In figures \ref{fig:beta_comparison_T3_routing_cost} and \ref{fig:beta_comparison_total_routing_cost}, we show how each decomposition method performs in terms of T3 routing costs and total routing costs, respectively, for different values of $\beta$. Once again, for T3 routing costs, we have a single column, named D3, for each value of $\beta$. This is because when T3 is solved first, we have not yet used the objective functions of T2, so we obtain the same solutions. For tier 3 routing costs, \textit{D3} performs the best since our instance sizes are small, and \textit{D2-Obj2} follows closely behind. \textit{D2-Obj2} performs the best among the decomposition techniques for the total routing costs as well. Thus, it emerges as a clear winner overall and provides solutions comparable with the best solutions obtained, regardless of the tiers or the values of $\beta$. This is because tier 2 is the linking tier, which balances the objectives of T1 and T3, even though we use approximations, thereby producing good quality solutions. On the other hand, when we prioritize either T1 or T3, the other tier is notably penalized, which increases the overall objective values. Other than the routing costs, there is no significant difference in the features of the solutions obtained compared to section \ref{computational_results}.

%%%%%%%%%%%%%%%%%%%%%%%%%%%%%%%%%%%%%%%%%%%%%%%%%%%%%%%%%%%%%
%%%%%%%%%%%%%%%%%%%%%%%%%%%%%%%%%%%%%%%%%%%%%%%%%%%%%%%%%%%%%

%%%%%%%%%%%%%%%%%%%%%%%%%%%%%%%%%%%%%%%%%%%%%%%%%%%%%%%%%%%%%
%%%%%%%%%%%%%%%%%%%%%%%%%%%%%%%%%%%%%%%%%%%%%%%%%%%%%%%%%%%%%

\FloatBarrier

\subsubsection{Service Costs Related to Tier 2}

So far in our experiments, we have assumed zero delivery costs related to the second tier. In this section, we discuss the impact of adding service costs related to delivery using public transit lines and drop-in and drop-out stops. These costs arise from the cost of personnel recruited at the stops for loading and unloading packages, daily maintenance of the infrastructure, and the cost of space, among others. To address them, we assume that these costs can be incorporated at the drop-in and drop-out stops, i.e., every time a truck visits a drop-in stop or a freighter leaves a drop-out stop, we add a cost to our objective function.
\begin{allowdisplaybreaks}
\begin{align}
    \textit{Minimize} \quad & \sum_{d \in \mathcal{D}} \; \sum_{u \in \mathcal{S}_{in} \cup \{o \} } \sum\limits_{ \substack{ v \in \mathcal{S}_{in} \\ u \neq v } } \left( C^{1}_{uvd} + \lambda_{1} \right) w_{uvd} \; + \; \sum_{d \in \mathcal{D}} \; \sum_{u \in \mathcal{S}_{in} } \sum\limits_{ v = o' } C^{1}_{uvd} w_{uvd} \nonumber \\
    & \; + \; \sum_{k \in \mathcal{K}} \; \sum_{i \in \mathcal{S}_{out} } \sum\limits_{ \substack{ j \in \mathcal{C} \\ i \neq j } } \left( C^{3}_{ijk} + \lambda_{3} \right) x_{ijk} \; + \; \sum_{k \in \mathcal{K}} \; \sum_{i \in \mathcal{C} } \sum\limits_{ \substack{ j \in \mathcal{S}_{out} \cup \mathcal{C} \\ i \neq j } } C^{3}_{ijk} x_{ijk}. \label{obj:full_model_service_costs}
\end{align}
\end{allowdisplaybreaks}
$\lambda_{1}$ and $\lambda_{3}$ are parameters that represent the service costs at tier 1 and tier 3, respectively. This does not imply that the costs are paid to the freighters or truck drivers exclusively; rather, it is a daily operating cost of delivery on the transit line. All the costs mentioned are integrated directly at the stops whenever there are exchange of packages.

We made computational experiments with the value of $\lambda_{1} = \mu \; *$ \textit{average routing cost of T1} and $\lambda_{3} = \mu \; *$ \textit{average routing cost of T3}, with $\mu = 0, 0.1, 0.25, 0.5, 2.0$, to see the impact of adding the service costs on the structure of the solutions. We solved $FULL$ on instances with ten customers since only those instances could be solved up to optimality. Over the instances described above, we did not find any significant change in the solution structure because the value of parameters already facilitates significant consolidation (we expected to see more consolidation as service costs increase). 

\begin{figure}[ht!]
    \centering
    \begin{subfigure}{.35\textwidth}
      \centering
      \includegraphics[width=1\linewidth]{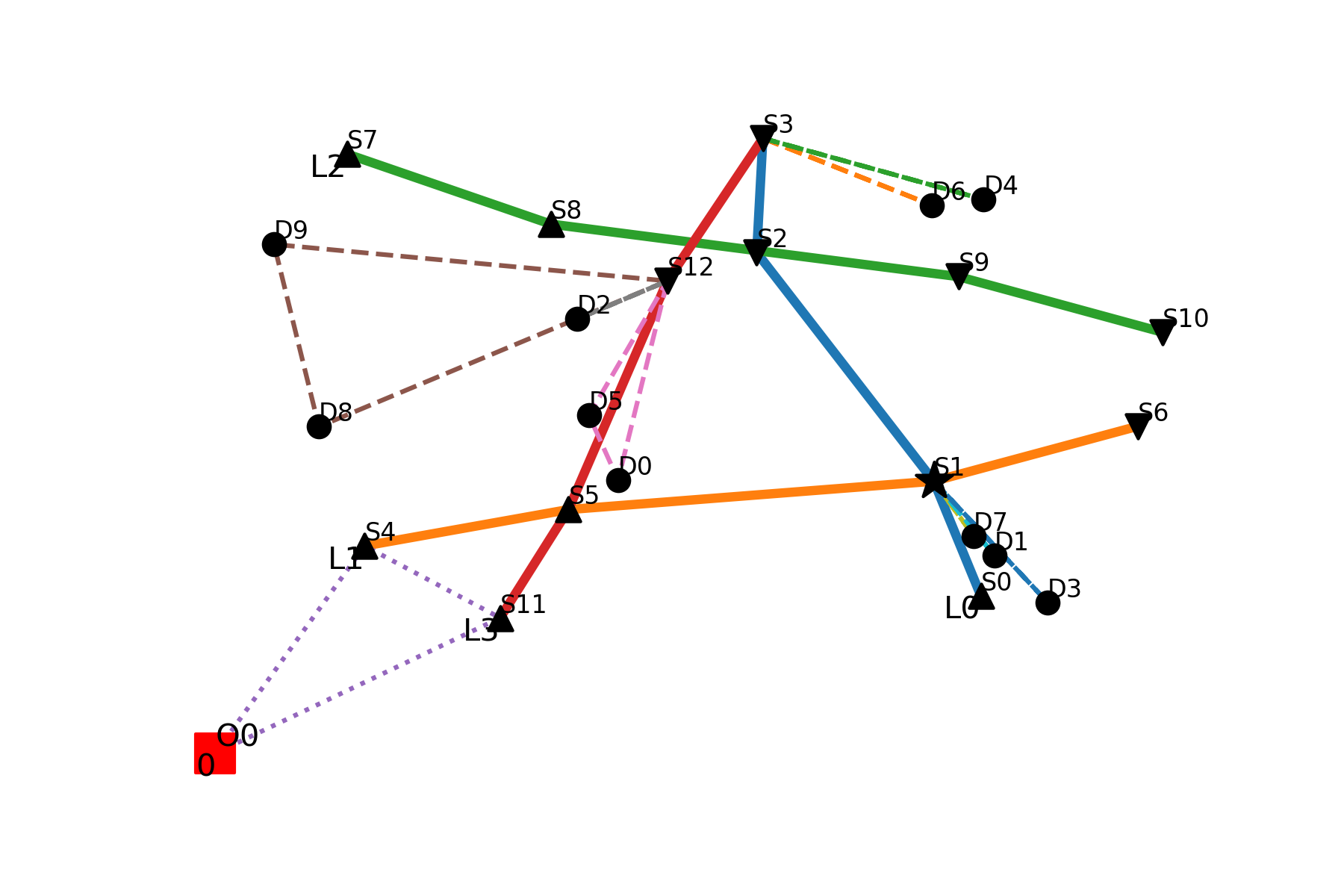}
    \caption{\textit{$\mu = 0.0$, $Q_{F}=20$}}
      \label{fig:SC_1}
    \end{subfigure}%
    %%%%%%%%%%%%%%%%%%%%%%%%%%%%%%%%%%%%
    \begin{subfigure}{.35\textwidth}
      \centering
      \includegraphics[width=1\linewidth]{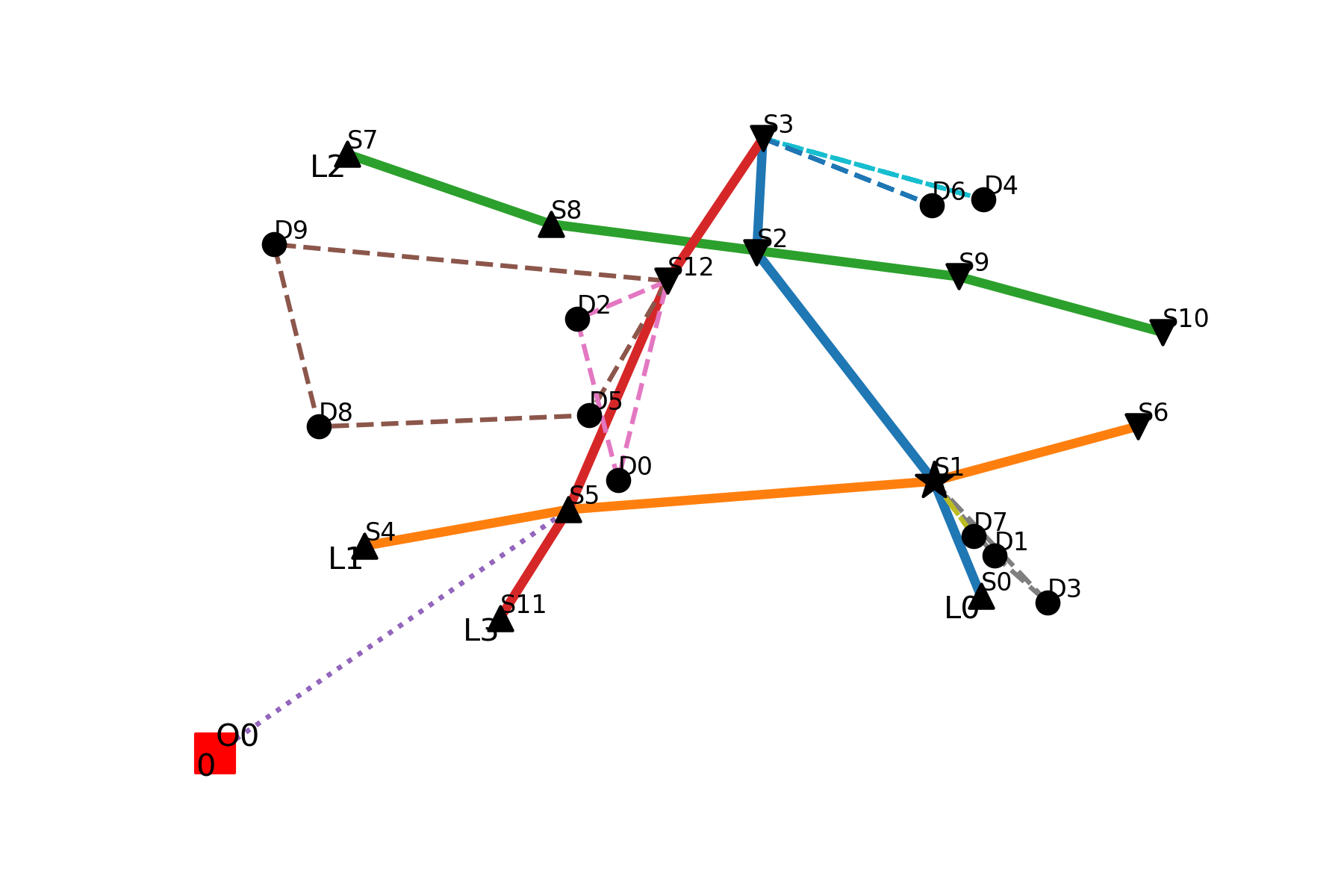}
      \caption{\textit{$\mu = 0.25$, $Q_{F}=30$}}
      \label{fig:SC_2}
    \end{subfigure}%
   %%%%%%%%%%%%%%%%%%%%%%%%%%%%%%%%%%%%
    \begin{subfigure}{.35\textwidth}
      \centering
      \includegraphics[width=1\linewidth]{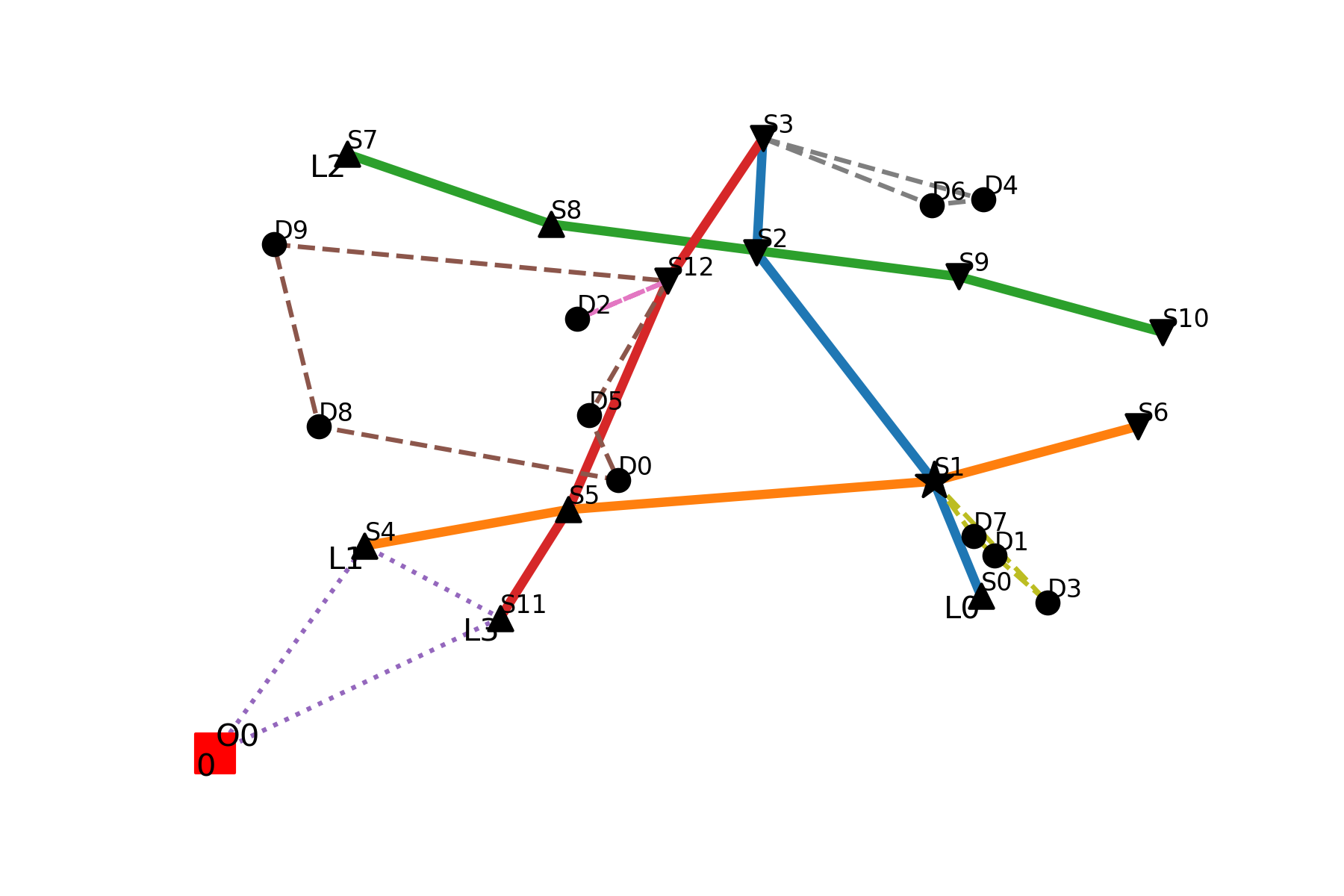}
      \caption{\textit{$\mu = 0.1$, $Q_{F}=50$}}
      \label{fig:SC_3}
    \end{subfigure}
    %%%%%%%%%%%%%%%%%%%%%%%%%%%%%%%%%%%%
    \begin{subfigure}{.35\textwidth}
      \centering
      \includegraphics[width=1\linewidth]{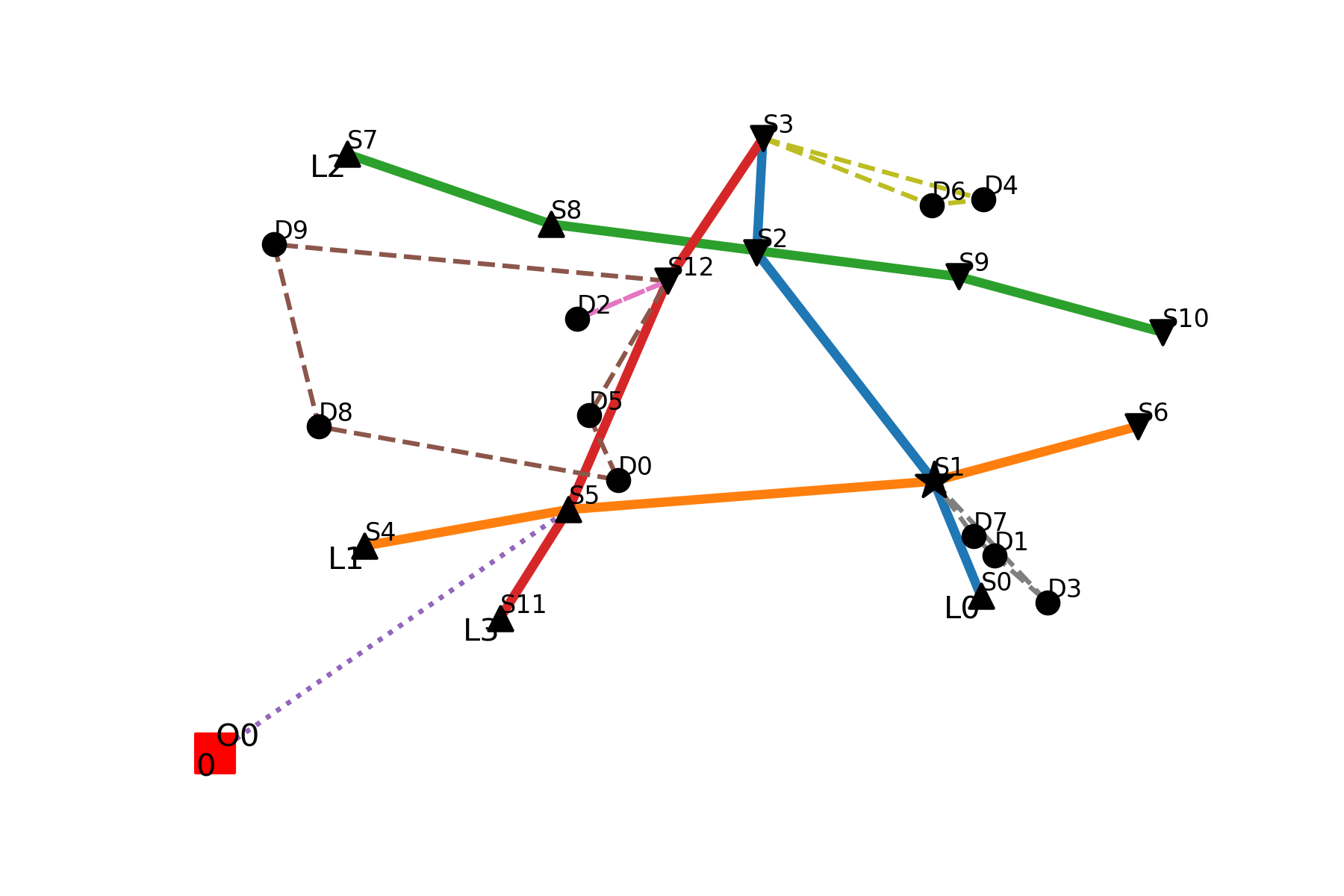}
      \caption{\textit{$\mu = 0.25$, $Q_{F}=50$}}
      \label{fig:SC_4}
    \end{subfigure}%
    %%%%%%%%%%%%%%%%%%%%%%%%%%%%%%%%%%%%
    \begin{subfigure}{.35\textwidth}
      \centering
      \includegraphics[width=1\linewidth]{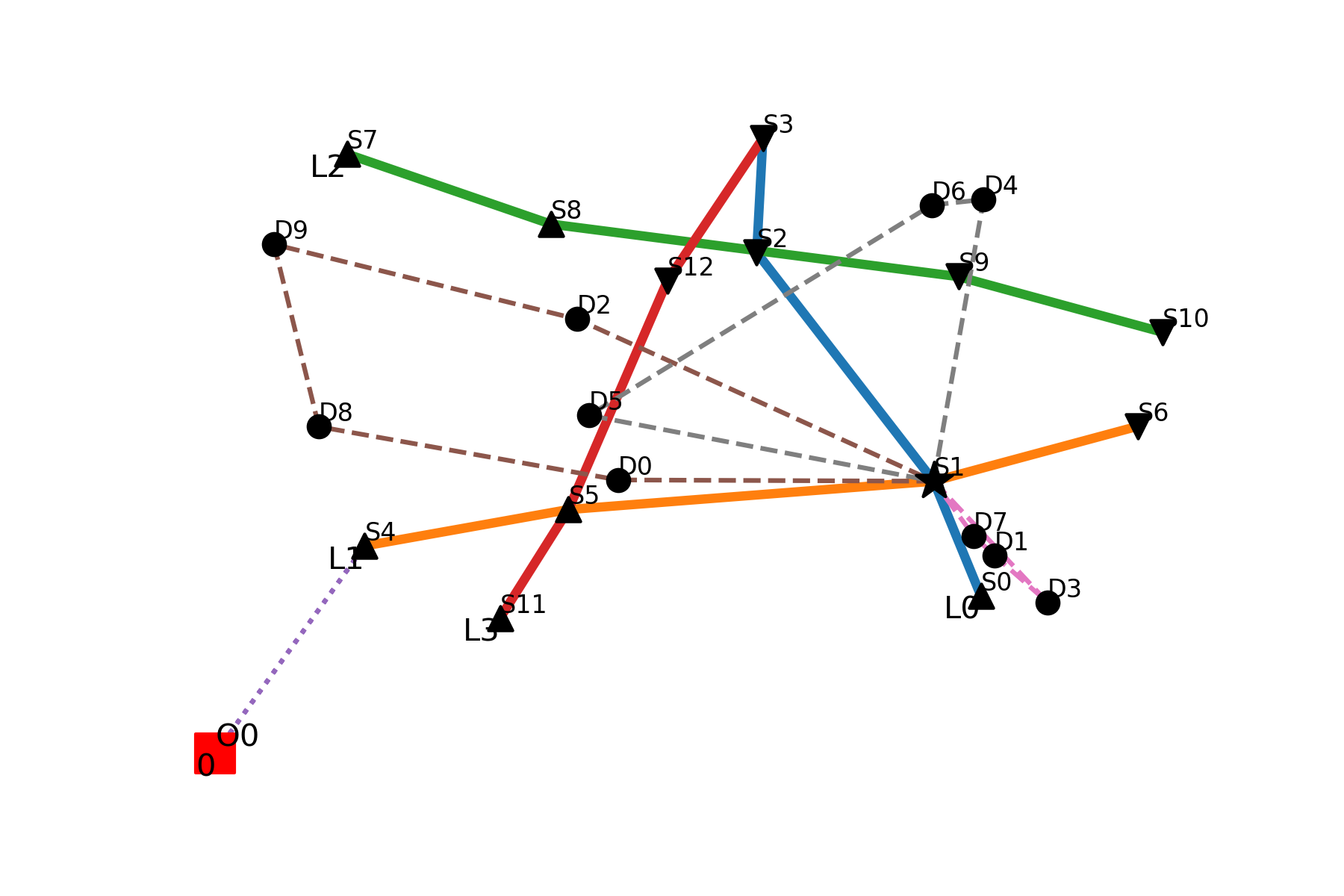}
      \caption{\textit{$\mu = 0.5$, $Q_{F}=50$}}
      \label{fig:SC_5}
    \end{subfigure}%
    %%%%%%%%%%%%%%%%%%%%%%%%%%%%%%%%%%%%
    \begin{subfigure}{.35\textwidth}
      \centering
      \includegraphics[width=1\linewidth]{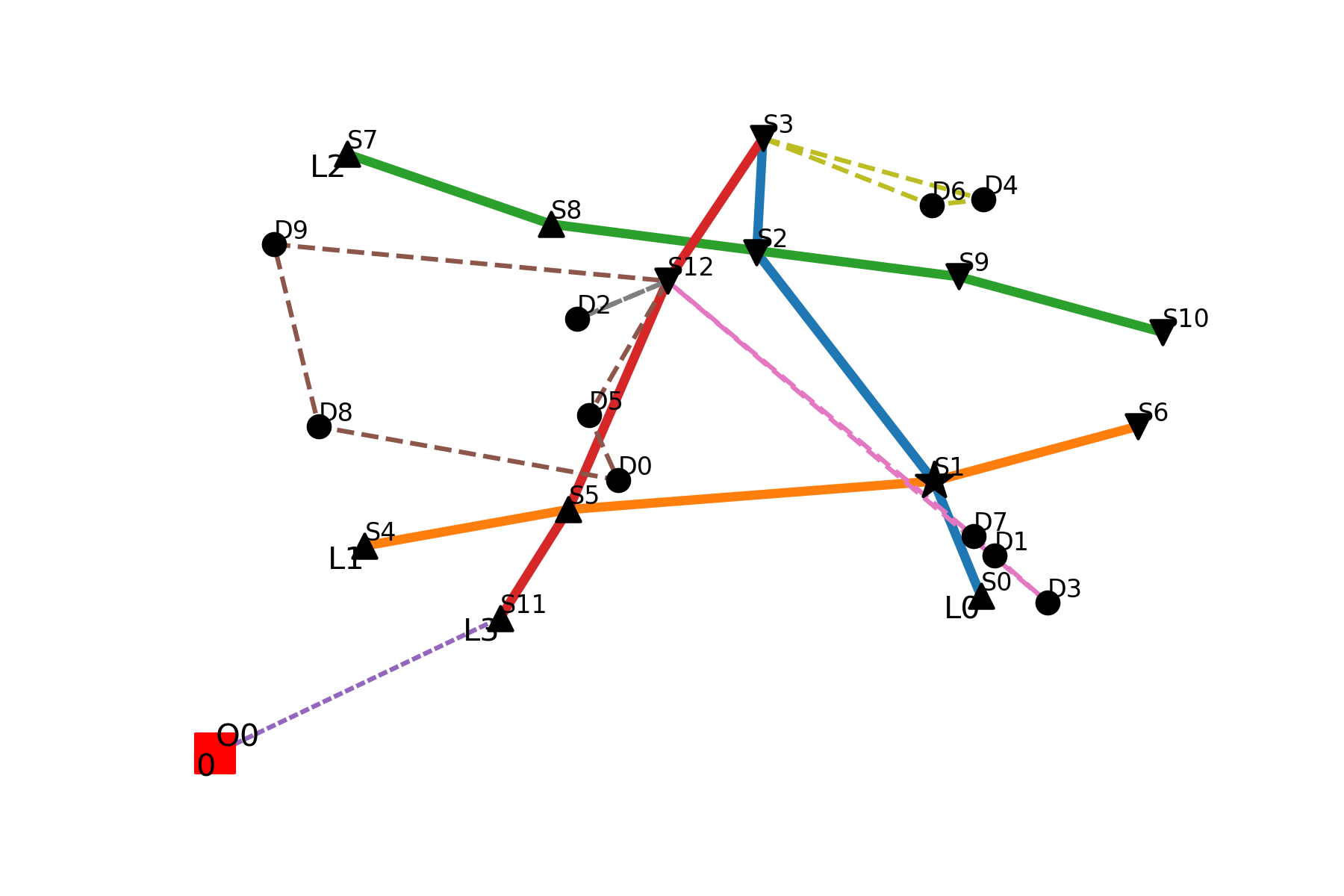}
      \caption{\textit{$\mu = 2.0$, $Q_{F}=50$}}
      \label{fig:SC_6}
    \end{subfigure}
    %%%%%%%%%%%%%%%%%%%%%%%%%%%%%%%%%%%%
    \vspace{-0.5em}
    \caption{\textit{Impact of adding service costs associated with the public transit network}}
    \label{fig:comparison_ServiceCosts}
\end{figure}

However, when the capacity of the freighters is increased, we did indeed observe a greater consolidation at the drop-in and drop-out stops, as shown in Table \ref{tab:service_cost_comparison} and Figure \ref{fig:comparison_ServiceCosts}. We increased the capacity of the freighters, $Q_{F}$, from 20 to 30 and 50 for these tests. Figure \ref{fig:comparison_ServiceCosts} demonstrates the change in the routes in tier 1 and tier 3 as we change the service cost parameter $\mu$ and the freighter capacity $Q_{F}$. \ref{fig:SC_1} shows the solution with the original model parameters, i.e., zero service costs and freighter capacity of 20 units. As the service costs are increased, the delivery truck tends to visit only one drop-in stop. The tier 3 routes also tend to originate from stops on a single line, even though the routing distances increase. In table \ref{tab:service_cost_comparison}, we report the number of drop-in stops, the number of drop-out stops, and the number of freighters used in the delivery for different values of $\mu$ and $Q_{F}$ on the instances with ten customers. We observe that fewer drop-in and drop-out stops are used when increasing service costs. This increase would be more pronounced in instances of larger size. The total cost, or the objective value described in \eqref{obj:full_model_service_costs}, increases proportionally to the service cost parameter based on the number of trucks and freighters used, which, in turn, depends on the number of packages that need to be delivered. If the service costs associated with some lines or stops are relatively high, for example, for stops with extremely high traffic, it is advisable to avoid the stops, even though the routing costs might be increase.

%%%%%%%%%%%%%%%%%%%%%%%%%%%%%%%%%%%%%%%%%%%%%%%%%%%%%%%%%%%%%%%
%%%%%%%%%%%%%%%%%%%%%%%%%%%%%%%%%%%%%%%%%%%%%%%%%%%%%%%%%%%%%%%

% Table generated by Excel2LaTeX from sheet 'Info Compare'
\begin{table}[htbp]
  \centering
  \caption{Impact of adding service costs on instances with 10 customers, solved with $FULL$}
    \resizebox{11cm}{!}{  
    \begin{tabular}{ccccccc}
    \toprule
    \multicolumn{1}{p{3em}}{} & \multicolumn{2}{p{9em}}{\textbf{Number of drop-in stops used}} & \multicolumn{2}{p{9em}}{\textbf{Number of drop-out stops used}} & \multicolumn{2}{p{9em}}{\textbf{Number of freighters used}} \\
    \cmidrule(lr){2-3} \cmidrule(lr){4-5} \cmidrule(lr){6-7}
          $\mu$ & $Q_{F}=30$ & $Q_{F}=50$ & $Q_{F}=30$ & $Q_{F}=50$ & $Q_{F}=30$ & $Q_{F}=50$ \\
    % \midrule
    \cmidrule{1-1} \cmidrule(lr){2-2} \cmidrule(lr){3-3} \cmidrule(lr){4-4} \cmidrule(lr){5-5} \cmidrule(lr){6-6} \cmidrule(lr){7-7}
    0.00     & 1.67  & 1.33  & 2.00  & 1.67  & 5.33  & 3.33 \\
    0.10   & 1.67  & 1.33  & 2.00  & 1.67  & 5.33  & 3.33 \\
    0.25  & 1.33  & 1.00  & 2.00  & 1.67  & 5.33  & 3.33 \\
    0.50   & 1.00  & 1.00  & 1.67  & 1.00  & 5.33  & 3.00 \\
    2.00     & 1.00  & 1.00  & 1.67  & 1.33  & 5.33  & 3.33 \\
    \bottomrule
    \end{tabular}%
    }
  \label{tab:service_cost_comparison}%
\end{table}%

\section{Conclusions} \label{Conclusion}

In this paper, we advocate using public transportation systems for package delivery in cities to reduce emissions and traffic-related issues caused by large delivery vehicles, and show the feasibility and advantages of utilizing such a delivery network. We propose a formulation for the 3T-DPPT, and provide solution methodologies by decomposing the problem into its three tiers and solving them individually. We have three decomposition approaches, along with three objective functions for T2, based on the sequence of solving the tiers. Among all the decomposition approaches, we find the performance of \textit{D2-Obj2} to be the best overall. The decomposition technique \textit{D3} performs well for smaller instances but is quickly impeded as the size of the instances increases. We find \textit{Obj2} to be the best-performing objective function among all. Extensive computational studies support the effectiveness of such an integrated system, which is more sustainable and provides economic opportunities for freight shipping companies as well as public transportation agencies.

Though our solution method is capable of handling small-sized instances, they can be used heuristically to design delivery plans. For example, cities are inherently divided into districts, and our approaches can be implemented in each district individually. The longest part of the delivery is intended to be performed on the public transit network, and the delivery routes in each district can be designed separately. Moreover, we might not have complete information about customers before the delivery plan is formulated, and thus the decisions can be subject to uncertainties. However, the deterministic model studied here can still be beneficial. If we only have an estimate of the customer demands each day instead of the actual demands, from previous implementations of our model, we would know the path that a package at a certain destination needs to follow-- in particular, the drop-in stop, the public vehicle line, and the drop-out stop. Thus, as packages arrive during the day, the delivery company already has an idea of what route a package must follow depending on the hour and the customer’s location. Thus, our study can serve as a viable tool in guiding operational decisions even when uncertainty is involved.

Although we achieved some promising results, a lot remains to be done on the 3T-DPPT. It would be interesting to study the change in the structure of the solutions as the setting of the problem evolves. For example, extensions to our problem could include multiple CDCs, transshipment of packages within the public transit networks, considering uncertainties in the second tier that arise from travel time or capacities, incorporating storage facilities at the stops, etc. Moreover, uncertainties in the problem due to unknown or dynamically arriving customer orders would also be worth investigating. Finally, and most importantly, there is a need to develop solution methodologies to handle larger-sized instances.

One of the limitations of our study is that we assume the system to have been set up before we implement our models. Though there are obvious benefits of the proposed delivery model, convincing authorities and garnering the required financing remain difficult. Amidst several successful examples of integrated delivery systems, some have also failed. The CityCargo project in Amsterdam, which used dedicated trams on existing tracks in the city, is one such example that was shut down, despite initial trial successes and proven potential for cost reductions \citep{marinov2013urban, loy2021combined}. The main reasons identified for its failure were lack of funds, disagreements, and conflicting objectives between the parties involved \citep{loy2021combined, cleophas2019collaborative}. \citet{cochrane2017moving} conducted a study using 34 transportation experts to identify the challenges and opportunities of using public transportation systems for freight movement and utilize expert opinions to find strategies for its implementation. While they acknowledge the benefits of moving FOT systems, they also find that organizational barriers and stakeholder disagreements would pose a more detrimental threat to this implementation than any technical issues. In order to minimize resistance, they suggest the involvement of all stakeholders in the decision-making process from the beginning, starting the operations slowly, and ensuring that the existing services are not disrupted. Thus, for the success of the system, it is imperative for the authorities to take notice and analyze, understand, co-operate, and co-ordinate with the freight delivery companies and public transportation agencies; so that initial investments and installation costs are minimized, all involved parties are satisfied, the delivery service reaches its full potential, and the primary objective of reducing emissions and congestion on roads are fulfilled.

\section*{Acknowledgement}
This work is part of the SISCO project, funded by CY Initiative of Excellence (grant “Investissements d’Avenir” ANR-16-IDEX-0008).

\bibliographystyle{apalike-ejor}
\bibliography{ref}

\newpage

\appendix
\section{Appendix} \label{Appendix}

\subsection{Instance Generation} \label{instance_generation_appendix}

We describe in detail here the generation of instances. As stated in subsection \ref{instance_generation} before, our public transportation networks are similar to the ones proposed in the paper by \citet{donne2021freight}. However, we have slightly modified the networks as discussed below.

\begin{figure}[ht!]
\centering
\begin{subfigure}{.5\textwidth}
  \centering
  \includegraphics[width=1\linewidth]{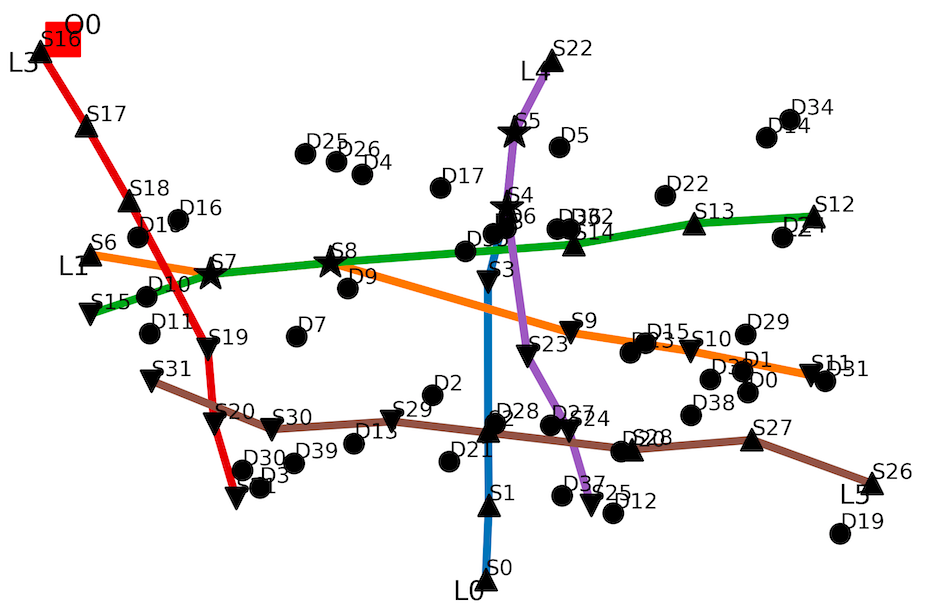}
  \caption{\textit{Original generated instance of the public vehicle network}}
  \label{fig:instance_manipulation_before}
\end{subfigure}%
\begin{subfigure}{.5\textwidth}
  \centering
  \includegraphics[width=1\linewidth]{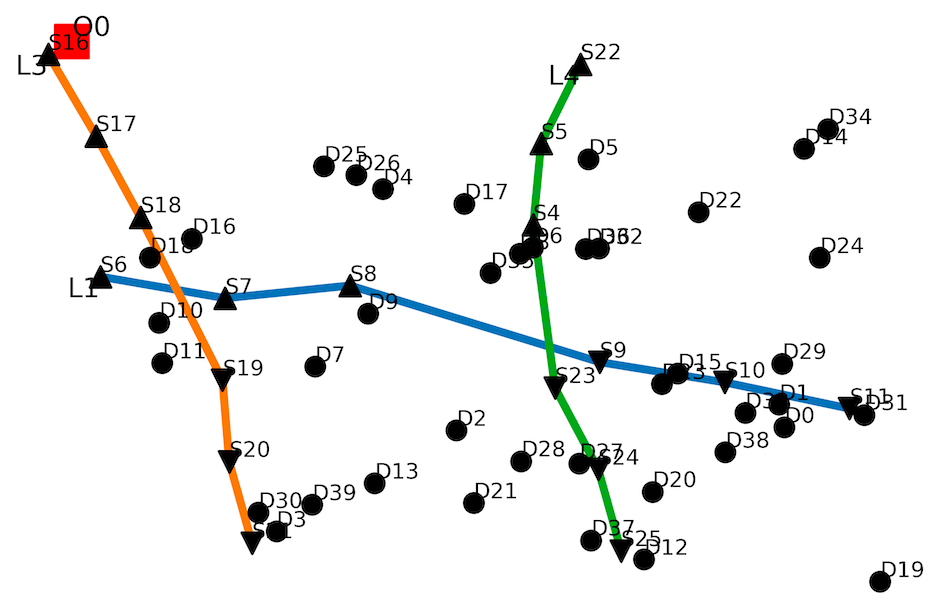}
  \caption{\textit{Modified network used in the paper}}
  \label{fig:instance manipulation after}
\end{subfigure}
\vspace{-1em}
\caption{\textit{Plot showing an instance before and after manipulation according to the criterion- ``the farthest drop-in is closer to the CDC than the nearest drop-out"}}
\label{fig:instance_manipulation}
\end{figure}

For each randomly generated instance, with a certain number of public transportation lines, we filter the lines so that it is profitable to send at least one package to a customer using the line. To do so, we consider the farthest drop-in stop and the nearest drop-out stop from the CDC. If the farthest drop-in stop is closer to the origin than the nearest drop-out stop, then we keep the line because it could still be profitable to use the line. In other words, we select the lines with at least one drop-in station closer to the CDC than its drop-out stops. \ref{fig:instance_manipulation_before} and \ref{fig:instance manipulation after} show an example of an instance before and after filtering out the lines using the criterion mentioned above. Figure \ref{fig:instance_manipulation_before} show a city public vehicle network with 6 lines. The assigned drop-in stops for lines L0, L2, and L5 are farther away from the CDC than the drop-out stops on these lines. Hence, we remove them from the instance to modify it, finally only keeping lines L1, L3, and L4, as shown in figure \ref{fig:instance manipulation after}.

%%%%%%%%%%%%%%%%%%%%%%%%%%%%%%%%%%%%%%%%%%%%%
%%%%%%%%%%%%%%%%%%%%%%%%%%%%%%%%%%%%%%%%%%%%%

Every customer has a set of drop-out stops that can serve them, and each drop-out stop can potentially deliver to at least one customer lying within a given radius. The radius is set as three times the average distance between the stops. If a customer has no drop-out stop within the given radius, we assign the three closest drop-out stops from that customer as their set of drop-out stops. Finally, sometimes a drop-out stop on the network is not assigned to any customer. In that case, we could either remove the drop-out stop from the instance altogether or generate a customer that the stop could potentially serve. We proceed with the latter and randomly generate a customer lying close to the stop. The set of drop-out stops that can serve this new customer is given by that stop, plus two other closest stops to the customer. We do this to ensure that each customer can be delivered from at least one drop-out stop; on the other hand, each drop-out stop can potentially serve at least one customer. Thus, it makes sense to include the entire network in the problem.

Next, we discuss the time windows of the customers. We generate the time windows such that there is enough time for the trucks and the public vehicles to make their deliveries. As discussed earlier, we have 30 periods, each 30 minutes long. We randomly generate the starting time window of a customer between the $3^{rd}$ and the $15^{th}$ period. We make sure to have at least three hours between the beginning and the end of each window. In other words, each time window is at least 3 hours long. Thus, the end of a delivery window is randomly generated between three hours after its beginning and the last period. This ensures we have enough time and resources to deliver each customer's package. It is also realistic because, in practice, some customers have narrower time windows, while others have wider time windows.

We set the number of delivery trucks for each instance such that there are enough trucks to serve all the customers. For the number of freighters, we proceed as follows. For each drop-out stop, we consider the number of customers that the stop could potentially serve. Then, we set the number of freighters at each stop to be the maximum of the above.

We summarize the description of the instances in Table \ref{tab:instance_description}.

\begin{table}[ht]
  \centering
  \caption{\textit{Table for description of instances}}
  \vspace{-1em}
  \resizebox{14cm}{!}{
    \begin{tabular}{cccccccccc}
    \toprule
    \multicolumn{1}{p{5em}}{\centering \textbf{Instance number}} & \multicolumn{1}{p{5.3em}}{\centering \textbf{Number of customers}} & \multicolumn{1}{p{4em}}{\centering \textbf{Total number of stops}} & \multicolumn{1}{p{4.5em}}{\centering \textbf{Number of drop-in stops}} & \multicolumn{1}{p{4.5em}}{\centering \textbf{Number of drop-out stops}} & \multicolumn{1}{p{3.5em}}{\centering \textbf{Number of lines}} & \multicolumn{1}{p{7.9em}}{\centering \textbf{Number of public vehicles per line (frequency)}} & \multicolumn{1}{p{6em}}{\centering \textbf{Total number of public vehicles}} & \multicolumn{1}{p{4em}}{\centering \textbf{Number of delivery trucks}} & \multicolumn{1}{p{4.2em}}{\centering \textbf{Number of freighters per stop}} \\
    \midrule
    1     & 10    & 4     & 2     & 2     & 1     & 15    & 15    & 5     & 10 \\
    2     & 10    & 8     & 4     & 4     & 2     & 15    & 30    & 5     & 7 \\
    3     & 10    & 13    & 7     & 7     & 4     & 15    & 60    & 5     & 7 \\
    4     & 20    & 22    & 10    & 12    & 4     & 15    & 60    & 5     & 10 \\
    5     & 20    & 23    & 11    & 12    & 4     & 15    & 60    & 5     & 9 \\
    6     & 20    & 17    & 8     & 9     & 3     & 15    & 45    & 5     & 11 \\
    7     & 30    & 23    & 12    & 12    & 4     & 15    & 60    & 8     & 16 \\
    8     & 30    & 21    & 11    & 11    & 4     & 15    & 60    & 8     & 18 \\
    9     & 30    & 22    & 12    & 10    & 4     & 15    & 60    & 8     & 15 \\
    10    & 40    & 18    & 9     & 9     & 3     & 15    & 45    & 8     & 25 \\
    11    & 40    & 27    & 14    & 15    & 5     & 15    & 75    & 8     & 21 \\
    12    & 40    & 27    & 15    & 12    & 5     & 15    & 75    & 8     & 20 \\
    13    & 50    & 32    & 16    & 16    & 7     & 15    & 105   & 10    & 23 \\
    14    & 50    & 22    & 10    & 12    & 4     & 15    & 60    & 10    & 25 \\
    15    & 50    & 30    & 17    & 14    & 6     & 15    & 90    & 10    & 27 \\
    16    & 60    & 33    & 17    & 17    & 7     & 18    & 126   & 14    & 44 \\
    17    & 60    & 28    & 15    & 14    & 5     & 18    & 90    & 14    & 26 \\
    18    & 60    & 23    & 14    & 10    & 5     & 18    & 90    & 14    & 23 \\
    19    & 70    & 28    & 16    & 12    & 4     & 18    & 72    & 14    & 51 \\
    20    & 70    & 31    & 15    & 16    & 4     & 18    & 72    & 14    & 37 \\
    21    & 70    & 35    & 18    & 17    & 5     & 18    & 90    & 14    & 23 \\
    22    & 80    & 40    & 21    & 20    & 6     & 18    & 108   & 14    & 29 \\
    23    & 80    & 43    & 24    & 20    & 6     & 18    & 108   & 14    & 30 \\
    24    & 80    & 44    & 24    & 20    & 6     & 18    & 108   & 14    & 40 \\
    \bottomrule
    \end{tabular}%
  }%
  \label{tab:instance_description}%
\end{table}%

\FloatBarrier

\subsection{Formulation Used for Simple VRP with Time Windows (VRPTW)} \label{Appendix_VRPTW}

We use very similar notations and variables to formulate a simple VRP with time windows, with only minor changes. The set of customers is denoted by $\mathcal{C}$. $o$ denotes the CDC, and $o'$ denotes a copy of the CDC. Let $\mathcal{N} = \mathcal{C} \cup \{ o \} \cup \{ o' \}$. Let $\mathcal{D}$ denote the set of delivery trucks. Let $r_{id}$ be a binary variable that takes the value 1 if the package for customer $i$ is assigned to delivery truck $d$ and 0 otherwise. Let $w_{uvd}$ be a binary variable if a truck $d$ traverses the arc $(u,v)$, where $u \in \mathcal{C} \cup \{ o \}, \; v \in \mathcal{C} \cup \{ o' \}$, and 0 otherwise, and let $T_{uvd}$ be the time taken to do so. Finally, we have the set of continuous variables $t_{ud}$ which denotes the time when the truck $d$ leaves the node $\mathcal{C} \cup \{ o \}$. Then, the VRPTW is given by:

\begin{allowdisplaybreaks}
\begin{align}
    \textit{Minimize} \quad & \sum_{d \in \mathcal{D}} \; \sum_{u \in \mathcal{C} \cup \{o \} } \sum\limits_{ \substack{ v \in \mathcal{C} \cup \{o' \} \\ u \neq v } } C^{1}_{uvd} w_{uvd} \label{VRPTW:obj} \\
    \textit{sub to:} \quad & \sum_{d \in \mathcal{D}} r_{id} = 1, \; \forall i \in \mathcal{C} \label{VRPTW:unique_customer_truck} \\
    & \sum_{i \in C} q_{i} r_{id} \leq Q^{1}_{d}, \; \forall d \in \mathcal{D} \label{VRPTW:truck_capacity} \\
    & \sum_{v \in \mathcal{C} \cup \{o'\}} w_{ovd} = 1, \; \forall d \in \mathcal{D} \label{VRPTW:truck_starting_node} \\
    & \sum_{v \in \mathcal{C} \cup \{o\}} w_{vo'd} = 1, \; \forall d \in \mathcal{D} \label{VRPTW:truck_ending_node} \\
    & \sum_{\mathclap{\substack{v \in \mathcal{C} \cup \{o'\} \\ v \neq u}}} \; w_{uvd} = \sum_{\mathclap{\substack{v \in \mathcal{C} \cup \{o\} \\ v \neq u}}} \; w_{vud}, \; \forall u \in \mathcal{C}, \; d \in \mathcal{D} \label{VRPTW:truck_node_balance} \\
    & \sum_{\mathclap{\substack{u \in \mathcal{C} \cup \{o'\} \\ v \neq u}}} \; w_{uvd} = r_{vd}, \; \forall v \in \mathcal{C}, \; d \in \mathcal{D} \label{VRPTW:customer_truck_route_link} \\
    & t_{vd} \geq t_{ud} + T_{uvd} + \widehat{T}_{v} - M (1 - w_{uvd}), \; \forall u \in \mathcal{C} \cup \{o\}, \; v \in \mathcal{C}, \; u \neq v, \; d \in \mathcal{D} \label{VRPTW:truck_time_balance} \\
    & t_{id} \geq \underline{T}_{i} - M (1 - r_{id}), \; \forall i \in \mathcal{C}, \; d \in \mathcal{D} \label{VRPTW:customer_time_limits_lower} \\
    & t_{id} \leq \overline{T}_{i} + M (1 - r_{id}), \; \forall i \in \mathcal{C}, \; d \in \mathcal{D} \label{VRPTW:customer_time_limits_upper} \\
    & \sum_{v \in \mathcal{C}} w_{ovd} \geq \sum_{v \in \mathcal{C}} w_{ovd+1}, \; \forall d \in |\mathcal{D}|-1 \label{VRPTW:symmetry_breaking_tier1_1} \\
    & \sum_{u \in \mathcal{C} \cup \{ o \} } \sum_{v \in \mathcal{C} \cup \{ o' \} } w_{uvd} \geq \sum_{u \in \mathcal{C} \cup \{ o \} } \sum_{v \in \mathcal{C} \cup \{ o' \} } w_{uvd+1} , \; \forall d \in |\mathcal{D}|-1 \label{VRPTW:symmetry_breaking_tier1_2} \\
    & r_{id} \in \{0,1\}, \; \forall i \in \mathcal{C}, \forall d \in \mathcal{D} \\
    & w_{uvd} \in \{0,1\}, \; \forall u \in \mathcal{C} \cup \{ o \}, \; v \in \mathcal{C} \cup \{ o' \}, \; d \in \mathcal{D} \\
    & t_{id} \geq 0, \; \forall i \in \mathcal{C} \cup \{ o \}, \; d \in \mathcal{D} 
\end{align}
\end{allowdisplaybreaks}

\end{document}